\documentclass[11pt]{article}
\usepackage{graphicx}
\usepackage{epstopdf}
\usepackage{amsmath, amssymb}
\usepackage{latexsym, color}
\usepackage{mathtools}
\usepackage{xcolor}
\usepackage{graphicx}

\def\la{\big\langle}
\def\ra{\big\rangle}

\def\ds{\displaystyle}
\def\forall{\hbox{for all}~}
\def\L{{\bf L}}

\def\bft{{\bf t}}
\def\bfw{{\bf w}}

\def\bfn{{\bf n}}

\def\ve{\varepsilon}

\def\Ups{\Upsilon}
\def\dint{\int\!\!\int}

\def\R{{\mathbb R} }

\def\implies{\Longrightarrow}
\def\vp{\varphi}

\def\P{{\cal P}}
\def\F{{\cal F}}

\def\vs{\vskip 2em}

\def\v{\vskip 1em}
\def\E{{\cal E}}
\def\A{{\cal A}}

\def\H{{\cal H}}
\def\T{{\cal T}}
\def\D{{\cal D}}
\def\O{{\cal O}}

\def\C{{\cal C}}
\def\caL{{\cal L}}
\def\J{{\cal J}}

\def\H{{\cal H}}
\def\rest{\llcorner}
\def\bega{\begin{array}}
\def\enda{\end{array}}
\def\begi{\begin{itemize}}
\def\endi{\end{itemize}}
\def\ov{\overline}

\def\Tilde{\widetilde}
\def\Hat{\widehat}

\def\bel{\begin{equation}\label}
\def\eeq{\end{equation}}
\def\sqr#1#2{\vbox{\hrule height .#2pt
\hbox{\vrule width .#2pt height #1pt \kern #1pt
\vrule width .#2pt}\hrule height .#2pt }}
\def\square{\sqr74}
\def\endproof{\hphantom{MM}\hfill\llap{$\square$}\goodbreak}
\definecolor{cadmiumgreen}{rgb}{0.0, 0.42, 0.24}

\newtheorem{theorem}{Theorem}[section]
\newtheorem{corollary}{Corollary}[section]

\newtheorem{proposition}{Proposition}[section]
\newtheorem{remark}{Remark}[section]
\newtheorem{definition}{Definition}[section]
\newtheorem{example}{Example}[section]

\begin{document}
\title{\bf Optimally Controlled Moving Sets\\
with Geographical Constraints}
\vs
\author{Alberto Bressan$^{(1)}$, Elsa M.~Marchini$^{(2)}$ and Vasile Staicu$^{(3)}$\\
\, \\
{\small $^{(1)}$Department of Mathematics, Penn State University,} \\{\small University Park, PA~16802, USA.}\\
{\small $^{(2)}$Dipartimento di Matematica, Politecnico di Milano,}\\{\small Piazza L.\,da Vinci, 32 - 20133 Milano, Italy.}\\
{\small $^{(3)}$ Center for Research and Development in Mathematics and Applications (CIDMA),}\\
{\small Department of Mathematics, University of Aveiro, 3810-193 Aveiro, Portugal.} \\
\, \\
{\small E-mails: axb62@psu.edu,~elsa.marchini@polimi.it,~vasile@ua.pt.}
}
\maketitle

\begin{abstract}  The paper is concerned with a family of geometric evolution problems, modeling the spatial control 
of an invasive population within a region $V\subset \R^2$ bounded by geographical barriers.   
If no control is applied, the contaminated set $\Omega(t)\subset V$ expands with unit speed in all directions.
By implementing a control, a region of area $M$ can be cleared up per unit time.

Given an initial set $\Omega(0)=\Omega_0\subseteq V$, 
three main problems are studied: (1)  Existence
of an admissible strategy $t\mapsto\Omega(t)$ which eradicates the contamination in finite time, so that 
$\Omega(T)=\emptyset$ for some $T>0$.
(2) Optimal strategies that achieve eradication in minimum time.
(3)  Strategies that minimize the average area of the contaminated set on a given time interval $[0,T]$.
For these optimization problems, a sufficient condition for optimality is proved, together with 
several necessary conditions.  Based on these conditions, optimal  set-valued motions $t\mapsto \Omega(t)$
are explicitly  constructed in a number of cases.
\end{abstract}

\section{Introduction}
\label{s:1}
\setcounter{equation}{0}
We consider a family of geometric evolution problems, modeling the spatial control 
of an invasive population 
\cite{BCS1}.
For each time $t\in [0,T]$, we denote by $\Omega(t)\subset\R^2$ a set moving in the plane.
This can be regarded as a ``contaminated region", which we would like to shrink as 
much as possible.   To control the evolution of this set, we assign the 
velocity $\beta=\beta(t,x)$ in the inward normal direction at every boundary point $x\in \partial \Omega(t)$.

\begin{figure}[ht]
\centerline{\hbox{\includegraphics[width=6cm]{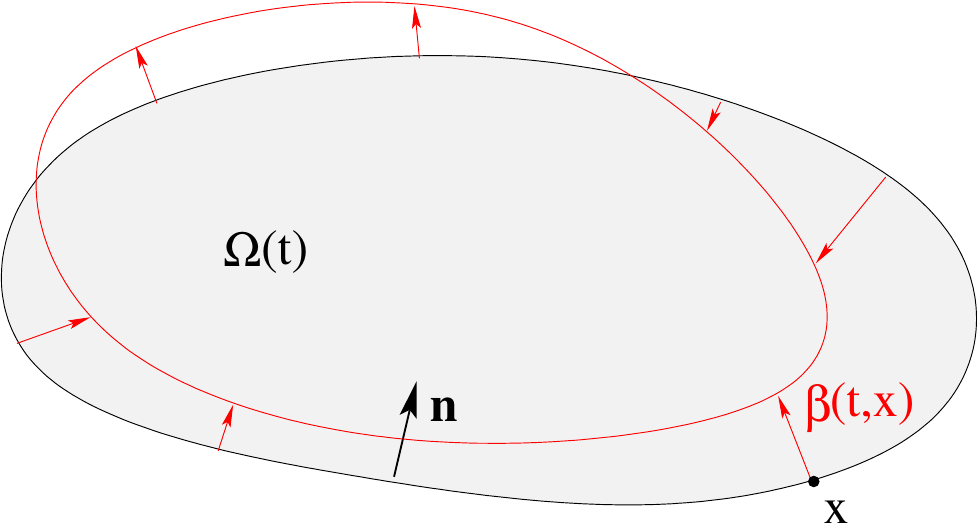}}
\qquad\qquad \hbox{\includegraphics[width=5cm]{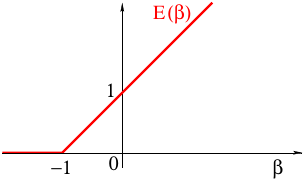}}}
\caption{\small Left: a moving set, where the evolution is determined by assigning the inward normal speed $\beta$
at each boundary point.
Right: the effort function $E(\beta)$ in (\ref{E}).
}
\label{f:sm345}
\end{figure}

A function $E(\beta)\geq 0$ is given, 
describing the {\bf effort} needed to push the boundary of $\Omega(t)$ inward, with speed $\beta$ in the 
normal direction (see Fig.~\ref{f:sm345}). 
 The {\bf total control effort} at time $t\in [0,T]$
is then defined as
\bel{Et} \E(t)~\doteq~
\int_{\partial\Omega(t)}E\bigl(\beta(t,x)\bigr)\, \H^1(dx),
\eeq
where the integral is computed w.r.t.~the 1-dimensional Hausdorff measure
along the boundary of $\Omega(t)$. 
In this paper we focus on the case where
\bel{E} E(\beta)~=~\left\{\bega{cl} 1+\beta\qquad&\hbox{if}~~\beta\geq -1,\\[1mm]
0 \qquad&\hbox{if}~~\beta< -1.\enda\right.\eeq
Given a constant $M>0$ accounting for the maximum control effort, we consider set motions
$t\mapsto \Omega(t)$
which satisfy the constraint
\bel{EM} \E(t)~\leq~M\qquad\forall ~t\in [0,T].\eeq

This models a situation where:
\begi
\item If the control effort is everywhere zero: $E(\beta)=0$, then the inward normal speed is $\beta= -1$
at every point. 
Hence the contaminated set $\Omega(t)$ expands with unit speed in all directions.
In particular, its area increases at a rate equal to the perimeter:
$${d\over dt} \caL^2\bigl(\Omega(t)\bigr)~=~- \int_{\partial \Omega(t)} \beta(t,x) \, \H^1(dx)~=~ 
\H^1\bigl(\partial \Omega(t)\bigr).$$
\item By implementing a control with total effort  $\E(t)= M$, we can clean up a region of area $M$ per unit time,
hence
\bel{darea2} \bega{rl} \ds{d\over dt} \caL^2\bigl(\Omega(t)\bigr)&\ds=~-\int_{\partial \Omega(t)}  \beta(t,x) \, \H^1(dx)~=~
\int_{\partial \Omega(t)} \Big[1 - E\bigl(\beta(t,x)\bigr)\Big] \, \H^1(dx)\\[4mm]
&=~\H^1\bigl(\partial \Omega(t)\bigr)-M.\enda\eeq
\endi
Here and in the sequel, $\H^1$ denotes the 1-dimensional Hausdorff measure while  $ \caL^2$ is the 2-dimensional 
Lebesgue measure. 
Given a set $S\subset\R^2$, we write
$$B_\delta(S)~\doteq~\bigl\{ x\in\R^2\,;~~d(x, S)\leq \delta\bigr\}$$
for the $\delta$-neighborhood around the set $S$.
Motivated by (\ref{darea2}), following \cite{BBC} we  introduce
\begin{definition}\label{def:11} Given a constant $M>0$, we say that the set motion $t\mapsto\Omega(t)$ is {\bf admissible} if 
\bel{a1} \Omega(t+h) ~\subseteq~
B_h\bigl(\Omega(t)\bigr)\qquad\qquad\qquad \forall ~0\leq t<t+h\leq T,\eeq
\bel{a2} \limsup_{h\to 0+} {\caL^2\Big( B_h\bigl(\Omega(t)\bigr)\setminus \Omega(t+h) \Big) \over h}
~\leq~M\qquad\forall~ 0\leq t<T.\eeq
\end{definition}

In the above setting, two basic problems can be formulated.

\begi
\item[{\bf (EP)}] {\bf Eradication Problem.}  
{\it   Let an initial set $\Omega_0\subset\R^2$ and a constant $M>0$
be given. Find an admissible set-valued function $t\mapsto \Omega(t)$  such that, for some $T>0$,
\bel{nco}\Omega(0)~=~\Omega_0,\qquad\qquad 
\Omega(T)~=~\emptyset,\eeq
}
\endi
\v
\begi
\item[{\bf (MTP)}]  {\bf Minimum Time Problem.}  {\it  Among all admissible strategies that satisfy
(\ref{nco}), find one which minimizes the time $T$.}
\endi

Notice that the constant $M$ puts an upper bound on the total control effort at each time $t$.
For example, in a pest eradication model, there will be an upper bound
on the amount of pesticides that can be sprayed per unit time.  
If $M$ is too small, compared with the size of the contaminated region, 
it may not be possible to completely eradicate the invasive population.  

More generally, even if the Eradication Problem does not have a solution,
one can consider an optimization problem, minimizing the size of the contaminated set over time.
\v
\begi
\item[{\bf (OP)}]  {\bf Optimization Problem.}  {\it  Given an initial set $\Omega_0$ and two constants $\kappa_1,\kappa_2>0$, find an admissible motion 
$t\mapsto \Omega(t)$ 
that minimizes the cost functional 
\bel{J}J(\Omega)~\doteq~\kappa_1\int_0^T \caL^2\bigl(\Omega(t)\bigr)\, dt + \kappa_2 \,\caL^2\bigl(\Omega(T)\bigr).
\eeq
with the initial data
$ \Omega(0)~=~\Omega_0$.
}
\endi

The main results in \cite{BCS2} provide the existence of optimal solutions,
together with necessary conditions for optimality.  
In the case where the initial set $\Omega_0\subset\R^2$ is convex, the optimal solution
to  {\bf (MTP)} and {\bf (OP)}  has been explicitly determined in 
\cite{BBC}.

In the present paper we study similar control problems for moving sets, but with 
geographical constraints.  These models describe an invasive biological population within an island,
where the sea provides a natural barrier to its expansion.
More precisely,
given an open set $V\subset\R^2$
we impose the additional constraint
$\Omega(t)\subseteq V$.
Since the population cannot propagate outside $V$,
 the instantaneous control effort (\ref{Et}) is now replaced by
\bel{CEt} \E(t)~\doteq~
\int_{\partial\Omega(t)\cap V}E\bigl(\beta(t,x)\bigr)\, \H^1(dx).
\eeq
Notice that in (\ref{CEt}) the effort is integrated only over the relative boundary of $\Omega(t)$,
contained inside the open set $V$. In this setting, Definition~\ref{def:11} is replaced by
\begin{definition}\label{def:12} Given a constant $M>0$ and a bounded open domain 
$V\subset\R^2$, we say that the set motion $t\mapsto\Omega(t)\subseteq V$ is {\bf admissible} if 
\bel{a1c} \Omega(t+h) ~\subseteq~
V\cap B_h\bigl(\Omega(t)\bigr)\qquad\qquad\qquad \forall ~0\leq t<t+h\leq T,\eeq
\bel{a2c} \limsup_{h\to 0+}~ { \caL^2\bigg( \Big( V\cap B_h\bigl(\Omega(t)\bigr)\Big) \setminus \Omega(t+h)  \bigg)\over h} 
~\leq~M\qquad\forall~ 0\leq t<T.\eeq
\end{definition}
The three problems {\bf (EP)}, {\bf (MTP)} and {\bf (OP)}  can now be formulated 
in the same way as before, replacing (\ref{Et}) with (\ref{CEt}).

Aim of this paper is to provide an extensive analysis of the above problems, in the presence of geographical constraints.
In the first part, several results are proved on the existence of solutions, sufficient conditions for optimality, 
and various necessary conditions. In the second part, we use the necessary conditions in order 
to explicitly construct optimal solutions in a variety of  cases.

More in detail, 
in Section~\ref{s:2} we give an equivalent notion of admissible strategy and
review the definition of weak solutions to the optimization problem,
introduced in~\cite{BCS2}.   Section~\ref{s:3} deals with the eradication problem {\bf (EP)}. 
Concerning the existence of an admissible eradication strategy, 
sufficient conditions as well as  necessary conditions
are proved in Theorem~\ref{t:31}.

In Section~\ref{s:4} we study the minimum time problem {\bf (MTP)}.
A sufficient condition,  for a strategy to eradicate the contamination in minimum time, is proved in 
Theorem~\ref{t:41}.  Examples are given, showing cases where such condition can be used. Here
the optimal strategy amounts to piecing together a family of solutions to the classical Dido's problem.
Namely, at each time $t\in [0,T]$, the relative boundary of the set $\Omega(t)\subseteq V$ has minimum length, 
compared with all other subsets $W\subseteq V$ of the same area.

In Section~\ref{s:5} we prove a general result on the existence of optimal solutions for the problem
{\bf (OSM)} and for the minimum time problem {\bf (MTP}), within a class of multifunctions with bounded variation.

Necessary conditions for optimality in a form similar to the Pontryagin maximum principle, are stated in Section~\ref{s:6}.
An intuitive argument is first given, 
motivating the result and explaining the role of the adjoint variable, which can be interpreted
as a ``shadow price for a cleaning service".  A detailed proof of the necessary conditions is worked out in
the following Sections~\ref{s:7} and~\ref{s:8}.   
Optimality conditions for the minimum time problem are then derived in Section~\ref{s:9}.

We should point out from the outset that the regularity of the set-valued map $t\mapsto \Omega(t)$ 
is a major issue.   At the present time, the existence of optimal strategies is proved within a class
of maps whose characteristic function
\bel{Omchar} {\bf 1}_{\Omega}(t,x_1, x_2)~=~\left\{\bega{rl} 1\quad &\hbox{if}~ (x_1, x_2)\in \Omega(t),\\[1mm]
0\quad &\hbox{if} ~(x_1, x_2)\notin \Omega(t),\enda\right.\eeq
has bounded variation \cite{AFP, EG, M}.
However, the optimality conditions obtained in \cite{BCS2} require a much stronger regularity assumption. 
Namely,  the boundaries $\partial \Omega(t)$ should be
$\C^2$ curves.
More recently, the explicit solutions constructed in \cite{BBC} show that these conditions are too restrictive.   
To handle the majority of relevant cases, is thus important
to relax these regularity assumptions.
Here we assume that  the boundary $\partial \Omega(t)$ is only $\C^1$ with Lipschitz continuous normal vector.
Its curvature $\omega(t,x)$ is defined  almost everywhere, and
the inward normal velocity $\beta(t,x)$ is continuous.
Because of these weaker assumptions, a more careful proof is needed.

In spite of this improvement, there still remains a major gap between the regularity provided by the 
existence theorem and the regularity required by the optimality conditions.  
To partially bridge this gap, in 
Sections~\ref{s:10} and \ref{s:11} we prove additional necessary conditions, which must be satisfied at  boundary points assuming only piecewise $\C^1$ regularity. Roughly speaking, Theorem~\ref{t:101} shows that
the boundary $\partial \Omega(t)$ cannot have corner points in the interior of $V$. 
Moreover, Theorem~\ref{t:111} 
implies that, at points where
$\partial\Omega(t)$ meets the boundary of the constraining set $V$, either the junction is perpendicular, or the 
control effort must vanish.

In the second part of the paper, Sections \ref{s:12} to \ref{s:16}, we study the problem of
how to construct admissible motions $t\mapsto\Omega(t)\subset V$  which satisfy 
all the previous optimality conditions.  Here the discussion is partly rigorous, partly heuristic.
We make no attempt to cover all possible situations.   Rather, we provide guidelines for constructing 
these set-motions.    Explicit formulas are obtained  in a number of specific cases.
The reader should be aware that, although these motions are the unique ones that 
satisfy all our necessary conditions for optimality,  this does not necessarily imply their optimality.
In principle, there may be other strategies, with only BV regularity, which achieve a lower cost.

We outline the main ideas.  At each time $t\in [0,T]$ the boundary $\partial
\Omega(t)$ is assumed to be the union of  finitely many ``controlled arcs", where the control is active, and ``free arcs",
where no control is applied and the set $\Omega(t)$ thus expands with unit speed. 
If at some intermediate time $\tau\in [0,T]$ a maximally extended free arc is identified, 
in turn this determines the free arcs
at all times $t$ in a neighborhood of $\tau$. The remaining portions of the boundary can now be constructed, 
since they must be arcs 
of circumferences, all with the same radius $r(t)$, whose endpoints either (i) meet tangentially a free arc, or (ii) 
cross perpendicularly the boundary $\partial V$.    Thanks to the identity (\ref{darea2}) we obtain an ODE 
for the radius $r(t)$.  In suitable cases, the boundaries $\partial\Omega(t)$ can thus be completely determined.

This general approach is described in Section~\ref{s:12}, together with some examples.
Details are discussed in the following sections. The initial stage of an optimal 
strategy is worked out in Section~\ref{s:13}, while in Section~\ref{s:14}  we derive a set of 
equations describing a
maximally extended free interface,  assuming that the set $V$ has a smooth boundary.

In Section~\ref{s:15} we assume that $V$ is a polygon, and  take up the more challenging task of constructing 
a maximal free arc through one of its vertices.
In this setting, maximal arcs are determined by a system of ODEs which is of second order,
implicit, and singular.  
In Section~\ref{s:16} a 1-parameter family of local solutions is constructed, in a neighborhood of a vertex.

Finally, Section~\ref{s:17} contains
some concluding remarks and a brief discussion of open problems.

Geometric optimization problem have been the subject of a rich literature.
We refer to \cite{BuBu, DZ, Mord} for comprehensive monographs, with further references.
More specifically, several different models related to the control of a moving set have recently been considered in
\cite{Bblock, Breview, BMN, BZ, CPo, CLP}.    
Eradication problems for invasive biological 
species are studied in \cite{ACD, ACM, ACS}.

\section{An equivalent formulation of admissible set motions}
\label{s:2}
\setcounter{equation}{0}
The concept of admissible set-motion, introduced is Definitions~\ref{def:11} and \ref{def:12}, is very general.
Indeed,  it
remains meaningful for any measurable set-valued map $t\mapsto\Omega(t)$. 
However, toward the analysis of optimal strategies, it is convenient to work within a class of more regular maps.
Following \cite{BCS2}, in this section we consider a weak formulation of the optimization problem
{\bf (OP)}, within the class of BV functions.

Let a bounded open set $V\subset\R^2$ be given, with finite perimeter.
Consider the family of subsets
\bel{Adef} \F~\doteq~\Big\{  \Omega\subset \,]0,T[\,\times V\,;~~\Omega ~\hbox{has
finite perimeter}\Big\}.\eeq  
Calling ${\bf 1}_\Omega$ the characteristic function of $\Omega$, 
this implies that  ${\bf 1}_\Omega\in BV$.
In other words, the distributional gradient $\mu_\Omega\doteq D\, {\bf 1}_\Omega$
is a finite $\R^3$-valued Radon measure:
\bel{div}\int_\Omega \hbox{div}\, \vp\, dx~=~-\int \vp\cdot d\mu_\Omega\qquad\qquad \forall \vp\in\C^1_c
\bigl(]0,T[\,\times\R^2\,;~\R^3\bigr).\eeq
Given a set $\Omega\in \F$, we consider the multifunction
\bel{slice} t~\mapsto~\Omega(t)~\doteq~\bigl\{ x\in\R^2\,;~(t,x)\in \Omega\bigr\}.\eeq
By possibly modifying ${\bf 1}_\Omega$ on a set of 3-dimensional measure zero,
the map $t\mapsto {\bf 1}_{\Omega(t)}$
has bounded variation from $\,]0,T[\,$ into $\L^1(\R^2)$.  In particular, for every $0<t<T$,
the one-sided limits
\bel{ompm}{\bf 1}_{\Omega(t+)}~\doteq~\lim_{t\to t+} {\bf 1}_{\Omega(t)}\,,
\qquad\qquad {\bf 1}_{\Omega(t-)}~\doteq~\lim_{t\to t-} {\bf 1}_{\Omega(t)}\,,\eeq
are well defined in $\L^1(\R^2)$.   This uniquely defines the sets $\Omega(t+)$,
$\Omega(t-)$, up to a set of 2-dimensional Lebesgue measure zero.
Throughout the following, we define
the sets $\Omega(0)$ and $\Omega(T)$ in terms of
\bel{om0T}{\bf 1}_{\Omega(0)}~\doteq~\lim_{t\to 0+} {\bf 1}_{\Omega(t)}\,,
\qquad\qquad {\bf 1}_{\Omega(T)}~\doteq~\lim_{t\to T-} {\bf 1}_{\Omega(t)}\,.\eeq
 In the following, by $B(y,r)$ we denote the open ball centered at $y$ with radius $r$, while
$S^2$ is the sphere of unit vectors in $\R^3$.
The distributional derivatives of $u=u(t,x)$
will be denoted by 
\bel{DDP}D_t u,\quad \qquad D_x u(t,\cdot) \,=\,\nu \cdot \H^1\llcorner \partial \Omega(t).\eeq
Here $\nu= \nu(t,x)\in\R^2$ is the unit inner normal to $\Omega(t)$ at a boundary point $x$, 
while the last expression indicates the restriction of the 1-dimensional Hausdorff measure
to the boundary $\partial \Omega(t)$.

For every set of finite perimeter 
$\Omega\in \F$, its reduced boundary  $\partial^*\Omega$ is defined to be the set of points
$y=(t,x)\in\,]0,T[\,\times\R^2$ such that 
\bel{nuo}
\nu_\Omega(y)~\doteq~\lim_{r\downarrow 0}~ {\mu_\Omega\bigl(B(y,r)\bigr)\over
|\mu_\Omega|\bigl(B(y,r)\bigr)}\eeq
exists in $\R^3$ and satisfies $|\nu_\Omega(y)|=1$. 
The function $\nu_\Omega: \partial^*\Omega\mapsto S^2$ is called the 
{\it generalized inner normal}
to $\Omega$.    A fundamental theorem of De Giorgi \cite{AFP, M} implies that
$\partial^*\Omega$ is countably 2-rectifiable and $|D {\bf 1}_\Omega| = \H^2\rest \partial^*\Omega$.

To formulate an admissibility condition for a set $\Omega\in\F$, 
we observe that, in the smooth case,  the (inward) normal velocity of the set 
$\Omega(t)$  at the point 
$(t,x)\in \partial^*\Omega$ is computed by
\bel{bdef}\beta~=~ {-\nu_0\over \sqrt{\nu_1^2 + \nu_2^2}}\,.\eeq
Recalling that $E(\beta)\doteq\max\{ 1+\beta, 0\}$,  the instantaneous effort is thus computed by
$$\E(t)
~=~\int_{\partial\Omega(t)\cap V} E\left({-\nu_0\over \sqrt{\nu_1^2 + \nu_2^2}}\right)
\,\H^1(dx)
~=~\int_{\partial\Omega(t)\cap V} \max\left\{ {-\nu_0+ \sqrt{\nu_1^2 + \nu_2^2}\over  \sqrt{\nu_1^2 + \nu_2^2}}\,,~0\right\}
\,\H^1(dx).$$
Therefore, integrating over a time interval $t\in \,]t_1,t_2[$ one finds
$$\int_{t_1}^{t_2} \E(t)\, dt ~=~\int_{\partial^*\Omega\cap\{(t,x);\, t_1<t<t_2, \,x\in V\}} 
\max\Big\{ -\nu_0 +\sqrt{\nu_1^2+\nu_2^2}\,,~0\Big\}\, d\H^2\,.$$

To  impose the requirement that $\E(t)\leq M$ for all $t$,    we consider the convex, positively homogeneous
function $L:\R^3\mapsto \R$, defined as
\bel{Ldef} L(v) ~=~L(v_0,v_1,v_2)~\doteq~\max\Big\{ -v_0 +\sqrt{v_1^2+v_2^2}\,,~0\Big\}.\eeq
\begin{definition}\label{def:21} Let  $V\subset\R^2$ be a bounded open set  with finite perimeter.
Given $M>0$, we say that  a set $\Omega\in\F$ in (\ref{Adef}) represents an  
{\bf admissible motion}, and write
$\Omega\in \A$, if for every $0\leq t_1<t_2\leq T$ one has
\bel{adm}
\int_{\partial^*\Omega\cap\{(t,x);\, t_1<t<t_2, \,x\in V\}} L(\nu)\, d\H^2~\leq~M(t_2-t_1).\eeq
\end{definition}

The original problem {\bf (OP)} can now be 
reformulated as
\begi
\item[{\bf (OSM)}] {\bf Optimal Set Motion problem.} {\it Given $T, \kappa_1,\kappa_2>0$ and an initial set $\Omega_0\subseteq V$,
find an admissible set  $\Omega\in\A$
which minimizes the functional
\bel{FU}
\J (\Omega)~\doteq~ \kappa_1 \caL^3(\Omega) + \kappa_2 \, {{\caL}^2}\bigl(\Omega(T)\bigr),\eeq
subject to  $\Omega(0)=\Omega_0$.}
 \endi
 Notice that the sets $\Omega(0)$ and $\Omega(T)$ are well defined in terms of (\ref{om0T}).

\section{Existence of eradication strategies}
\label{s:3}
\setcounter{equation}{0}
We consider here a
null-controllability problem where all sets $\Omega(t)$ remain within a given domain $V\subset \R^2$.
This models the problem of eradicating
a pest population over an island.

\begi
\item[{\bf (CEP)}] {\bf Constrained Eradication Problem.}  
{\it   Consider a bounded open set $V\subset\R^2$  with finite perimeter.
Find a set-valued function $t\mapsto \Omega(t)\subseteq V$  such that, for some $T>0$,
\bel{ncoo}\Omega(0)~=~V,\qquad\qquad 
\Omega(T)~=~\emptyset,\eeq
\bel{ebb}
\int_{V\cap \partial\Omega(t)}E\bigl(\beta(t,x)\bigr)\, \H^1(dx)~\leq~ M,\qquad \forall t\in [0,T].\eeq
}
\endi

Note that in (\ref{ebb}) the effort $E(\beta)$ is integrated only along the portion of the boundary $\partial \Omega(t)$ 
which is contained in the interior of $V$.   Since we are thinking of $V$ as an island, the contamination will never spread outside $V$.

We seek conditions that imply the existence (or non-existence)
of an admissible strategy which eradicates the contamination.
Toward this goal, two geometric invariants associated with the set $V$ play a key role.

\begi\item[(i)] For any given $\lambda\in [0,1]$, we choose a set
$V_\lambda\subset V$ with area ${{\caL}^2}(V_\lambda)= \lambda \,{{\caL}^2}(V)$, so that the length of
its relative boundary is as small as possible.   Then we take the supremum of these lengths over all  $\lambda$.
This leads to the constants
\bel{kl} \kappa(V,\lambda)~\doteq~\inf\Big\{\H^1 \bigl( \partial V_\lambda\cap 
 V\bigr)\,;~~~V_\lambda\subseteq V, ~~{{\caL}^2}(V_\lambda)= \lambda \,{{\caL}^2}(V)\Big\},\eeq
\bel{ko}\kappa(V)~\doteq~\sup_{\lambda\in [0,1]} \kappa(V,\lambda)
.\eeq
\v
\item[(ii)] Next, we   slice the set $V$  in terms of a continuous map $\phi:V\mapsto [0,1]$.
Here the slices are the pre-images $\phi^{-1}(\lambda)$, $\lambda\in [0,1]$.
We choose $\phi$ such that the maximum length of these slices is as small as possible.
This yields the constant
\bel{KO}K(V)~\doteq~\inf_{\phi\,:\,V\mapsto [0,1]} \left(
\sup_{\lambda\in [0,1]}  \H^1\bigl( \phi^{-1}(\lambda)\bigr)\right).
\eeq
\endi

The invariant $K(V)$ in \eqref{KO}  
is known as $1$-$width$ in dimension 2, see \cite{Gromov, Guth}.
It is easy to see that $\kappa(V)\leq K(V)$. In several cases a 
strict inequality   holds.    

\begin{figure}[ht]
\centerline{\hbox{\includegraphics[width=14cm]{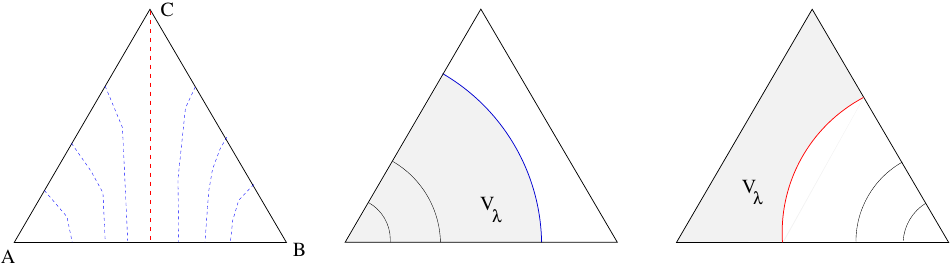}}}
\caption{\small The two invariants  (\ref{ko}) and (\ref{KO}) in the case of an equilateral triangle with unit side. Left: one of the level sets of the function $\phi$  in (\ref{KO}) 
must go through the vertex $C$. Hence it will have length $\geq \sqrt 3/2$.
Center: for any $\lambda\in [0, 1/2]$ we can cut a sector $V_\lambda$ with area 
$\lambda {{\caL}^2}(V)$, so that its boundary is an arc with length $\leq\sqrt{3\sqrt 3/4\pi}$.
Right: for  $\lambda\in [1/2, 1]$ we  can take $V_\lambda$ to be the complement 
of a sector with the same property.
}
\label{f:sc28}
\end{figure}

\begin{example}\label{e:51} {\rm Let $V$ be an equilateral triangle with side of unit length (Fig.~\ref{f:sc28}).
Then 
\bel{Ktri} K(V)~=~ {\sqrt 3\over 2}~\approx~0.866\ldots\eeq
coincides with the height of the triangle.  On the other hand, for 
any $\lambda\in [0,1]$ we have
\bel{ktri}\kappa(V, \lambda) ~\leq~\kappa\Big(V,{1\over 2}\Big)~=~     \sqrt {3\sqrt 3\over 4\pi}~=~\kappa(V)~\approx~0.643\ldots\eeq
According to Theorem~\ref{t:31} below,
the contamination can be eradicated from $V$
if $M>{\sqrt 3\over 2}$, and it cannot be eradicated if $M<   \sqrt {3\sqrt 3\over 4\pi}$.
}
\end{example}

%
\begin{theorem}\label{t:31} {\bf (existence of eradication strategies).}
Let $V\subset\R^2$ be a compact set with finite perimeter. Consider the following statements:
\begi
\item[(i)] There exists a continuous map $\vp: V\mapsto [0,1]$ whose level sets satisfy
\bel{slic}\sup_{s\in [0,1]}~\H^1\Big(\{ x\in V\,;~~\vp(x)=s\}\Big)~<~{M}.\eeq
\item[(ii)] The constrained eradication problem {\bf (CEP)} on $V$ has a solution which satisfies the 
additional monotonicity property
\bel{monp}
s<t\qquad\implies \qquad\Omega(s)\supseteq \Omega(t).\eeq
\item[(iii)] The  problem {\bf (CEP)} on $V$ has a solution. 
\item[(iv)] For every $\lambda\in [0, 1]$ there exists a subset $V_\lambda\subset V$
such that 
\bel{cut}{{\caL}^2}(V_\lambda)~=~\lambda {{\caL}^2}(V),\qquad \H^1(V\cap
\partial V_\lambda)~\leq~{M}\,.\eeq
\endi
Then we have the implications~ (i) $\implies$ (ii) $\implies$ (iii) $\implies$ (iv).
\end{theorem}

\begin{figure}[ht]
\centerline{\hbox{\includegraphics[width=12cm]{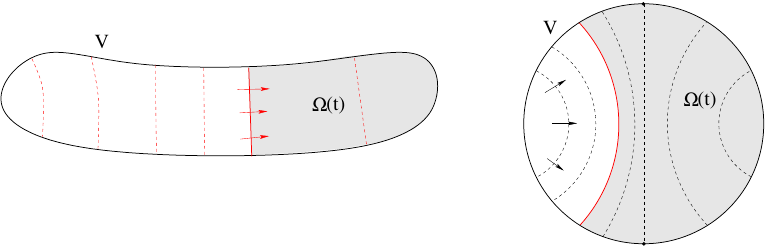}}}
\caption{\small Left: in the case where $M>K(V)$, an  eradication strategy can be constructed 
Right: if (\ref{m2b}) holds and $M<\kappa(V)$, then the area of the set $\Omega(t)$ can never become smaller than
$\lambda^*\caL^2(V)$.
}
\label{f:sc42}
\end{figure}

{\bf Proof.} {\bf 1.} We begin by proving the implication (i)$\implies$(ii).
Assume that  (i) holds (see Fig.~\ref{f:sc42}, left).
As a first step, we observe that the function
\bel{m2s}s~\mapsto~f(s)~\doteq~{{\caL}^2}\Big(\bigl\{ x\in V\,;~\vp(x) \geq  s\bigr\}\Big)\eeq
is continuous.   Indeed, the function $f$ is clearly decreasing.  
If there exists $s_0$ such that
$$f(s_0)~<~\inf_{s< s_0-} f(s),$$
then
$${{\caL}^2}\Big(\bigl\{ x\in V\,;~\vp(x) =  s_0\bigr\}\Big)~=~\inf_{s<s_0} ~{{\caL}^2}\Big(\bigl\{ x\in V\,;~\vp(x) \geq s\bigr\}\Big)
-{{\caL}^2}\Big(\bigl\{ x\in V\,;~\vp(x) \geq  s_0\bigr\}\Big)~>~0.$$
Hence 
$$\H^1\Big(\bigl\{ x\in V\,;~\vp(x) =  s_0\bigr\}\Big)~=~+\infty,$$
reaching a contradiction.
A similar argument rules out the possibility that
$$f(s_0)~>~\sup_{s> s_0-} f(s).$$
\v

{\bf 2.} 
To construct a strategy satisfying (\ref{ncoo})-(\ref{ebb}) and  (\ref{monp}), we proceed as follows.
By (\ref{slic}) 
we can find $b_1>0$ small enough such that 
\bel{m19}
\H^1\Big(\bigl\{ x\in V\,;~~\vp(x)=s\bigr\}\Big)~<~{M\over 1 + b_1}
\qquad\qquad \forall s\in [0,1].
\eeq
We now define
\bel{omet}\Omega(t) ~=~\bigl\{ x\in V\,;~\vp(x) \geq  s(t)\bigr\},\eeq
for a suitable function $t\mapsto s(t)$, with $s(0)=0$, $s(T)=1$ for some time $T>0$.
This function $s(\cdot)$ is chosen so that
\bel{simp}{\caL}^2\bigl(\Omega(t)\bigr)~\doteq~{{\caL}^2}\Big(\bigl\{ x\in V\,;~\vp(x) \geq  s(t)\bigr\}\Big)~=~{{\caL}^2}(V)- {b_1M\over 1+b_1} t.\eeq
\v
{\bf 3.} It remains to prove that the above multifunction  $t\mapsto\Omega(t)$ satisfies all requirements.

By construction, the property (\ref{monp}) is trivially satisfied.   This implies $\beta(t,x)\geq 0$
for all $t$ and all $x\in V\cap\partial \Omega(t)$.  For any $0\leq t_1<t_2\leq 1$, using (\ref{m19}) 
and then (\ref{simp}) we  compute
$$\bega{l}\ds\int_{t_1}^{t_2} \int_{V\cap\partial \Omega(t)} E\bigl(\beta(t,x)\bigr)\, \H^1(dx)\, 
dt ~=~\int_{t_1}^{t_2}
\int_{V\cap\partial \Omega(t)}\Big(
1+\,\beta(t,x)\Big) \H^1(dx)\, dt\\[4mm]
\qquad \ds \leq~{M\over 1+b_1} (t_2-t_1)+  \Big[ {{\caL}^2}\bigl(\Omega(t_1)\bigr) -  {{\caL}^2}\bigl(\Omega(t_2)\bigr)
\Big]\\[4mm]
\qquad\ds  =~{M\over 1+b_1} (t_2-t_1) + {b_1 M\over 1+b_1} (t_2-t_1) ~=~M(t_2-t_1).
\enda $$
Since this is true for all $t_1<t_2$, this yields (\ref{ebb}).
\v
{\bf 4.} The implication  (ii)$\implies$(iii) is trivial.
\v
{\bf 5.} To prove the implication  (iii)$\implies$(iv), assume that, on the contrary,  (iv) fails
 (see Fig.~\ref{f:sc42}, right).  Hence there exists $0< \lambda^*<1$ and $\ve>0$ such that
\bel{m2a}
\H^1(V\cap\partial W)~>~M+2\ve \qquad\forall W\subset V~~\hbox{with}~~{{\caL}^2}(W)= \lambda^* {{\caL}^2}(V).
\eeq
We claim that there exists $\delta>0$ such that, 
\bel{m2b}
\H^1(V\cap \partial W)~\geq ~M+\ve
\qquad\forall W\subset V~~\hbox{with}~~{{\caL}^2}(W)= \lambda {{\caL}^2}(V), 
~~|\lambda-\lambda^*|\leq\delta.\eeq
Indeed, if (\ref{m2b}) fails, we could find a sequence of subsets
$W_n\subset V$ such that 
\bel{m2c} \lim_{n\to\infty} {{\caL}^2}(W_n)~=~ \lambda^* {{\caL}^2}(V),
\qquad\quad \limsup_{n\to\infty}~
\H^1(V\cap\partial W_n)~\leq~M+\ve.\eeq
Since all these sets have uniformly bounded perimeters, taking a  subsequence we can assume the 
$\L^1$ convergence of the characteristic functions: ${\bf 1}_{W_n}\to {\bf 1}_{W^*}$, for a set
$W^*\subset V$ such that 
$$\H^1(V\cap\partial W^*)~\leq ~{M+\ve}\qquad~~{{\caL}^2}(W^*)= \lambda^* {{\caL}^2}(V).$$
This yields a contradiction with (\ref{m2a}).
\v
{\bf 6.} Using (\ref{m2b}), we claim that, for any admissible strategy 
$t\mapsto \Omega(t)$ one has the implication
$${{{\caL}^2}(\Omega(t))\over {{\caL}^2}(V)}~\in~[\lambda^*-\delta, \lambda^*+\delta]\qquad  \implies
\qquad
{d\over dt} {{\caL}^2}(\Omega(t))~>~0.$$
Indeed, 
$$\bega{rl} \ds{d\over dt} {{\caL}^2}(\Omega(t))&\ds=~ -\int_{\partial \Omega(t)\cap V}  \beta(t,x) \,\H^1(dx)~=~
\int_{\partial \Omega(t)\cap V}  \Big[1-E\bigl(\beta(t,x)\bigr)\Big] \,\H^1(dx)\\[4mm]
\qquad\qquad\ds &\geq ~\H^1\bigl(\partial \Omega(t)\cap V
\bigr)- M~\geq~\ve.\enda$$
%
%
As a consequence, the area ${{\caL}^2}(\Omega(t))$ can never become smaller than $\lambda^* {{\caL}^2}(V)$.
In particular, this area cannot decrease to zero.
\endproof

\v
\section{A sufficient condition for optimality}
\label{s:4}
\setcounter{equation}{0}
In this section we consider the minimum time eradication problem {\bf (MTP)}, constrained to a
set $V\subset\R^2$
with finite perimeter. The following result provides a sufficient condition for a set motion 
$t\mapsto \Omega(t)$ to be optimal for the minimum time problem.
\begin{theorem} 
\label{t:41} 
Let $t\mapsto \Omega(t)\subseteq V$ be an admissible set motion
with $\Omega(0)=V$, $\Omega(T)=\emptyset$, and such that, for all $0<t<t'<T$,
\bel{da1} \caL^2\bigl(\Omega(t')\bigr)\,<\,\caL^2\bigl(\Omega(t)\bigr),\qquad\qquad \Omega(t')\,\subseteq\, B_{t'-t} \bigl(\Omega(t)\bigr),\eeq
\bel{da2}
 {d\over dt} \caL^2\bigl(\Omega(t)\bigr)\,=\,\H^1\bigl(\partial \Omega(t)\cap V\bigr)-M.\eeq
Assume that
\bel{opc}\H^1 \bigl(V\cap \partial \Omega(t)\bigr)\,=\,\min\Big\{ \H^1 (V\cap \partial W)\,;
\qquad
W\subseteq V,\quad {{\caL}^2}(W) = {{\caL}^2}\bigl(\Omega(t)\bigr)\Big\}\eeq
 for all $t\in [0,T]$.   Then this motion is optimal for the minimum time problem.
\end{theorem}
\v
{\bf Proof.} 
Let $A= {{\caL}^2}(V)$ be the area of the set $V$.  For every $a\in [0, A]$ define
\bel{gdef}g(a)~\doteq~\min\{ \H^1(V\cap \partial W)\,;~~~W\subseteq V,\quad 
{{\caL}^2}(W)= a\Big\}.\eeq
Consider any admissible strategy $t\mapsto \Tilde\Omega(t)$, $t\in [0, \Tilde T]$.  We then have
\bel{a4}{d\over dt} {{\caL}^2}(\Tilde \Omega(t))~\geq~\H^1\bigl(V\cap \partial \Tilde 
\Omega(t)\bigr)-M~\geq~g\Big( {{\caL}^2}\bigl(\Tilde\Omega(t)\bigr)\Big) -M.\eeq
Calling $\Tilde a(t)\doteq {{\caL}^2}\bigl(\Tilde\Omega(t)\bigr)$, for every $t$ one has
$${d\over dt} \Tilde a(t)~\geq~g\bigl(\Tilde a(t)\bigr) -M,$$
hence
\bel{a5}
t~\geq~\int_{\Tilde a(t)}^A {da\over M-g(a)}\,.\eeq
In particular, if $\Tilde\Omega(0)=V$ and $\Tilde \Omega(\Tilde T)=\emptyset$, this implies
\bel{a6}
\Tilde T~\ge ~\int_0^A {da\over M-g(a)}\,.\eeq
On the other hand, if the strategy $t\mapsto \Omega(t)$ satisfies (\ref{da1})--(\ref{opc}),
then 
\bel{a7}
T~=~\int_0^A {da\over M-g(a)}~\leq~\Tilde T.\eeq
Hence $\Omega(\cdot)$ is optimal for the minimum time problem.
\endproof
\v
\begin{figure}[ht]
\centerline{\hbox{\includegraphics[width=9cm]{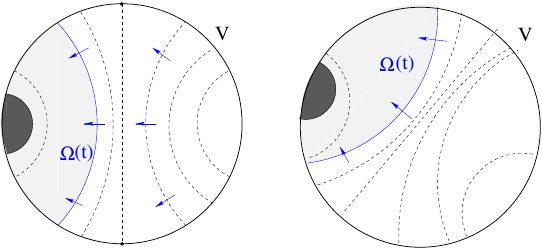}}}
\caption{\small Two optimal solutions for the minimum time problem on a disc. For every $t\in \,]0,T[\,$,
the  boundary $\partial \Omega(t)$  should be an arc of circumference, crossing the boundary 
$\partial V$ perpendicularly.}
\label{f:sc40}
\end{figure}

\begin{remark} {\rm
Notice that the minimality condition in (\ref{opc}) implies that, for each
$t\in \,]0,T[\,$, the set $\Omega(t)$ should provide a solution to 
the classical Dido's problem
\cite{Dido}.  As a consequence, the boundary $\partial \Omega(t)$ should be an arc of circumference, crossing the  boundary 
$\partial V$ perpendicularly at both endpoints
(see Fig.~\ref{f:sc40}).}
\end{remark}

More generally, the  arguments used in the proof of Theorem~\ref{t:41} can be used to obtain necessary and
sufficient conditions for a family of {\bf minimum time transfer problems} with prescribed initial and terminal conditions.

\begi
\item[{\bf (MTTP)}] {\it 
Given a bounded open set $V\subset \R^2$ and two subsets
$$\Omega_{final}~\subset ~\Omega_{initial}\subseteq V,$$
consider the family of all admissible motions $t\mapsto \Omega(t)\subseteq V$ such that 
\bel{infin}\Omega(0)=\Omega_{initial},\qquad\qquad \Omega(T)=\Omega_{final}\eeq
for some $T>0$.   Among all such motions, minimize the time $T$.}
\endi

\begin{corollary}\label{c:41}
In connection with the optimization problem {\bf (MTTP)}, assume that there exists an admissible
motion $t\mapsto\Omega^*(t)$, $t\in [0,T^*]$  that satisfies (\ref{infin}) together with (\ref{da2})-(\ref{opc}), for every
$t$.    Then $\Omega^*$ is optimal.   Moreover, any other optimal motion must satisfy (\ref{opc}) as well.
\end{corollary}

\begin{example}\label{e:61} {\rm In the case where $V$ is  disc with diameter $\hbox{diam}(V)<M$, the two constants $\kappa(V)$ and $K(V)$ in (\ref{ko})-(\ref{KO}) both coincide with the diameter of $V$.
In this setting,  two  strategies  satisfying the 
sufficient conditions stated in Theorem~\ref{t:41} are shown in Fig.~\ref{f:sc40}.  
On the other hand, if $\hbox{diam}(V)>M$, by Theorem~\ref{t:31} the contamination cannot be eradicated. 

In the borderline case where $\hbox{diam}(V)=M$,  the contamination still cannot be eradicated. 
Indeed, recalling the function $g$ at (\ref{gdef}),  
for any admissible strategy $t\mapsto \Omega(t)$, if $ \Omega( T)=\emptyset$, this implies
\bel{a8}
 T~\ge ~\int_0^{ \pi M^2/ 4} {da\over M-g(a)}~=~{2\over M}\int_0^\pi  {d\theta\over
 \ds  1-\frac{\pi-\theta}{2}\tan\Big({\theta\over 2}\Big)}~=~+\infty\,.\eeq}
 \end{example}

\section{Existence of optimal solutions}
\label{s:5}
\setcounter{equation}{0}

In the case without geographical constraints, the existence of solutions for a general class of
optimal set motion problems was proved in \cite{BCS2}. The same arguments can be adapted to 
the present setting.

\begin{theorem}\label{t:51}
Let $V\subset \R^2$ be a bounded open set with finite perimeter. 
Then for any  set $\Omega_0\subseteq V$ with finite perimeter and any $T, \kappa_1,\kappa_2>0$,
the set motion problem {\bf (OSM)} has an optimal solution.
\end{theorem}

{\bf Proof.} 
{\bf 1.} We observe that the functional $\J(\Omega)$ at (\ref{FU}) is non-negative.   Moreover,
the strategy $t\mapsto \Omega(t)\doteq B_t (\Omega_0)\cap V$ is admissible and has finite cost.
Indeed, setting
\bel{omo}\Omega~\doteq~\Big\{ (t,x)\,;~x\in B_t(\Omega_0)\cap V, ~~t\in [0,T]\Big\},\eeq
one has
$$\J(\Omega)~\leq~\kappa_1\int_0^T \caL^2\bigl(B_t(\Omega_0)\bigr)\, dt + \kappa_2 \caL^2\bigl(B_T(\Omega_0)\bigr)
~<~+\infty.$$
Moreover, by (\ref{omo}) the inner unit normal satisfies
$$-\nu_0 +\sqrt{\nu_1^2+\nu_2^2} ~=~0$$
at a.e.~boundary point $(t,x)\in \Omega\cap \bigl(]0,T[\times V\bigr)$.  Hence the  admissibility condition in Definition~\ref{def:21} is trivially satisfied: $\Omega\in \A$.

We can thus consider 
a minimizing sequence of sets  $\Omega_n\in \A$ 
such that, as $n\to\infty$,
$$\J(\Omega_n)~\to~\J_{min}~\doteq~\inf_{\Omega\in \A} \J (\Omega)\,.$$
Without loss of generality we can assume that
\bel{52}\Omega_n(t)~\doteq~\bigl\{x\in \R^2\,; ~(t,x)\in\Omega_n\bigr\}~\subseteq~B_t( \Omega_0)\cap V\eeq
for all $t\in [0,T]$ and $n\geq 1$.
Otherwise, we can simply replace 
each set $\Omega_n(t)$ with the intersection
$\Omega_n(t)\cap B_t( \Omega_0)$, without increasing 
the total cost.
\v
{\bf 2.} We now prove a uniform bound on the perimeters of the 
sets $\Omega_n\subset \R^3$.  
For every $n\geq 1$ we split the reduced boundary  in the form
\bel{spn}
\partial^*\Omega_n \cap \bigl(]0,T[\times V\bigr)~=~\Sigma_n^-\cup\Sigma_n^+\,,\eeq
so that the following holds.
Calling $\nu=(\nu_0, \nu_1,\nu_2)$ the unit inner normal vector at the point
$(t,x)\in \partial^*\Omega_n\cap V$, and defining 
 the inner normal velocity $\beta_n=\beta_n(t,x)$ as in 
(\ref{bdef}), one has
\bel{Snpm}
\left\{ \bega{rl} \beta_n(t,x)~\leq ~-1/2\qquad &\hbox{if}\quad (t,x)\in \Sigma_n^-\,
,\\[2mm] \beta_n(t,x)~>~-1/2\qquad &\hbox{if}\quad (t,x)\in \Sigma_n^+\,.
\enda\right.\eeq
Integrating the effort over $\Sigma_n^+$ one obtains 
\bel{Sn+}
MT~\geq~\dint_{\Sigma_n^+} \max\Big\{ -\nu_0 +\sqrt{\nu_1^2+\nu_2^2}\,,~0\Big\} \,d\H^2 ~\geq~\dint_{\Sigma_n^+} {1\over 2}  \,d\H^2 \,.\eeq
Hence
\bel{H2S+}\H^2\bigl( \Sigma_n^+\bigr)~\leq~2MT.\eeq
On the other hand,
on the domain $\Sigma_n^-$ 
one has the lower bound
\bel{n00}
\nu_0~\geq~{1\over 2}.\eeq
We can now write
\bel{Sm1}\bega{l}\ds \caL^2\Big( B_T(\Omega_0)\Big)~\geq~
\caL^2\bigl(\Omega_n(T)\bigr) - \caL^2\bigl(\Omega_n(0)\bigr) ~=~\int_0^T \int_{\partial \Omega(t)} -\beta(t,x) \H^1(dx)\, dt   \\[4mm]
 \quad \ds=~\int_0^T \int_{\partial \Omega(t)} {\nu_0\over\sqrt{\nu_1^2+\nu_2^2}} \H^1(dx)\, dt  
=~\int_{\Sigma_n^+\cup\Sigma_n^-} \nu_0\,d\H^2
~ \geq~{1\over 2}  \int_{\Sigma_n^-}
d\H^2 - \int_{\Sigma_n^+} d\H^2.
\enda
\eeq
Combining (\ref{Sm1}) with (\ref{H2S+}) one obtains
\bel{Sm2}\H^2\bigl(\Sigma_n^-\bigr) ~\leq~2  \caL^2\Big( B_T(\Omega_0)\Big) + 4MT.
\eeq
Together, the two inequalities (\ref{H2S+}) and (\ref{Sm2})  yield a
uniform bound on the 2-dimensional measure $\H^2\Big(\partial^* \Omega_n \cap \bigl(]0,T[\times V\bigr)\Big)$
of the relative boundary of $\Omega_n$ inside $V$.   

Since $V$ has finite perimeter, this implies that the characteristic functions ${\bf 1}_{\Omega_n}$ have uniformly 
bounded variation.
\v
{\bf 3.} 
By  the uniform BV bound proved in the previous step, by possibly taking a subsequence,
a compactness argument (see Theorem~12.26 in \cite{M}) yields the existence of a 
set with finite perimeter $\Omega\in \F$ such that the following holds.
As $n\to\infty$, one has the convergence 
\bel{cvg}\Big\| {\bf 1}_{\Omega_n} - {\bf 1}_\Omega\Big\|_{\L^1\bigl(]0,T[\,\times\R^2\bigr)}~\to~0,\eeq
together with the weak convergence of measures
\bel{wco}
\mu_{\strut\Omega_n}~{^*}\!\!\!\!\!\rightharpoonup~\mu_{\strut\Omega}\,.\eeq
\v
{\bf 4.} In this step we check that the limit strategy $\Omega\in\F$ is admissible.
Since the function 
$L$ in (\ref{Ldef}) is convex, for any $0\leq t_1<t_2<T$
we can use a lower semicontinuity result for anisotropic functionals
 (see Theorem 20.1 in \cite{M})  and conclude 
\bel{lsc}\bega{l} \ds
\int_{\partial^*\Omega\cap\{ t_1<t<t_2, x\in V\}}   L(\nu)\, d\H^2
~\leq~\liminf_{n\to\infty} 
\int_{\partial^*\Omega_n\cap\{ t_1<t<t_2, x\in V\}}  L(\nu_n)\, d\H^2~\leq~M(\tau'-\tau).
\enda \eeq
Therefore $\Omega\in\A$.
\v
{\bf 5.} It remains to check that the limit set $\Omega$ is optimal. As $n\to\infty$, the convergence (\ref{cvg}) immediately implies
\bel{cm3}\caL^3(\Omega_n)~\to~\caL^3(\Omega).\eeq
    Moreover, since the map $t\mapsto {\bf 1}_{\Omega(t)}\in \L^1(\R^2)$ has bounded variation, given $\ve>0$ we can find $\delta>0$ such that 
\bel{m2o}
\left| \caL^2(\Omega(T)\big)-  {1\over\delta} \int_{T-\delta}^T \caL^2(\Omega(t)\big)\, dt
\right|~<~\ve.\eeq
On the other hand,  for every $n\geq 1$ and $t<T$ we have the one-sided estimate 
$$\caL^2\bigl(\Omega_n(t)\setminus\Omega_n(T)\bigr)~\leq~ M(T-t).
$$
Therefore, choosing $\delta < \ve/M$ we obtain
$$\caL^2\bigl(\Omega_n(T)\bigr)~\geq~{1\over \delta} \int_{T-\delta}^T \Big[\caL^2
\bigl(\Omega_n(t)\bigr)-M(T-t)\Big]\, dt ~\geq~{1\over \delta} \int_{T-\delta}^T \caL^2
\bigl(\Omega_n(t)\bigr)\, dt -\ve$$
for every $n\geq 1$.  Using the convergence
$${1\over \delta} \int_{T-\delta}^T \caL^2
\bigl(\Omega_n(t)\bigr)\, dt ~\to~ {1\over\delta} \int_{T-\delta}^T \caL^2
\bigl(\Omega(t)\bigr)\, dt \,,$$
we obtain
\bel{lsom}
\caL^2\bigl(\Omega(T)\bigr) ~\leq~\liminf_{n\to\infty}  \caL^2\bigl(\Omega_n(T)\bigr)  + 2\ve.\eeq
Combining (\ref{cm3}) and (\ref{lsom}),
since $\ve>0$ was arbitrary we conclude 
\bel{linf}
\J(\Omega)~\doteq~
\kappa_1\, \caL^3(\Omega) + \kappa_2 \, \caL^2\bigl(\Omega(T)\bigr) ~\leq~\liminf_{n\to \infty} \J(\Omega_n).\eeq
\v
{\bf 6.} It remains to prove that the  initial 
condition $\Omega(0)=\Omega_0$ is satisfied.   
Since $\Omega_0$ has bounded perimeter, by (\ref{52}) we immediately have
\bel{li4}\lim_{t\to 0+} \caL^2\bigl(\Omega(t)\setminus\Omega_0\bigr)~\leq~\lim_{t\to 0+} 
\caL^2\Big(B_t(\Omega_0) \setminus\Omega_0\Big) ~=~0.\eeq
On the other hand,  for every $n\geq 1$ and $t>0$ the admissibility condition $\Omega_n\in\A$
implies
$$\caL^2\bigl( \Omega_0\setminus \Omega_n(t)\bigr)  ~\leq~Mt.
$$
Taking the limit as $n\to\infty$ one obtains
\bel{li9} 
\caL^2\bigl(\Omega_0\setminus \Omega(t)\bigr)~ \leq~Mt.\eeq
Together, (\ref{li4}) and (\ref{li9}) yield the convergence ${\bf 1}_{\Omega(t)}\to 
{\bf 1}_{\Omega_0}$ in $\L^1(\R^2)$, completing the proof.
\endproof

Entirely similar arguments yield the existence of an optimal solution for the minimum time problem.

\begin{theorem}\label{t:52}
Let $V\subset \R^2$ be a bounded open set with finite perimeter, and let $M>0$ be given.
If the constrained 
eradication problem {\bf (CEP)} has a solution, then the minimum time problem {\bf (MTP)} has an optimal solution.
\end{theorem}

{\bf Proof.} By the assumption, there exists a minimizing sequence $(\Omega_n)_{n\geq 1}$, with $\Omega_n\in\A$,
$\Omega_n(0)=V$,~$\Omega_n(T_n)=\emptyset$   for all $n\geq 1$. Here $T_n\to T$ converge to the infimum
among all eradication times.

Since all sets $\Omega_n$ are admissible,
by possibly taking a subsequence, the same arguments used in the proof of Theorem~\ref{t:51} yield the convergence
(\ref{cvg})-(\ref{wco}).  By (\ref{lsc}) we again conclude that $\Omega\in\A$.
 
Since $\Omega_n(t)\subseteq V$ for all $t\geq 0$, from  the inequality
$$\caL^2\bigl(V\setminus \Omega_n(t)\bigr)~\leq~Mt, $$
letting $n\to \infty$ one obtains $\Omega(0)=V$.

Finally, using the fact that  $T<T_n$ and $\Omega_n(T_n)=\emptyset$, since every $\Omega_n$ is admissible 
we obtain
$$\caL^2\bigl(\Omega_n(T)\setminus \Omega_n(T_n)\bigr)~=~\caL^2\bigl(\Omega_n(T)\bigr)~\leq~M(T_n-T).$$
Letting $n\to\infty$ we conclude that $\caL^2\bigl(\Omega(T)\bigr)=0$, hence $\Omega(T)=\emptyset$.
This completes the proof.
 \endproof

\section{Necessary conditions for optimality}
\label{s:6}
\setcounter{equation}{0}

In order to derive a set of necessary conditions for optimality, some regularity
conditions on the domain $V$ and on the moving sets $\Omega(t)$
will be needed.
A key assumption is  that the unit outer normal $\bfn(x)$ to a boundary
point $x\in \partial \Omega(t)\cap V$ depends Lipschitz continuously on  $x$.
This allows us to parameterize the relative boundary $\partial \Omega(t)\cap V$ 
as $\xi\mapsto x(t,\xi)$, in such a way that the maps $t\mapsto x(t,\xi)$ describe 
perpendicular curves.

In the following, $S^1\doteq\{\xi\in\R^2\,;~|\xi|=1\}$ denotes the unit circumference.
Given a vector $v=(v_1, v_2)\in\R^2$, we write
$v^\perp= (-v_2, v_1)$ for  the perpendicular vector, obtained by a rotation of $\pi/2$.
Moreover, $\C^{1,1}$ denotes a space of continuously differentiable functions
 with locally Lipschitz continuous partial derivatives.

\begin{figure}[ht]
\centerline{\hbox{\includegraphics[width=13cm]{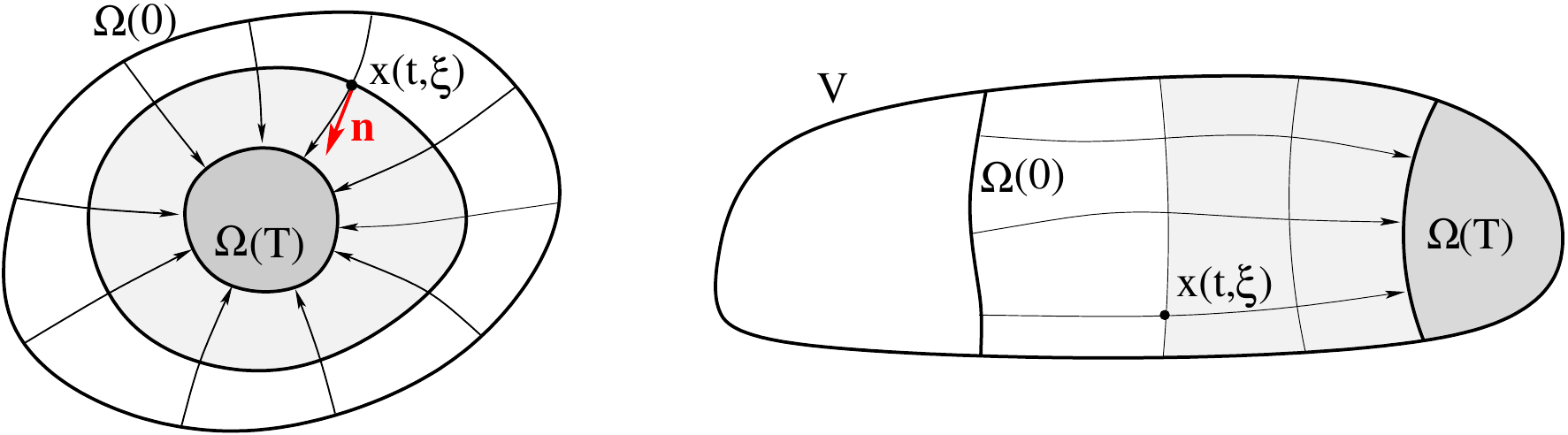}}}
\caption{\small Two examples where the boundaries of the sets $ \Omega(t)$ can be 
parameterized as $\xi\mapsto x(t,\xi)$, according to {\bf (A1)--(A3)}. Left:  a case where $(t,\xi)\in W=[0,T]\times S^1$.  Right: a case  with
geographical constraints, where  $(t,\xi)\in W=[0,T]\times ]0,1[\,$.}
\label{f:sm56}
\end{figure}

We shall assume (see Fig.~\ref{f:sm56}): 
\begi
\item[{\bf (A1)}]  {\it Either $V=\R^2$, or  $V\subset \R^2$ is a bounded open set with 
piecewise $\C^1$  boundary.}
\endi
\begi
\item[{\bf (A2)}]  {\it  The relative boundaries of the sets $\Omega(t)$, $t\in [0,T]$, admit a $\C^{1,1}$ 
parameterization of the form
\bel{xtxi} W~\ni~(t,\xi)~\mapsto ~x(t,\xi)~\in ~\partial \Omega(t)\cap V\eeq 
 such that the following holds.
 \begi
\item[(i)] Either $W = [0,T]\times S^1$ if $V=\R^2$,  or $W=[0,T]\times \,]0,1[\,$ in the case with geographical constraints.
 \item[(ii)]
At each time $0<t<T$, the map 
\bel{part}
\xi~\mapsto ~x(t,\xi)~\in~\partial \Omega(t)\cap V\eeq
 is one-to-one.  Its range covers all the relative  boundary $\partial \Omega(t)\cap V$.
 \item[(iii)]  For every $(t,\xi)\in W$ with $t>0$, the partial derivative
$ x_\xi(t,\xi)$
is a nonzero tangent vector to the boundary $\partial\Omega(t)$ at the point $x(t,\xi)$.
For any given $\xi$ and any $\tau>0$ there holds
\bel{ixi} \min_{t\in [\tau,T]} ~\bigl| x_\xi(t,\xi)\bigr|~>~0.\eeq
\item[(iv)]
The perpendicular vector 
\bel{ndef}
\bfn(t,\xi)~\doteq~\left({x_\xi(t,\xi)\over \bigl| x_\xi(t,\xi)\bigr|}\right)^\perp\eeq
yields
the unit inner normal  to the set $ \Omega(t)$ at the boundary point $x(t,\xi)$.
\item[(v)] For each $\xi$,  the trajectory
$t\mapsto x(t,\xi)$ is orthogonal to the boundary $\partial\Omega(t)$ at every time $t\in [0,T]$.
Namely, there exists a continuous, scalar function $\beta:W\mapsto \R$ such that 
\bel{xt3}x_t(t,\xi)~=~\beta(t,\xi)\,\bfn(t,\xi)\qquad\qquad\forall (t,\xi)\in W.\eeq
\endi
}
\endi
We denote by 
\bel{curv}
\omega(t,\xi)~\doteq~{\la \bfn(t,\xi),\, x_{\xi\xi}(t,\xi)\ra\over \bigl|x_\xi(t,\xi)\bigr|^2}\eeq
the curvature of the boundary $\partial \Omega(t)$ at the point $x(t,\xi)$.
Notice that this curvature is well defined for a.e.~$(t,\xi)\in W$. Indeed, the functions
$x_\xi$ and $\bfn$ are locally Lipschitz continuous,  while $|x_\xi|$ is a continuous, strictly positive function.
By Rademacher's theorem~\cite{EG}, $x_{\xi\xi}$ exists almost everywhere.
The following regularity property will be assumed:
\begi
\item[{\bf (A3)}] {\it For each $\xi$, the curvature function $t\mapsto \omega(t,\xi)$ 
is measurable and bounded. Moreover,}
\bel{lepo} \lim_{\ve\to 0} ~\int_0^T \sup_{|\zeta-\xi|<\ve} \bigl|\omega(t,\zeta)-\omega(t,\xi)\bigr|\,dt~=~0.\eeq
 \endi

\begin{figure}[ht]
\centerline{\hbox{\includegraphics[width=10cm]{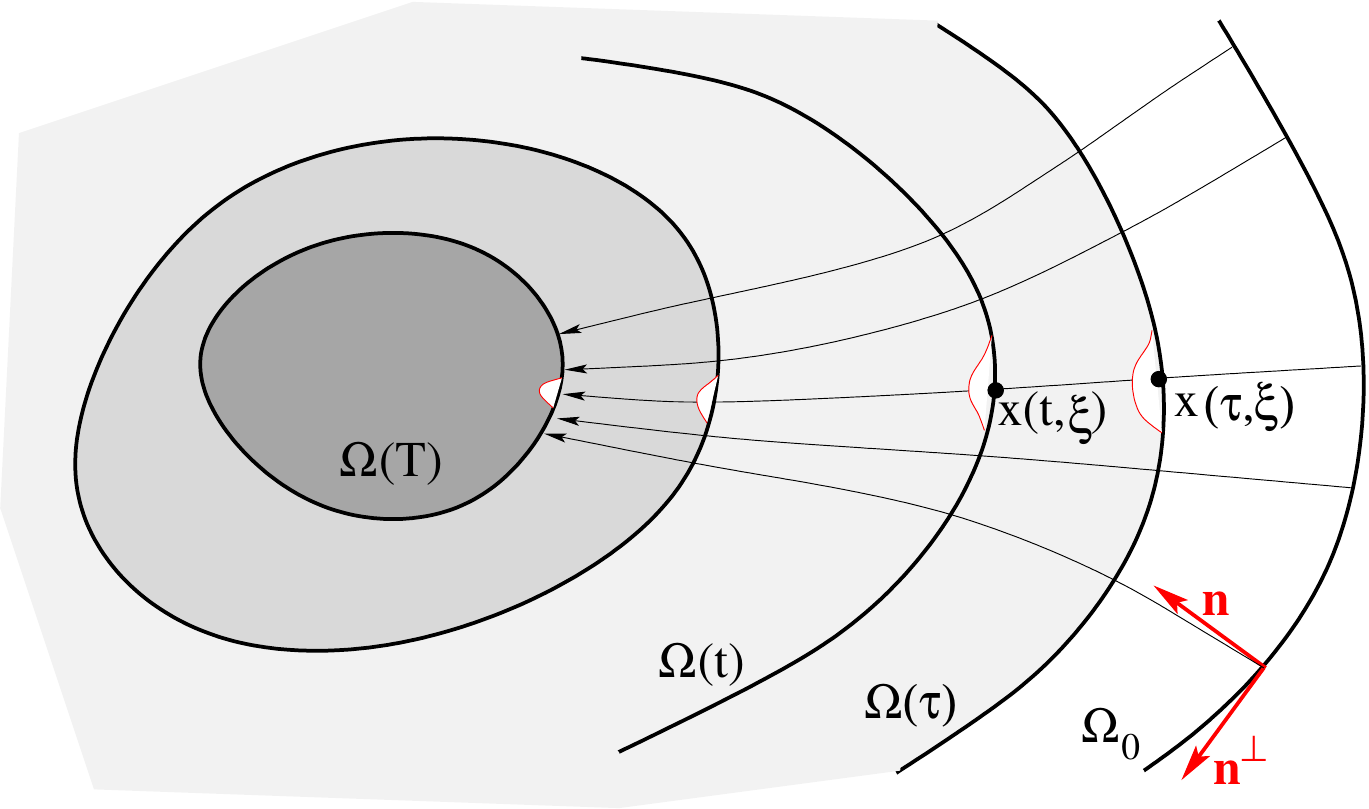}}}
\caption{\small   At each time $t\in [0,T]$, the boundary of the set $\Omega(t)$
is parameterized by $\xi\mapsto x(t,\xi)$. The curves $t\mapsto
x(t,\xi)$ cross the boundaries perpendicularly.   }
\label{f:sm42}
\end{figure}

Before stating a precise set of necessary conditions, we 
outline the main ideas. 
As shown in Fig.~\ref{f:sm42}, assume that at some time $\tau\in [0,T]$
we change our control $\beta(\cdot)$ in a neighborhood of a point $(\tau, \xi)$.
This corresponds to a ``needle variation" used in the proof of the Pontryagin Maximum Principle
\cite{BPi, Ce, FR}.
However, in the present case the perturbation is localized in time and also in space. 
As a result of this perturbation, at time $\tau$ the set $ \Omega(\tau)$ will be replaced by 
a (possibly smaller) set $\Omega^\ve(\tau)$.  Its relative boundary will be parameterized by
$$\xi\mapsto x^\ve(\tau,\xi)~\in~\partial \Omega^\ve(\tau)\cap V.$$

Afterwards, for  $t\in [\tau, T]$, we keep the same control effort
as before at all points of the boundary $ \partial \Omega^\ve(t)$.
More precisely, calling 
$$\beta^\ve(t,x)~\doteq~\la x^\ve_t(t,\xi),\,\bfn^\ve(t,\xi)\ra$$
the inward speed of the boundary point  $x^\ve(t,\xi)\in \partial \Omega^\ve(t)$,
we impose
\bel{sameff} \beta^\ve(t,x)\, \bigl| x^\ve_\xi(t,\xi)\bigr|~=~\beta(t,x)\, \bigl| x_\xi(t,\xi)\bigr|.\eeq
This will guarantee that the total effort in the perturbed strategy is admissible.  Namely, for every $t\in \,]\tau,T[\,$,
\bel{effe}\E^\ve(t)~=~\int  E\bigl(\beta^\ve(t,x)\bigr) \, \bigl| x^\ve_\xi(t,\xi)\bigr|\, d\xi~=~
\int  E\bigl(\beta(t,x)\bigr) \, \bigl| x_\xi(t,\xi)\bigr|\, d\xi~=~\E(t)~=~M.\eeq

\begin{figure}[ht]
\centerline{\hbox{\includegraphics[width=8cm]{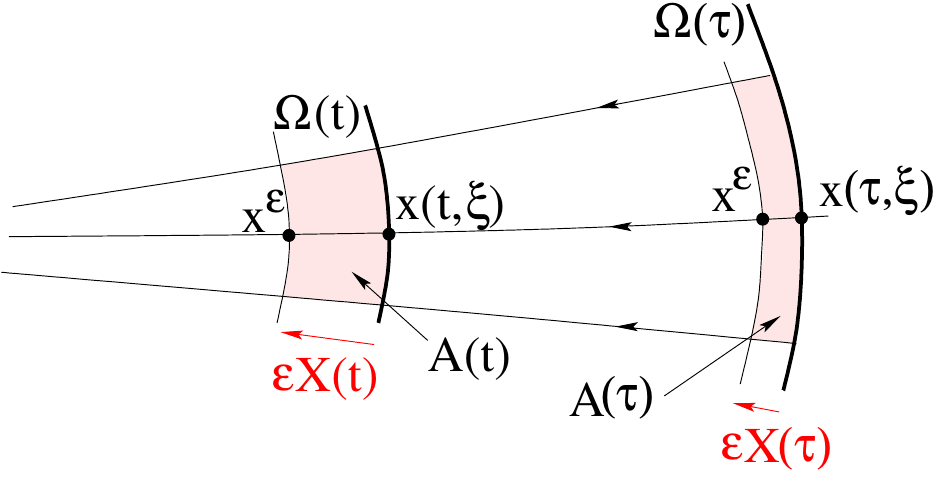}}}
\caption{\small  
In the perturbed motion, at time $t$ the boundary $\partial\Omega(t)$ is pushed inward in the amount $\ve X(t)$.
To leading order, 
the  change $A(t)$ in the area is described by the ODE  (\ref{pat}). }
\label{f:sm77}
\end{figure}

We wish to estimate the reduction in the total cost, due to the fact that all sets $\Omega^\ve(t)$,
$t\in [\tau, T]$ are now smaller than in the original optimal solution.
Toward this goal, consider the first order approximation
\bel{pe1}x^\ve(t,\xi)~=~x(t,\xi)+ \ve X(t,\xi)\bfn(t,\xi) + o(\ve).\eeq
 This implies 
\bel{pe2} \bigl| x^\ve_\xi(t,\xi)\bigr| ~\approx~ \bigl| x_\xi(t,\xi)\bigr| \bigl(1-\omega(t,\xi)\ve X(t,\xi)\bigr),\eeq
where $\omega$ is the boundary curvature.
In view of (\ref{sameff}) one obtains
\bel{bep}\beta^\ve(t,\xi)~\approx~{\bigl( 1+\beta(t,\xi)\bigr) \bigl|x_\xi(t,\xi)\bigr| \over 
\bigl|x^\ve_\xi(t,\xi)\bigr|}-1  ~\approx~\bigl(1+\beta(t,\xi)\bigr)  \bigl(1+\omega(t,\xi)\ve X(t,\xi)\bigr)-1.\eeq
To leading order we thus find
\bel{doX} \partial_t X(t,\xi)~=~{\beta^\ve(t,\xi)-\beta(t,\xi)\over\ve}~\approx~ \bigl(1+\beta(t,\xi)\bigr)\omega(t,\xi) X(t,\xi).\eeq
Calling $A(t,\xi)= X(t,\xi)\, \bigl|x_\xi(t,\xi)\bigr|$ the infinitesimal decrease in the area produced by the perturbation
(see Fig.~\ref{f:sm77}), differentiating w.r.t.~time $t\in [\tau,T]$, we find
$${\partial_t}\, \bigl|x_\xi(t,\xi)\bigr|~=\,-\beta(t,\xi) \omega(t,\xi)
\bigl|x_\xi(t,\xi)\bigr| ,$$
\bel{pat} 
\bega{rl}{\partial_t}\,  A(t,\xi)&\approx~ \bigl(1+\beta(t,\xi)\bigr)\omega(t,\xi) X(t,\xi) \bigl|x_\xi(t,\xi)\bigr| - X(t,\xi) \beta(t,\xi)
\omega(t,\xi)
\bigl|x_\xi(t,\xi)\bigr| \\[3mm]
&=~\omega(t,\xi) A(t,\xi).\enda\eeq
The linear ODE (\ref{pat}) is the key to understanding the optimality condition.
If at time $\tau$ we are able to clean up an additional area $\bar a$ in a neighborhood of the point $x(\tau,\xi)$, 
the total cost would decrease in the amount $\bar a\cdot Y(\tau,\xi)$, where
\bel{YA} Y(\tau,\xi)~\doteq~\int_\tau^T \kappa_1 A(t,\xi)\, dt + \kappa_2 A(T,\xi).\eeq
Here $A$ is the solution to the ODE (\ref{pat}) with initial data $A(\tau,\xi)=1$. 
The adjoint function $Y$ can be computed as the unique solution to the backward, linear ODE
\bel{Ydt}
\partial_t Y(t,\xi) ~=~ -\  \omega(t,\xi) \,Y(t,\xi)  - \kappa_1\,,\qquad\qquad  Y(T,\xi)~=~\kappa_2\,.
\eeq
Indeed, (\ref{pat}) and (\ref{Ydt}) together imply
$${d\over dt} (A Y)~=~-\kappa_1 AY,$$
$$Y(T,\xi) A(T,\xi)- Y(\tau,\xi) A(\tau,\xi)~=~-\int_\tau^T \kappa_1 A(t,\xi)\, dt,$$
and this yields (\ref{YA}) because $Y(T,\xi)=\kappa_2$, $A(\tau,\xi)=1$.

\begin{remark}\label{r:21} {\rm 
The adjoint variable $Y>0$ introduced at 
(\ref{Ydt}) can be interpreted as a ``shadow price".
Namely (see Fig.~\ref{f:sm42}),  think of $\Omega(t)$ as the contaminated set, 
and assume that  at time $\tau$ an external contractor offered to ``clean up"
a neighborhood of the point $x(\tau, \xi)$, thus replacing the set $\Omega(\tau)$ with a smaller set 
$\Omega^\ve(\tau)$, 
 at a price of $Y(\tau,\xi)$ per unit area.   
In this case, accepting or refusing the offer would make no difference in the total cost.   
}
\end{remark}

To obtain an optimality condition, consider  two boundary points $P_i=x(\tau,\xi_i)\in \partial
\Omega(\tau)$, $i=1,2$.   Assume that the control is active in a neighborhood of both points:
$\beta(\tau,\xi_i)>-1$.  If $Y(\tau, \xi_1)>Y(\tau,\xi_2)$ we can then rule out optimality.
Indeed, we can
increase the effort (i.e., the value of $\beta$) in a neighborhood of $P_1$ and decrease it by the same amount
in a in a neighborhood of $P_2$.  This will produce an admissible perturbation with a strictly lower total cost.

We conclude that, at each time $\tau$, the control $\beta(\tau,\cdot)$ can be active only at points where
$Y(\tau,\cdot)$ attains its global maximum.  Otherwise stated, one has the implication
\bel{max1}\beta(\tau,\xi)\,>\,-1\qquad\qquad \implies\qquad\qquad Y(\tau,\xi)~=~\max_\zeta \,Y(\tau,\zeta)
~\doteq~Y^*(\tau).\eeq

We are now ready to state our main result, providing necessary conditions for optimality.

\begin{theorem} \label{t:61} 
Assume that $t\mapsto \Omega(t)$ provides an optimal solution to {\bf (OP)}.   Let $(t,\xi)\mapsto x(t,\xi)$ be a  parameterization of the boundaries
of the sets $\Omega(t)$, satisfying the assumptions  {\bf (A1)-(A2)}.
Let $Y=Y(t,\xi)$ be the adjoint function introduced at  (\ref{Ydt}).

Then, for every $(t,\xi)\in W$ 
the inward normal velocity $\beta=\beta(t,\xi)$  satisfies (\ref{max1}).
\end{theorem}

A proof will be worked out in the next two sections.  

\begin{remark}\label{r:62}  {\rm The condition (\ref{max1}) does not explicitly determine the values
$\beta(\tau,\xi)$. However, it implies that at any time $t$ the portion of the boundary $\partial \Omega(t)$ where the control is active
must be the union of arcs with constant curvature $\omega(t)$.

Indeed,  by the assumptions  {\bf (A2)-(A3)}, both $\beta(t,\xi)$ and the solution $Y=Y(t,\xi)$ of the ODE (\ref{Ydt}) are continuous functions of  $(t,\xi)$.   Consider the open subset where the control is active:
\bel{Wact} W^{active}~\doteq~\bigl\{ (t,\xi)\,;~~\beta(t,\xi)>-1\bigr\}.\eeq
Assume $(t,\xi_1),(t,\xi_2)\in W^{active}$ for all $t\in \,]a,b[\,$. Then 
$$Y(t,\xi_1)~=~Y(t,\xi_2)~=~\max_{\zeta} Y(t,\zeta).\qquad\qquad  \forall t\in \,]a,b[\,.$$
In view of (\ref{Ydt}), this implies 
\bel{sameom} \omega(t,\xi_1)~=~\omega(t,\xi_2)\qquad\quad\hbox{for a.e.}~~ t\in \,]a,b[\,.\eeq
Since $\xi_1,\xi_2$ are arbitrary, we conclude that 
$\omega(t,\xi)=\omega(t)$ is constant for all $\xi$ such that $(t,,\xi)\in W^{active}$.
In particular: if $\omega(t)\not=0$, this means that at time $t$ the control effort is concentrated on the 
union of arcs of circumferences with radius $r(t) = \bigl|\omega(t)\bigr|^{-1}$.    If $\omega(t)=0$, at time $t$ 
the effort is concentrated on
a family of  straight lines.

In view of (\ref{Ydt}), one expects that the effort should be concentrated precisely on the portion of the boundary $\partial \Omega(t)$ where the curvature is maximum.  
The recent paper \cite{BBC} proves that this is indeed the case when $V=\R^2$ (i.e., without geographical constraints),
assuming that  the initial set $\Omega_0$ is   convex.}
\end{remark}

\section{A family of admissible perturbations}
\label{s:7}
\begin{figure}[ht]
\centerline{\hbox{\includegraphics[width=9cm]{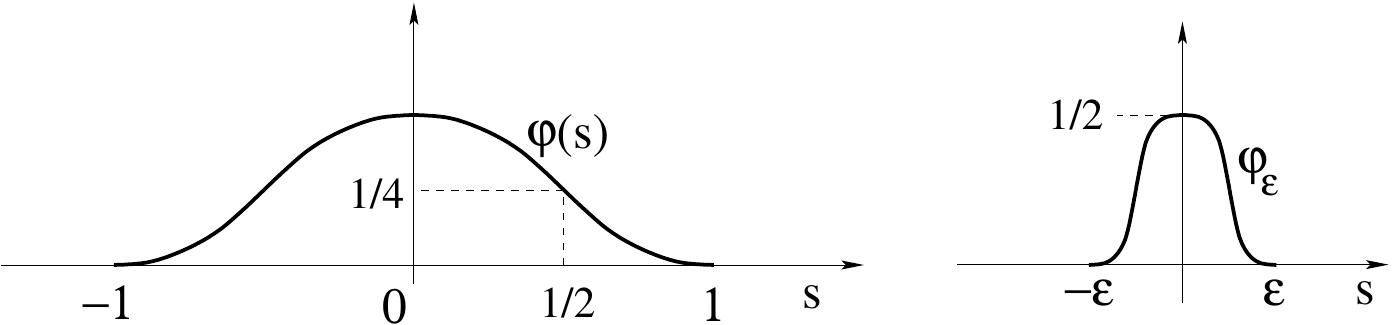}}}
\caption{\small The functions $\vp$ and $\vp_\ve$ introduced at (\ref{presc}). }
\label{f:sm76}
\end{figure}

\setcounter{equation}{0}
Toward a proof of Theorem~\ref{t:61} we will construct a family of admissible perturbations
of a set motion $t\mapsto\Omega(t)$, satisfying the regularity conditions {\bf (A1)-(A2)}.
\v
{\bf 1.}
Define the functions (see Fig.~\ref{f:sm76})
\bel{presc}\vp(s)~\doteq~\left\{ \bega{cl}{1\over 2} - s^2 &\hbox{for}~~|x|\leq {1\over 2},\\[2mm]
\bigl(1-|s|\bigr)^2&\hbox{for}~~ {1\over 2}\leq |s|\leq 1,\\[2mm]
0&\hbox{for}~~ |s|\geq 1,\enda\right.\qquad\qquad \vp_\ve(s)\doteq \vp\left( s\over\ve\right).\eeq
Notice that these are $\C^{1,1}$ functions, continuously differentiable with Lipschitz derivative and with compact support.
Moreover, they
satisfy
\bel{vpr1}
\int \vp(s)\, ds~=~{1\over 2}\,,\qquad\qquad \int \vp_\ve(s)\, ds~=~{\ve\over 2}\,,\eeq
\bel{vprop}\bigl[\vp'(s)\bigr]^2~\leq~4 \vp(s),\qquad \bigl[\vp_\ve'(s)\bigr]^2~\leq~4\ve^{-2} \vp_\ve(s)
\qquad\qquad\forall s\in \R.\eeq

Assume that, at a time $\tau\in \,]0,T[\,$, we can perturb the boundary of the set $\Omega(\tau)$
in a neighborhood of a point $x(\tau, \xi_0)$, so that  the new boundary is described by
\bel{pertau}
x^\ve(\tau,\xi)~=~x(\tau,\xi) \pm \ve^\gamma  \vp_\ve (\xi-\xi_0) \bfn(\tau,\xi) ,\eeq
for some exponent $\gamma>0$ whose value will be determined later.
As a consequence, the boundary $\partial \Omega(t)$ will be perturbed for all $t\in [\tau,T]$.
We seek a function $\sigma: [\tau,T]
\mapsto\R_+$, with 
\bel{Xtau}  \sigma(\tau)\,=\, \pm 1,  \eeq
and 
such that the perturbed motion described by
\bel{pert}
x^\ve(t,\xi)~=~x(t,\xi) + \ve^\gamma \sigma(t) \vp_\ve(\xi-\xi_0) \bfn(t,\xi) \eeq 
is admissible. 

Notice that, since the normal vector $\bfn$ is inward pointing, the perturbed sets
$\Omega^\ve(t)$ will be smaller in the case $\sigma(t)>0$, and larger when $\sigma(t)<0$.

To estimate the total effort $\E^\ve(t)$ for the perturbed motion, we write
$$\bfn^\ve(t,\xi)~\doteq ~{\bigl[x^\ve_\xi(t,\xi)\bigr]^\perp\over\bigl| x^\ve_\xi(t,\xi)\bigr|},
\qquad\qquad\beta^\ve(t,\xi)~\doteq~\la \bfn^\ve(t,\xi), \, x^\ve_t(t,\xi)\ra.$$
Since $x^\ve(t,\xi) =x(t,\xi)$ for $|\xi-\xi_0|\geq \ve$, it suffices to find a function $\sigma$ such that 
\bel{Eep} 
E\bigl(\beta^\ve(t,\xi)\bigr) \bigl|x^\ve_\xi(t,\xi)\bigr|~\leq ~E\bigl(\beta(t,\xi)\bigr)
\bigl|x_\xi(t,\xi)\bigr|\eeq
for all $(t,\xi)\in [\tau,T]\times [\xi_0-\ve,\,\xi_0+\ve]$.
\v
{\bf 2.} Observing that $\langle x_\xi, \bfn\rangle=\langle \bfn_\xi, \bfn\rangle=0$, recalling the formula
(\ref{curv}) for the curvature, one obtains
\bel{cu2} \la \bfn_\xi(t,\xi),\, x_{\xi}(t,\xi)\ra ~=~-
\la \bfn(t,\xi),\, x_{\xi\xi}(t,\xi)\ra ~=~-\omega(t,\xi) \bigl|x_\xi(t,\xi)\bigr|^2.\eeq
Differentiating (\ref{pert}) w.r.t.~$\xi$ we compute
\bel{xxi1}
x_\xi^\ve~=~x_\xi + \ve^{\gamma} \sigma\vp'_\ve \bfn + \ve^\gamma \sigma \vp_\ve\bfn_\xi\,,\eeq
\bel{xxi2}\bega{rl}
|x_\xi^\ve|^2
&=~
|x_\xi|^2 + 2 \la x_\xi ,\, \ve^\gamma \sigma\vp_\ve \bfn_\xi \ra+ \ve^{2\gamma} \sigma^2 (\vp'_\ve)^2
+  \ve^{2\gamma} \sigma^2\vp_\ve^2| \bfn_\xi|^2 \\[3mm]
&=~
\bigl(1-2\omega \ve^\gamma \sigma \vp_\ve\bigr) \,|x_\xi|^2 + \ve^{2\gamma} \sigma^2 (\vp'_\ve)^2
+  \ve^{2\gamma} \sigma^2\vp_\ve^2| \bfn_\xi|^2.
\enda\eeq
Taking square roots, since we are assuming that $|x_\xi|$ remains uniformly positive, we obtain
\bel{xxe}
|x^\ve_\xi| ~=~\bigl(1-\omega \ve^\gamma \sigma \vp_\ve\bigr) \,|x_\xi| 
+ \O(1)\cdot  \ve^{2\gamma}\sigma^2\bigl[(\vp'_\ve)^2 + \vp_\ve^2\bigr].\eeq
\v
{\bf 3.} Next, we need to estimate the new inward speed $\beta^\ve= \langle x^\ve_t, \bfn^\ve\rangle$.
For this purpose we compute
\bel{xet}
x^\ve_t~=~x_t +\ve^\gamma \dot \sigma \vp_\ve \bfn + \ve^\gamma \sigma \vp_\ve \bfn_t\,.\eeq
\bel{bfnt}
\bfn^\ve~=~{(x_\xi^\ve)^\perp\over |x^\ve_\xi|}\,.\eeq
It is convenient to write
\bel{xne} \bega{rl} \beta^\ve &=~\langle x^\ve_t, \bfn^\ve\rangle~=~ \langle x^\ve_t, \bfn\rangle +  \langle x^\ve_t, \bfn^\ve-\bfn\rangle
\\[2mm]
&=~\beta + \ve^\gamma \dot\sigma \vp_\ve + \langle x^\ve_t, \bfn^\ve-\bfn\rangle.\enda\eeq
Since $|\bfn^\ve|=|\bfn|=1$, we have
\bel{nne}  |\bfn^\ve-\bfn|~=~\O(1)\cdot \la \bfn,\, (x^\ve_\xi-x_\xi)^\perp\ra~=~\O(1)\cdot  \la x_\xi,\, (x^\ve_\xi-x_\xi)
\ra~=~\O(1)\cdot \la x_\xi,\, \ve^\gamma \sigma\vp_\ve \bfn_\xi
\ra.
\eeq
Moreover, since $x_t$ is parallel to $\bfn$,
assuming that $|x_\xi|$ remains uniformly positive, 
the last term on the right hand side of (\ref{xne}) can be bounded as
\bel{xn2}\bega{rl}
 \langle x^\ve_t, \bfn^\ve-\bfn\rangle&=~\O(1)\cdot\bigl( |x_t| +\ve^\gamma | \dot\sigma|\vp_\ve \bigr)\cdot 
|\bfn^\ve-\bfn|^2  + \O(1)\cdot \ve^\gamma \sigma \vp_\ve \, |\bfn_t|\, |\bfn^\ve-\bfn|\\[3mm]
&=~\O(1)\cdot \ve^{2\gamma} \vp_\ve^2\,.
\enda
\eeq
\v
{\bf 4.} 
Combining (\ref{xxe}), (\ref{bfnt}) and (\ref{nne}), and then using (\ref{vprop}), we obtain 
\bel{eff1}\bega{l} 
(1+\beta^\ve) |x_\xi^\ve|- (1+\beta) |x_\xi| \\[3mm]
\qquad = ~ \Big( 1+\beta+ \ve^\gamma \dot\sigma \vp_\ve+
 \O(1)\cdot \ve^{2\gamma}\sigma^2\vp_\ve^2  \Big)
 \cdot \Big( \bigl(1-\omega \ve^\gamma \sigma\vp_\ve\bigr) \,|x_\xi| 
+ \O(1)\cdot  \ve^{2\gamma}\sigma^2 \bigl[\vp_\ve^2 + (\vp_\ve')^2\bigr]\Big)\\[3mm]
\qquad\qquad - (1+\beta) |x_\xi|\\[3mm]
\qquad =~ \ve^\gamma\bigl( \dot\sigma - (1+\beta) \omega \sigma\bigr)\vp_\ve |x_\xi| + 
\O(1)\cdot \ve^{2\gamma}\sigma^2 \vp_\ve^2 + \O(1)\cdot (1+\beta+\ve^\gamma \dot\sigma\vp_\ve) \,  \ve^{2\gamma-2}
\vp_\ve\\[3mm]
\qquad =~ \ve^\gamma\bigl( \dot\sigma - (1+\beta) \omega \sigma\bigr)\vp_\ve |x_\xi| +\O(1)\cdot  \ve^{2\gamma-2}
\vp_\ve \,. \enda
\eeq
\v
{\bf 5.} We now choose
\bel{agchoose}
\alpha\,=\, 1,\qquad\gamma\,=\,4,\eeq
\bel{dschoose}
\dot \sigma(t)~\doteq ~\min_{|\zeta-\xi_0|\leq\ve} \bigl(1+\beta(t,\zeta)\bigr) \omega(t,\zeta) \sigma(t) - \ve^\alpha. \eeq
We claim that, with the above choices, the left hand side of  (\ref{eff1}) 
is non-positive, for all $t,\xi\in [\tau,T]\times [\xi_0-\ve, \xi_0+\ve]$ and all $\ve>0$ sufficiently small.
Indeed, 
from (\ref{eff1}),  it now follows
\bel{eff2}
\bega{l} 
\bigl(1+\beta^\ve(t,\xi)\bigr) \bigl|x_\xi^\ve(t,\xi)\bigr|- \bigl(1+\beta(t,\xi)\bigr) \bigl|x_\xi(t,\xi)\bigr| \\[3mm]
\qquad \leq~ -\ve^{\gamma+\alpha} \vp_\ve(\xi-\xi_0)  \bigl|x_\xi(t,\xi)\bigr| +\O(1)\cdot  \ve^{2\gamma-2}
\vp_\ve (\xi-\xi_0)\\[3mm]
\qquad = ~\Big( - \ve^5  \bigl|x_\xi(t,\xi)\bigr|  + \O(1) \cdot \ve^6\Big) \vp_\ve (\xi-\xi_0)~\leq ~0\,, \enda
\eeq
since we are assuming that $|x_\xi|$ remains uniformly positive.

\section{Proof of Theorem~\ref{t:61}}
\label{s:8}
\setcounter{equation}{0}

In this section we give a proof of Theorem~\ref{t:61}, in several steps.
\v
{\bf 1.} 
Let $t\mapsto\Omega(t)$ be an optimal strategy for the problem {\bf (OP)}, and 
let $Y=Y(t,\xi)$ be the adjoint function, 
obtained by solving  the backward Cauchy problem (\ref{Ydt}).
 Assume that at some time $\tau\in [0,T]$ the condition
(\ref{max1}) is violated. Namely, there exist $(\tau, \xi_1), (\tau, \xi_2)\in W$ such that
\bel{bY}
\beta(\tau,\xi_1)\,>\,-1,\qquad\qquad Y(\tau, \xi_1) \,<\, Y(\tau, \xi_2).\eeq
By continuity, without loss of generality we can assume that $0<\tau<T$.
In this case, by reducing the control effort in a neighborhood of $\xi_1$ (where it is less effective)
and increasing the effort near $\xi_2$ (where it is more effective),
we will construct another admissible strategy $t\mapsto \Omega^\ve(t)$ 
with strictly smaller total cost, thus achieving a contradiction.
\v
{\bf 2.} For $\ve>0$ small we let $\Omega^\ve(t)=\Omega(t)$ for $t\in [0, \tau-\ve]$.
Calling $x^\ve(t,\xi)$ the parameterization of the boundary $\partial \Omega^\ve(t)$.
On the interval $t\in [\tau-\ve,\tau]$,  we define
\bel{Oe1} x^\ve(t,\xi)~=~x(t,\xi) - \ve^\gamma
\sigma_1(t)  \bfn(t,\xi)+ \ve^\gamma\sigma_2(t)\bfn(t,\xi),
\eeq
for suitable functions $\sigma_1, \sigma_2$.  These must be chosen so that the perturbed motion is
also admissible.
Recalling (\ref{eff1}), for every $t\in [\tau-\ve, \tau]$ this requires
\bel{eff3}\bega{l} 0~\geq~\ds\sum_{i=1,2}\int_{|\xi-\xi_i|<\ve} \Big\{
(1+\beta^\ve) |x_\xi^\ve|- (1+\beta) |x_\xi|\Big\}\,d\xi \\[4mm]
 \qquad \ds =~ \sum_{i=1,2}\int_{|\xi-\xi_i|<\ve} \Big\{ \ve^\gamma\bigl( \dot\sigma_i - (1+\beta) \omega \sigma\bigr)\vp_\ve |x_\xi| +\O(1)\cdot  \ve^{2\gamma-2}\, 
\vp_\ve \Big\}\,d\xi\\[4mm]
 \qquad \ds =~ \sum_{i=1,2}\int_{|\xi-\xi_i|<\ve} \Big\{ \ve^\gamma \dot\sigma_i \vp_\ve |x_\xi| +\O(1)\cdot  \ve^{\gamma}\, 
\vp_\ve \Big\}\,d\xi , \enda
\eeq
To achieve (\ref{eff3}), for $i=1,2$ we first define
\bel{etae} \eta_i(\ve)~\doteq~\sup\left\{ {\bigl|x_\xi(t,\xi)-x_\xi(\tau,\xi_i)\bigr|\over \bigl|x_\xi(\tau,\xi_i)\bigr|}\,;\quad
t\in [\tau-\ve,\tau],~~|\xi-\xi_i|<\ve\right\}.\eeq
Then we let $\sigma_1,\sigma_2:[\tau-\ve,\tau]\mapsto\R$ be the affine functions such that
\bel{dsi} 
\left\{\bega{rl}\dot\sigma_i(t)&\ds =\, {-1\over  \ve\,\bigl|x_\xi(\tau,\xi_1)\bigr|} -{\eta_1(\ve)\over\ve} -{1\over \sqrt\ve}\,,
\\[4mm]
\dot\sigma_2(t)&\ds =~ {1\over  \ve\,\bigl|x_\xi(\tau,\xi_2)\bigr|} -{\eta_2(\ve)\over\ve} -{1\over \sqrt\ve}\,,\enda\right.
\qquad\qquad \sigma_1(\tau-\ve)\,=\,\sigma_2(\tau-\ve)\,=\,0.
\eeq
Inserting (\ref{dsi}) in the right hand side of (\ref{eff3}) and observing that the leading order terms cancel,  
for every $\ve>0$ sufficiently small we obtain
\bel{eff4}\bega{l} \ds\sum_{i=1,2}\int_{|\xi-\xi_i|<\ve} \bigg\{ \ve^\gamma \dot\sigma_i(t) 
\Big[\bigl|x_\xi(\tau,\xi_i)\bigr| +\Big(  \bigl|x_\xi(t,\xi)\bigr| -\bigl|x_\xi(\tau,\xi_i)\bigr| \Big) \Big]
+ \O(1)\cdot  \ve^{\gamma}\bigg\}
\vp_\ve(\xi-\xi_i) \,d\xi \\[4mm]
\ds\quad\leq~\sum_{i=1,2}\int_{|\xi-\xi_i|<\ve} \bigg\{ -\ve^{\gamma-{1\over 2} }
+ \O(1)\cdot  \ve^{\gamma}\bigg\}
\vp_\ve(\xi-\xi_i) \,d\xi ~<~0, \enda
\eeq
showing that the perturbed motion is admissible.
\v
{\bf 3.} On the remaining interval $[\tau,T]$, recalling (\ref{dschoose}) we let the functions $\sigma_i$
be the solutions to the linear ODEs
\bel{dsii}
\dot \sigma_i(t)~=~\min_{|\zeta-\xi_i|\leq\ve}~\bigr(1+\beta(t,\zeta)\bigr) \omega(t,\zeta)\sigma_i(t)-\ve^\alpha,
\eeq
with initial data given at $t=\tau$ corresponding to the function constructed in step {\bf 2}, namely
\bel{idsi} 
\sigma_1(\tau)\,=\,\bar\sigma_1\,\doteq\, {-1\over \bigl|x_\xi(\tau,\xi_1)\bigr|} -{\eta_1(\ve)} - \sqrt\ve\,,
\qquad\qquad \sigma_2(\tau)\,=\,\bar\sigma_2\,\doteq\, {1\over \bigl|x_\xi(\tau,\xi_2)\bigr|} -{\eta_2(\ve)} - \sqrt\ve\,.\eeq
As shown by the analysis in Section~\ref{s:7}, by (\ref{eff2}) this motion is admissible.
\v
{\bf 4.} It remains to prove that, for $\ve>0$, the new strategy achieves a strictly lower cost.
Let the functions $\sigma_i^\sharp:[\tau,T]\mapsto \R$, $i=1,2$, be the solutions to
\bel{sshi}
\dot \sigma_i^\sharp(t)~=~\bigl( 1+\beta(t,\xi_1)\bigr) \sigma_1^\sharp(t),\qquad\qquad \sigma_i^\sharp(\tau)~=~1.\eeq
Thanks to the assumption (\ref{lepo}), our construction implies the uniform convergence
\bel{scon}
\sigma_{1,\ve}(t)~\to~-{\sigma^\sharp(t)\over \bigl|x_\xi(\tau,\xi_1)\bigr|}\,,\qquad\qquad \sigma_{2,\ve}(t)~\to~{\sigma_2^\sharp(t)\over \bigl|x_\xi(\tau,\xi_2)\bigr|},\eeq
uniformly for $t\in [\tau,T]$.

Recalling (\ref{vpr1}) and then using (\ref{YA})-(\ref{Ydt}), the difference in the total cost over the interval 
$[\tau,T]$ can be estimated as
\bel{cdiff}\bega{l} \ds \kappa_1 \int_\tau^T \Big[\caL^2\bigl(\Omega^\ve(t)\bigr)- \caL^2\bigl(\Omega(t)\bigr)\Big]\, dt+ 
\kappa_2 \Big[\caL^2\bigl(\Omega^\ve(T)\bigr)- \caL^2\bigl(\Omega(T)\bigr)\Big]\\[4mm]
\qquad \ds =~\kappa_1
\int_\tau^T \ve^\gamma \cdot {\ve\over 2} \big[ \sigma_1^\sharp(t) - \sigma_2^\sharp(t)\bigr]\, dt + 
\kappa_2 \ve^\gamma \cdot {\ve\over 2} \big[ \sigma_1^\sharp(T) - \sigma_2^\sharp(T)\bigr] + o(\ve^{\gamma+1})
\\[4mm]
\qquad \ds =~{\ve^{\gamma+1}\over 2}\Big( Y(\tau,\sigma_1)- Y(\tau,\sigma_2)\Big)+  o(\ve^{\gamma+1})~<~0,
\enda
\eeq
As usual, here $o(\ve^{\gamma+1})$ denotes an infinitesimal of higher order, as $\ve\to 0$.

On the other hand, during  the time interval $[\tau-\ve,\tau]$ the difference in cost can be bounded as
\bel{cdtau}
 \int_{\tau-\ve}^\tau \Big[\caL^2\bigl(\Omega^\ve(t)\bigr)- \caL^2\bigl(\Omega(t)\bigr)\Big]\, dt~=~\O(1)\cdot \ve^{\gamma+2}.\eeq
Combining (\ref{cdiff})-(\ref{cdtau}) we conclude
$$J(\Omega^\ve)-J(\Omega)~<~0$$
for all $\ve>0$ sufficiently small. 
This contradicts optimality, proving the theorem.
\endproof

\section{Optimality conditions for the minimum time problem}
\label{s:9}
\setcounter{equation}{0}

\begin{figure}[ht]
\centerline{\hbox{\includegraphics[width=9cm]{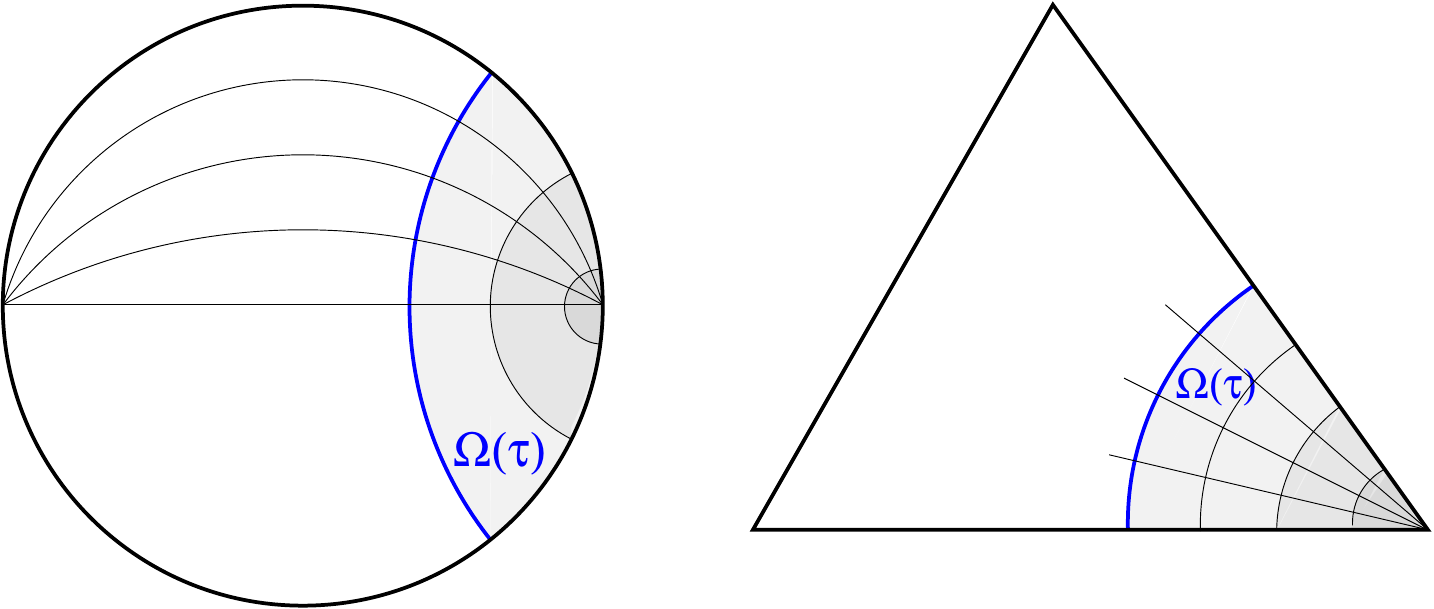}}}
\caption{\small As $t\to T-$ and the eradication is nearly completed, the length of boundary  $\partial \Omega(t)$
approaches zero.   Hence the assumption that $\bigl|x_\xi(t,\xi)\bigr|$ remains uniformly positive as $t\to T-$  is never satisfied.   On the other hand, in most cases the optimal set motion satisfies the necessary and sufficient conditions 
in Corollary~\ref{c:41}, on some terminal interval $[\tau,T]$.
}
\label{f:sm78}
\end{figure}

If one tries to apply the optimality conditions derived in the previous sections to the minimum time eradication 
problem, the construction invariably breaks down.
Indeed, a key assumption used in the construction of admissible perturbations is that 
the arc-length $\bigl|x_\xi(t,\xi)\bigr|$ remains uniformly positive for $t\in [\tau,T]$.  However, in all examples
(see Fig.~\ref{f:sm78})
one finds $\bigl|x_\xi(t,\xi)\bigr|\to 0$ for all $\xi$, as $t\to T-$.
To obtain necessary conditions,  a different approach is needed.  

We start with the simple observation that, if $t\mapsto \Omega(t)$, $t\in [0,T]$, is an optimal solution to the minimum time eradication problem, then for any intermediate time $0<T_1<T$ one has:
\begi 
\item[{\bf (1)}] Restricted to the subinterval $[0, T_1]$, the strategy $t\mapsto \Omega(t)$  is optimal for the minimum time transfer
problem {\bf (MTTP)}, with
\bel{minit1}\Omega_{initial}= V,\qquad\qquad \Omega_{final} = \Omega(T_1),\eeq
\item[{\bf (2)}]  Restricted to the subinterval $[ T_1,T]$, the strategy $t\mapsto \Omega(t)$   is optimal for the minimum time transfer
problem {\bf (MTTP)}, with
\bel{minit2}\Omega_{initial}= \Omega(T_1),\qquad\qquad \Omega_{final} = \emptyset.\eeq
\endi

In many situations, one can choose $T_1<T$ so that the optimal transfer problem 
with data (\ref{minit2}) has a solution satisfying the assumptions of Corollary~\ref{c:41}. This provides 
necessary and sufficient conditions for optimality, restricted to the subinterval $[T_1, T]$.

It thus remains to derive necessary conditions for the problem (\ref{minit1}).  
This will require an additional assumption on the parameterization $(t,\xi)\mapsto x(t,\xi)$ of the boundaries 
$\partial \Omega(t)$, for $t\in[0,T_1]$.
\begi
\item[{\bf (A4)}] {\it
At time $T_1$ the effort is uniformly positive along the boundary $\partial \Omega(T_1)\cap V$.
Namely, there  exists $\delta_0>0$ such that 
\bel{poseff} 
\beta(T_1,\xi)~\geq~-1+\delta_0\qquad\qquad\forall  \xi\,.\eeq
}
\endi
The next result shows that, if $t\mapsto\Omega(t)$ is optimal for the minimum time transfer problem
(\ref{minit2}) and the additional assumptions {\bf (A4)} hold, then it satisfies the necessary conditions 
for the problem of minimizing the terminal area at time $t=T_1$.
\bel{minarea} \hbox{Minimize:}~~\caL^2\bigr(\Omega(T_1)\bigr)\qquad\hbox{subject to:}~~\Omega(0)=V,\eeq
among all admissible motions.
Notice that this is precisely the problem {\bf (OSM)}, in the special case where $\kappa_1=0$, $\kappa_2=1$ in (\ref{FU}).  The corresponding adjoint equations (\ref{Ydt}) reduce to
\bel{Yeq} 
\partial_t Y(t,\xi)~=~-\omega(t,\xi) Y(t,\xi),\qquad\qquad Y(T_1)\,=\,1.\eeq

\begin{theorem} \label{t:91} 
Assume that $t\mapsto \Omega(t)$ provides an optimal solution to the minimum time transfer problem {\bf (MTTP)}
with data (\ref{minit1}).   Let $(t,\xi)\mapsto x(t,\xi)$ be a  parameterization of the boundaries
$\partial\Omega(t)$, satisfying the regularity assumptions in  {\bf (A1)-(A3)} with $T=T_1$.
Moreover, assume that {\bf (A4)} holds.
Let $Y=Y(t,\xi)$ be the adjoint function at  (\ref{Yeq}).

Then, for every $(t,\xi)\in \,]0,T_1[\,\times\,] 0,1[\,$,
the inward normal velocity $\beta=\beta(t,\xi)$  satisfies the implication
\bel{bop}
\beta(t,\xi)~>~-1\qquad\implies\qquad Y(t,\xi)\,=\,\max_{\zeta\in \,]0,1[} Y(t,\zeta).\eeq
\end{theorem}

{\bf Proof.} 
{\bf 1.} Assume that at some time $\tau\in [0,T]$ the condition
(\ref{bop}) is violated. Namely, there exist $(\tau, \xi_1), (\tau, \xi_2)$ with $0<\tau<T_1$,  such that
$$
\beta(\tau,\xi_1)\,>\,-1,\qquad\qquad Y(\tau, \xi_1) \,<\, Y(\tau, \xi_2).$$
 For $\ve>0$ small we then let $\Omega^\ve(t)=\Omega(t)$ for $t\in [0, \tau-\ve]$.
Calling $x^\ve(t,\xi)$ the parameterization of the boundary $\partial \Omega^\ve(t)$.
On the interval $t\in [\tau-\ve,\tau]$,  we define
\bel{pertm} x^\ve(t,\xi)~=~x(t,\xi) - \ve^\gamma
\sigma_1(t)  \bfn(t,\xi)+ \ve^\gamma\sigma_2(t)\bfn(t,\xi),
\eeq
for suitable functions $\sigma_1, \sigma_2$, constructed as in steps {\bf 2} and {\bf 3} of the proof of Theorem~\ref{t:61}.
This guarantees that the perturbed strategy (\ref{pertm}) is admissible.

\begin{figure}[ht]
\centerline{\hbox{\includegraphics[width=13cm]{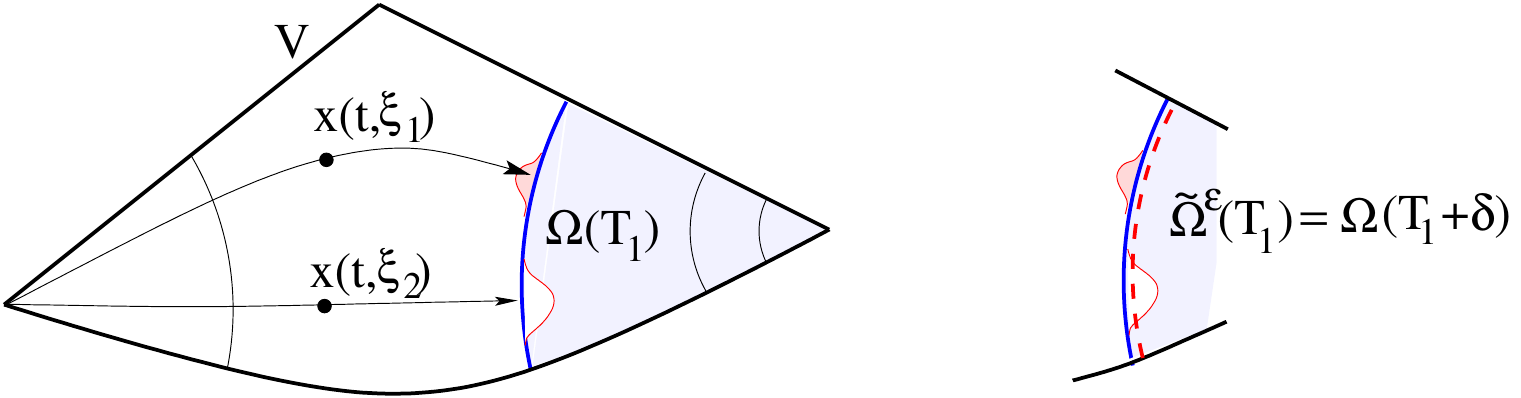}}}
\caption{\small Left: The perturbed set $\Omega^\ve(T_1)$, whose boundary is parameterized as
in (\ref{pertm}),  is obtained from $\Omega(T_1)$ by adding a small region near $x(T_1, \xi_1)$ and 
removing a larger region near $x(T_1, \xi_2)$.   Right: by a further modification, we obtain a 
family of sets such that $\Tilde \Omega^\ve(t) = \Omega(t')$, with $T_1\leq t<t'$.
This rules out optimality for the minimum time problem.
}
\label{f:sm79}
\end{figure}

The same analysis as in (\ref{cdiff}) now yields
\bel{ca}\bega{l} \ds\caL^2\bigl(\Omega^\ve(T_1)\bigr)- \caL^2\bigl(\Omega(T_1)\bigr)
~=~{\ve^{\gamma+1}\over 2}\Big( Y(\tau,\sigma_1)- Y(\tau,\sigma_2)\Big)+  o(\ve^{\gamma+1})~<~0.
\enda
\eeq
This shows that $\Omega(\cdot)$ is not optimal for problem of minimizing the terminal area at time $T_1$.
In particular, the previous construction yields
\bel{area2} \sigma_2(T_1) \,\bigl| x_\xi(T_1, \xi_2)\bigr| -\sigma_1(T_1) \,\bigl| x_\xi(T_1, \xi_1)\bigr| ~\geq
c_1~>~0,\eeq
for some constant $c_1>0$ independent of $\ve>0$.  
As shown in Fig.~\ref{f:sm79},  this means that
the perturbed set  $\Omega^\ve(T_1)$ contains an additional small region
in a neighborhood of the point $x(T_1,\xi_1)$, while a larger region is removed in a neighborhood of $x(T_1,\xi_2)$.

\v
{\bf 2.} By a further modification, we now construct a set motion $t\mapsto \Hat\Omega^\ve(t)$ such that
\bel{HOME} 
\Hat\Omega^\ve(T_1-\ve) \,=\,\Omega^\ve(T_1-\ve), \qquad\qquad \Hat\Omega^\ve(T_1)\,=\,\Omega(T_1),\eeq
while, for all $t\in [T_1-\ve, T_1]$, a strictly smaller effort is required:
\bel{HET} \Hat \E^\ve(t)~=~\int_0^1 \big(1+\Hat \beta^\ve(t,\xi)\big)\,\bigl|\Hat x^\ve_{\xi}(t,\xi)\bigr|\, d\xi~\leq~\E^\ve(t)-\delta_\ve~\leq~M-\delta_\ve\,,\eeq
for some $\delta_\ve>0$.

This new motion is obtained replacing (\ref{pertm}) with
\bel{hpert}  \Hat x^\ve(t,\xi)~=~x(t,\xi) - \ve^\gamma
\Hat \sigma_1(t)  \bfn(t,\xi)+ \ve^\gamma\Hat \sigma_2(t)\bfn(t,\xi),
\eeq
where
\bel{Hsi} \Hat\sigma_i(t)~=~
\left\{\bega{cl}  \sigma_i(t)\quad &\hbox{if}\quad t\leq T_1-\ve,\\[2mm]
{T_1-t\over\ve} \, \sigma_i(T_1-\ve)\quad &\hbox{if}\quad t\in [T_1-\ve, T_1].
\enda\right.\eeq
This clearly implies (\ref{HOME}). 

Concerning the total effort, by the assumption {\bf (A4)} for $\ve>0$ small by continuity we can assume
$\Hat \beta^\ve>-1$,  hence $E(\Hat \beta^\ve) = 1+\Hat\beta^\ve$.
Moreover, for $t\in [T_1-\ve, T_1]$, similar computations as in (\ref{cdiff})
yield
\bel{difef}\bega{rl} \Hat \E^\ve(t)-\E^\ve(t)&\ds =~\ve^\gamma \int \Big(\dot \sigma_2(t) \bigl|x_\xi(t,\xi_2)\bigr| \vp_\ve(\xi-\xi_2) -
 \dot \sigma_1(t) \bigl|x_\xi(t,\xi_1)\bigr| \vp_\ve(\xi-\xi_1) \Big)d\xi + o(\ve^\gamma)\\[4mm]
 &\ds=~{\ve^\gamma \over 2} \Big( -\sigma_2(T_1-\ve)  \bigl|x_\xi(t,\xi_2)\bigr|+
 \sigma_1(T_1-\ve)  \bigl|x_\xi(t,\xi_1)\bigr|
 \Big) + o(\ve^\gamma)\\[4mm]
 &\ds=~{\ve^\gamma \over 2} \Big( -\sigma_2(T_1)  \bigl|x_\xi(T_1,\xi_2)\bigr|+
 \sigma_1(T_1)  \bigl|x_\xi(T_1,\xi_1)\bigr|
 \Big) + o(\ve^\gamma)\\[4mm]
 &\ds \leq ~- {c_1 \ve^\gamma \over 2} + o(\ve^\gamma)~<~0,
 \enda
\eeq
where we used continuity and finally (\ref{area2}). 
\v
{\bf 3.} In the previous step we constructed a set-motion that reaches the same terminal 
configuration $\Omega(T_1)$ but with a strictly smaller effort (\ref{HET}). 
In this step we use the additional available effort $M-\Hat \E^\ve(t)>0$ to construct a further motion
$t\mapsto \Tilde \Omega^\ve(t)$ that reaches $\Omega(T_1)$ in a strictly shorter time.
This  motion will have the form
\bel{TOT}
t~\mapsto ~\Tilde \Omega^\ve(t)~=~\Hat \Omega^\ve\bigl(\tau(t)\bigr)\qquad \qquad t\in [T_1-\ve, \, T_1],\eeq
with $\tau(t)>t$ for $t>T_1-\ve$.   

For this purpose, notice that the corresponding inner normal velocity for the motion (\ref{TOT}) is
$$\Tilde \beta^\ve(t,\xi)~=~\dot \tau(t) \cdot \Hat \beta^\ve\bigl(\tau(t),\xi\bigr).$$ 
By the assumption (\ref{poseff}) we can assume that $E(\beta^\ve), E(\Hat \beta^\ve)$ remain uniformly positive
for $T_1-\delta_1<t<T$.  The instantaneous total efforts are thus computed by
\bel{TTE}\bega{rl} \Tilde \E^\ve(t)&=~\ds
\int_0^1 \Big( 1+ \Tilde \beta^\ve(t,\xi)\Big) \bigl|\Tilde x^\ve_\xi (t,\xi)\bigr|\, d\xi~=~
\int_0^1 \Big( 1+\dot \tau(t) \cdot \Hat \beta^\ve\bigl(\tau(t),\xi\bigr)\Big) \bigl|\Hat x^\ve_\xi (t,\xi)\bigr|\, d\xi\\[4mm]
&=~\ds \Hat\E^\ve(t) + 
\int_0^1 \bigl( \dot \tau(t)-1\bigr) \cdot \Hat \beta^\ve\bigl(\tau(t),\xi\bigr)\, \bigl|\Hat x^\ve_\xi (t,\xi)\bigr|\, d\xi
\,.\enda
\eeq
In view of (\ref{HET}), the right hand side is $\leq M$ provided we choose
\bel{dott}
1~<~\dot\tau(t)~<~1+\delta_\ve \left[ \int_0^1\Hat \beta^\ve\bigl(\tau(t),\xi\bigr)\, \bigl|\Hat x^\ve_\xi (t,\xi)\bigr|\, d\xi
\right]^{-1}.\eeq
Let $T_1'<T_1$ be such that $\tau(T_1')=T_1$.    Then the motion $t\mapsto \Tilde \Omega^\ve(t)$ 
is admissible and satisfies $\Tilde\Omega^\ve(T_1') = \Omega(T_1)$, showing that $\Omega$ is not time optimal.
This contradiction proves the theorem.
\endproof

\begin{remark}
{\rm In the above proof, the assumption (\ref{poseff}) was only used to conclude that the effort 
$E\bigl(\beta(T_1,\xi)\bigr)$ is strictly positive in a neighborhood of the point $x(T_1, \xi_2)$.
A sharper version of the above result is as follows.
Let $t\mapsto \Omega(t)$ be optimal for the minimum time problem. Assume 
\bel{AA1} \beta(T_1,\xi_2)\,>\,-1  \eeq
and let $Y=Y(t,\xi)$ be as in (\ref{Yeq}).   
Then, at any point $(\tau,\xi_1)$ where $\beta(\tau,\xi_1)>-1$ one has
\bel{AA2} Y(\tau, \xi_1)~\geq~Y(\tau, \xi_2).\eeq
}\end{remark}

\section{Optimality conditions at junctions}
\label{s:10}
\setcounter{equation}{0}

In Theorem~\ref{t:41}, the existence of optimal solutions  was proved within a class of 
functions with BV regularity.  On the other hand,  the necessary conditions for optimality 
derived in Theorem~\ref{t:61} require that the sets $\Omega(t)$ have $\C^{1,1}$  boundary. 
Indeed, we used this assumption to uniquely define the perpendicular curves $t\mapsto x(t,\xi)$.

\begin{figure}[ht]
\centerline{\hbox{\includegraphics[width=9cm]{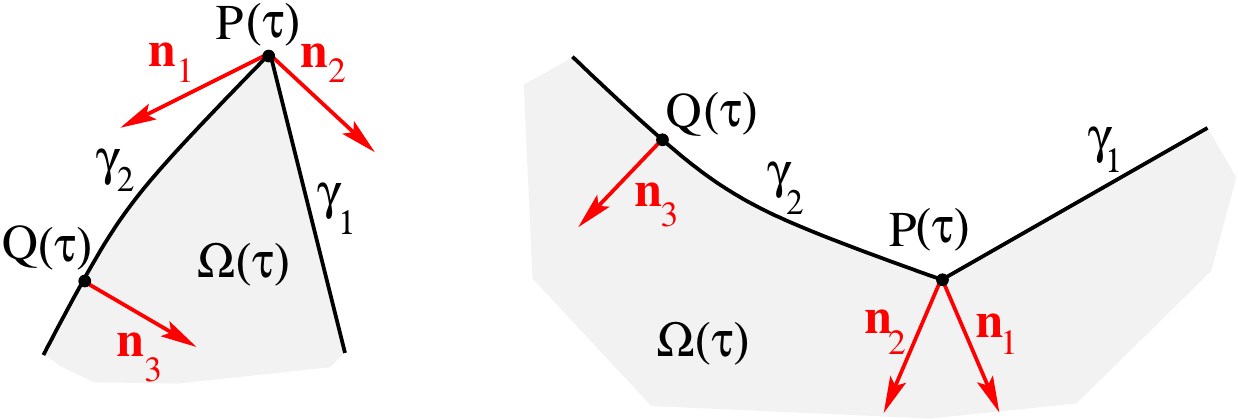}}}
\caption{\small Left: the case of outward corner. Right: the case of an inward corner.}
\label{f:sm81}
\end{figure}

Aim of this section is to partially fill this regularity gap, ruling out certain configurations where  the sets $\Omega(t)$ have 
corners.   The following situation will be considered:
\begi
\item[{\bf (A5)}] {\it
There exists $\tau,\delta_0>0$ such that, for $|t-\tau|<\delta_0$,  the boundary $\partial\Omega(t)$ contains two adjacent arcs $\gamma_1(t,\cdot)$, $\gamma_2(t,\cdot)$ joining at a point $P(t)$ 
at an angle $\theta(t)$.    Each of these arcs admits a $\C^1$  parameterization of the form
\bel{angle}\left\{ \bega{l}s\mapsto \gamma_1(t,s),\qquad s\leq 0,\\[2mm]
s\mapsto \gamma_2(t,s),\qquad s\geq 0,\enda
\right.\qquad\qquad \gamma_1(t,0) = \gamma_2(t,0)
=P(t).\eeq
}\endi
The tangent vectors to the curves $\gamma_1(\tau,\cdot)$ and $\gamma_2(\tau,\cdot)$ at the intersection point
$P(\tau)$ will be denoted by
\bel{bw}\bfw_1~=~\gamma_{1,s}(\tau, 0-)\,,\qquad\qquad \bfw_2~=~\gamma_{2,s}(\tau, 0+)\,.\eeq
Moreover, we call $\bfw^\perp_1,\bfw_2^\perp$ the orthogonal vectors (rotated by $90^o$ counterclockwise).
Notice that the two curves $\gamma_1,\gamma_2$ form an outward corner at $P(\tau)$ if the vector product
satisfies (see Fig.~\ref{f:sm81}, left)
$$\bfw_1\times \bfw_2~\doteq~\la \bfw_1^\perp, \bfw_2\ra~>~0.$$
On the other hand, if $\bfw_1\times \bfw_2<0$ one has an inward corner, as shown in Fig.~\ref{f:sm81}, right.

As before, we say that the control is {\em active} on a portion of the boundary $\partial \Omega(t)$ if
the inward normal speed is $\beta>-1$, hence the effort is strictly positive: $E(\beta)= 1+\beta>0$.
The main result of this section shows that, for an optimal
motion $t\mapsto \Omega(t)$, non-parallel junctions  cannot be optimal if the control is active on at least
one of the adjacent arcs.  

\begin{remark}  {\rm If the control is not active along any of the two arcs $\gamma_1, \gamma_2$,
then in a neighborhood of $P(t)$ the set $\Omega(t)$ expands with unit speed all along the boundary.
This implies that for $t>\tau$ the set $\Omega(t)$ satisfies an interior ball condition, 
hence it can only have
inward corners.
}
\end{remark}

\begin{figure}[ht]
\centerline{\hbox{\includegraphics[width=15cm]{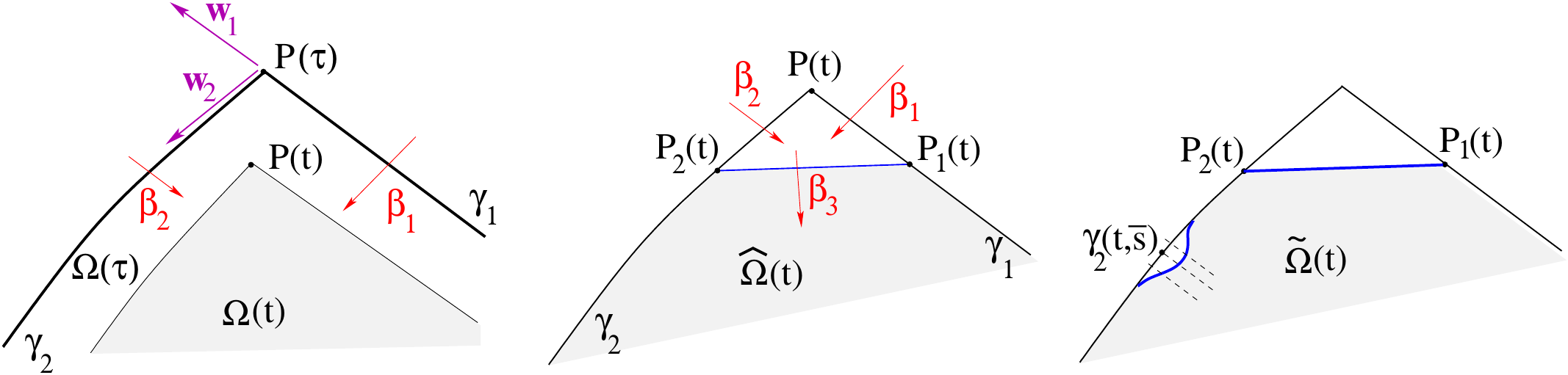}}}
\caption{\small The case of an outward corner, with $\bfw_1\times \bfw_2>0$. 
Left: the original set $\Omega(t)$, for $t\in [\tau, \tau+\delta]$. 
Center: the set $\Hat \Omega(t)$ obtained from $\Omega(t)$ by  removing the triangular
region $\Hat{P_1PP_2}$. Right: the set $\Tilde\Omega(t)$ obtained by removing an additional region
in a neighborhood of $\gamma_2(t,\bar s)$.}
\label{f:sm83}
\end{figure}

\begin{figure}[ht]
\centerline{\hbox{\includegraphics[width=16cm]{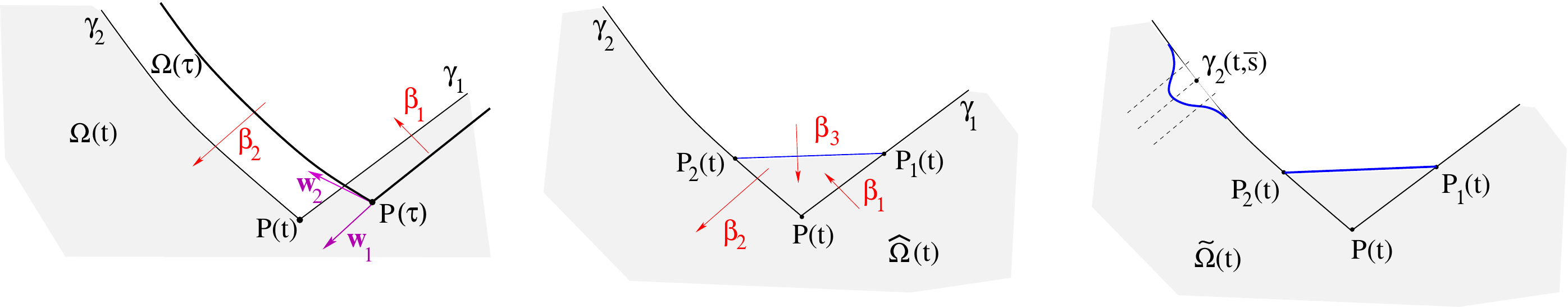}}}
\caption{\small The case of an inward corner, with $\bfw_1\times \bfw_2<0$. 
Left: the original set $\Omega(t)$, for $t\in [\tau, \tau+\delta]$. 
Center: the set $\Hat \Omega(t)$ obtained from $\Omega(t)$ by  adding the triangular
region $\Hat{PP_1P_2}$. Right: the set $\Tilde\Omega(t)$ obtained by removing an additional region
in a neighborhood of $\gamma_2(t,\bar s)$.}
\label{f:sm84}
\end{figure}

The next result shows that non-tangential junctions between two arcs cannot be optimal.
An intuitive explanation comes from the basic formula (\ref{darea2}).   
By ``cutting the corners", as shown in Figures~\ref{f:sm83} and \ref{f:sm84}, we can reduce the 
perimeter $\partial\Omega(t)$ in a neighborhood of the junction point $P(t)$.
After this perturbation,  the control effort spent in this neighborhood of $P(t)$ is reduced,
hence $\E(t)<M$.   We can then use the additional available effort $M-\E(t)$
to shrink the set $\Omega(t)$ in the neighborhood of some other point $Q(t)$ where the boundary
is smooth.  This will produce a bigger decrease in the area, thus decreasing the total 
cost $J(\Omega)$ at (\ref{J}).

%

\begin{theorem}\label{t:101} 
Assume that $t\mapsto \Omega(t)$ provides an optimal solution to {\bf (OP)}, with $\kappa_1>0$. 

In the setting described at {\bf (A5)}, if the control is active along at least one of the two arcs $\gamma_1, \gamma_2$,
 then these two arcs must be tangent at $P(t)$.
\end{theorem}

{\bf Proof.} {\bf 1.} Assume that the  two arcs $\gamma_1,\gamma_2$ are not tangent, forming an angle $\theta(t)\neq \pi$ for $|\tau-t|<\delta$. We 
will then construct an admissible strategy  $t\mapsto\Omega^\ve (t)$ with smaller total cost.

\begin{figure}[ht]
\centerline{\hbox{\includegraphics[width=9cm]{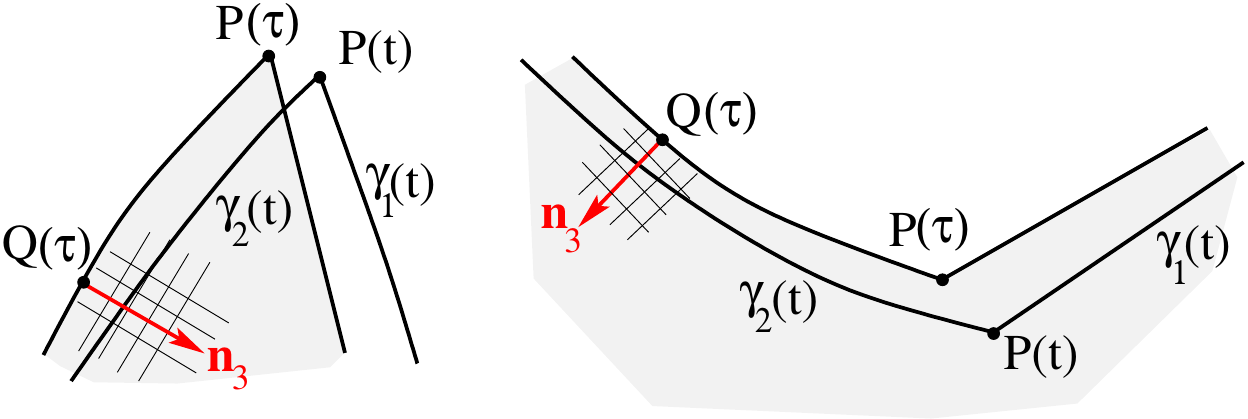}}}
\caption{\small In a neighborhood of the point $Q(\tau)$, the curves $\gamma_2(t,\cdot)$ are 
parameterized according to a set of coordinates with axes parallel to $\bfn_3, \bfn_3^\perp$.}
\label{f:sm82}
\end{figure}

To fix ideas, assume that the control is active along $\gamma_2$ and choose a point $Q(\tau)= \gamma_2(\tau,\bar s)$, with $\bar s>0$.  Consider the tangent vectors
$$\bfw_1\,\doteq\,\gamma_{1,s}(\tau,0),\qquad\quad \bfw_2\,\doteq\,\gamma_{1,s}(\tau,0),$$
and the three
inward-pointing unit normal vectors
(see Fig.~\ref{f:sm81})
\bel{unv}
\bfn_1\,\doteq\, {\gamma_{1,s}(\tau,0)^\perp\over \bigl|\gamma_{1,s}(\tau,0)\bigr|}\,,
\qquad\quad \bfn_2\,\doteq\, {\gamma_{2,s}(\tau,0)^\perp\over \bigl|\gamma_{2,s}(\tau,0)\bigr|}\,,
\qquad\quad \bfn_3\,\doteq\, {\gamma_{2,s}(\tau,\bar s)^\perp\over \bigl|\gamma_{2,s}(\tau,\bar s)\bigr|}\,.\eeq
To simplify the computations, without loss of generality we can assume that, for some $\rho>0$ small enough:
\begi
\item For $|t-\tau|\leq\rho$ and $s\in [-\rho, 0]$, the arc $s\mapsto \gamma_1(t,s)$ 
is parameterized by arc-length.
\item For $|t-\tau|\leq\rho$ and $s\in [0, \rho]$, the arc $\gamma_2$ is parameterized by arc-length.
\item For $|t-\tau|\leq\rho$ and $s\in [\bar s-\rho, \bar s+\rho]$, the arc $\gamma_2$ is parameterized according to 
a system of coordinates with coordinate axes parallel to $\bfn_3, \bfn_3^\perp$.
More precisely (see Fig.~\ref{f:sm82}), this means:
\bel{ncoord}
\langle \gamma_{2,t} ,\bfn_3^\perp\rangle\,=\,0~~\hbox{and}~~
 \langle \gamma_{2,s} ,\bfn_3^\perp\rangle\,=-1\qquad \hbox{for}~ |t-\tau|\leq\rho,~|s-\bar s|\leq\rho.
\eeq
\endi
Notice that the above arc-length parameterization implies 
$$|\bfw_1|= |\bfw_2|=1,\qquad \bfn_1=\bfw_1^\perp,\quad \bfn_2=\bfw^\perp.$$

Given $0<\ve<\!<\delta$, for $t\in [\tau-\ve, \tau+\delta+\ve]$ consider the points
\bel{At} P_1(t) ~=~\left\{ \bega{cl} \gamma_1\bigl(t, \, \tau-t-\ve\bigr) \quad & \hbox{if} ~~t\in [\tau-\ve,\, \tau],\cr
\gamma_1(t, -\ve) \quad &\hbox{if}~~ t\in [\tau, \,\tau+\delta],\cr
\gamma_1\bigl(t, \, t-\tau -\delta -\ve) \bigr)\quad& \hbox{if}~~ t\in [\tau+\delta,\, \tau+\delta+\ve],\enda\right.\eeq
\bel{Bt}
 P_2(t) ~=~\left\{ \bega{cl} \gamma_2\bigl(t, \,t-\tau+\ve\bigr) \quad & \hbox{if} ~~t\in [\tau-\ve, \,\tau],\cr
\gamma_2(t, \ve) \quad &\hbox{if}~~ t\in [\tau,\, \tau+\delta],\cr
\gamma_2\bigl(t,\, \tau+\delta+\ve-t \bigr)\quad& \hbox{if}~~ t\in [\tau+\delta,\, \tau+\delta+\ve].\enda\right.\eeq
Observe that $P_1(t) = P_2(t)=P(t)$ for $t=\tau-\ve$ and for $t=\tau+\delta+\ve$. 
We now construct a modified set motion $t\mapsto\Hat \Omega(t)$ by the following rules.
\begi
\item[(i)]  For $t\notin [\tau-\ve, \tau+\delta+\ve]$, one has $\Hat\Omega(t)~=~\Omega(t)$.
\item[(ii)] If $\bfw_1\times \bfw_2>0$ (an outward corner), then for $t\in [\tau-\ve, \tau+\delta+\ve]$
the set $\Hat \Omega(t)$ is obtained from $\Omega(t)$ by removing the triangular region with vertices $P_1(t), P(t), P_2(t)$,
as shown in Fig.~\ref{f:sm83}.
\item[(iii)] If $\bfw_1\times \bfw_2<0$ (an inward corner), then for $t\in [\tau-\ve, \tau+\delta+\ve]$
the set $\Hat \Omega(t)$ is obtained from $\Omega(t)$ by adding the triangular region with vertices $P_1(t), P(t), P_2(t)$,
as shown in Fig.~\ref{f:sm84}.
\endi

\v
{\bf 2.} To estimate the change in the total effort $\Hat\E(t)$ required by the strategy $\Hat \Omega$, 
we need to compute the inward normal speed $\beta_3$ along  the segment $\ov{P_1 P_2}$. 

\begin{figure}[ht]
\centerline{\hbox{\includegraphics[width=8cm]{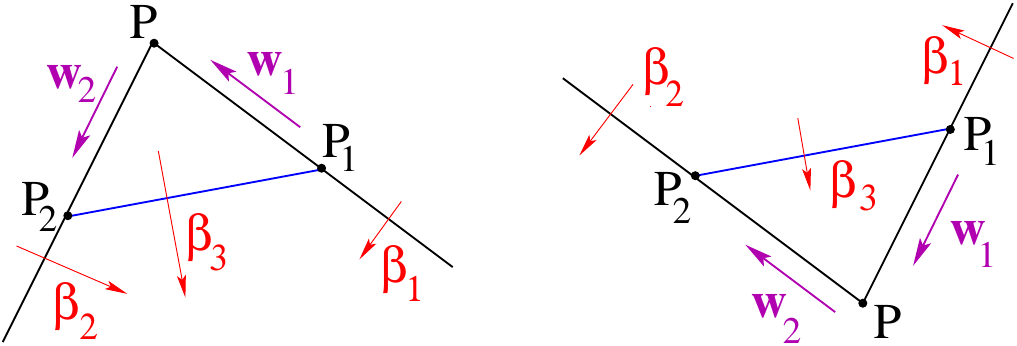}}}
\caption{\small Computing the normal speed $\beta_3$ of the side $AB$, as in (\ref{b123}).}
\label{f:sm85}
\end{figure}

Referring to Fig.~\ref{f:sm85}, left, 
consider a triangle with vertices 
$$P(t),\qquad  P_1(t)\,=\, P(t) - \bfw_1,\qquad  P_2(t)\,=\, P(t) + \bfw_2\,.$$ 
Call $\beta_1, \beta_2,\beta_3$ respectively the normal speeds of the three sides $\ov{P_1P}$,
 $\ov{P_2P}$ and $\ov{P_1P_2}$. 
Knowing the velocity  $\dot P= dP/dt$, these are computed by
\bel{b123}\beta_1~=~\la\dot P, \bfw_1^\perp\ra,\qquad \beta_2~=~\la\dot P, \bfw_2^\perp\ra,
\qquad \beta_3~=~\left\langle \dot P, {(\bfw_1 + \bfw_2)^\perp\over |\bfw_1+\bfw_2|} \right\rangle ~=~{  \beta_1 + \beta_2\over |\bfw_1 +\bfw_2| }.\eeq
During the time interval $[\tau, \tau+\delta]$, the change in the total effort is thus 
computed as
\bel{save}
\bega{rl} \Hat \E(t)-\E(t)& =~\ve\Big( |\bfw_1 +\bfw_2| \, E(\beta_3) -  E(\beta_1) -  E(\beta_2)+\O(1)\cdot (\delta+\ve)\Big)
+o(\ve)\\[3mm]
&=\ds ~\ve \left( |\bfw_1 +\bfw_2| \, E\left( {  \beta_1 + \beta_2\over |\bfw_1 +\bfw_2| }
\right) -  E(\beta_1) -  E(\beta_2)
\right)+\O(1)\cdot \delta\ve+o(\ve).\enda\eeq
We claim that, for $\ve,\delta$ sufficiently small, the right hand side of (\ref{save}) is strictly negative.
Indeed, set 
\bel{Ein}\ov \beta~\doteq~{\beta_1+\beta_2\over 2}~>~-1 \,\qquad\qquad 
\lambda~\doteq ~{2\over |\bfw_1 +\bfw_2|}~>~1 \,.\eeq
Notice that the first inequality follows from the assumption that at least one of the normal speeds $\beta_1,\beta_2$ 
is strictly larger than $-1$.  The second inequality is trivially true because $\bfw_1,\bfw_2$ are non-parallel unit vectors.
By (\ref{Ein}) and the convexity of effort function $E(\beta)\doteq\max\{ 0, \,1+\beta\}$ it now follows
\bel{bes}\bega{l}  \ds {1\over 2} \left[{|\bfw_1 +\bfw_2|}  \, E\left( {  \beta_1 + \beta_2\over |\bfw_1 +\bfw_2| }\right) -  E(\beta_1) -  E(\beta_2)\right] =~{1\over \lambda}  E\left( \lambda\cdot { \beta_1 + \beta_2\over 2}\right) - 
{ E(\beta_1) +  E(\beta_2)\over 2} \\[4mm]
\qquad <~\ds E\left( { \beta_1 + \beta_2\over 2}\right) - 
{ E(\beta_1) +  E(\beta_2)\over 2}~\leq~0.
\enda
\eeq
By choosing $0<\ve<\!<\delta$ sufficiently small in (\ref{save}),  our claim is proved.

By the above analysis, for $\tau-\ve\leq t\leq  \tau+\delta+\ve$ we have the bounds
\bel{HEE} \left\{\bega{rll} 
\Hat \E(t)-\E(t)&\leq \,- c_1\ve\qquad&\hbox{for}~~ t\in [\tau, \tau+\delta],\\[2mm]
\Hat \E(t)-\E(t)&\leq ~ c_2\ve\qquad&\hbox{for}~~ t\in [\tau-\ve, \tau]\cup [\tau+\delta,\tau+\ve+\delta],
\enda\right.
\eeq
for suitable constants $c_1, c_2>0$.
\v
{\bf 3.} In this step we construct a further set $\Tilde \Omega(t)$, adding a perturbation
in a neighborhood of the point $Q(\tau)=\gamma_2(\tau, \bar s)$ where the boundary is smooth
(see Figures \ref{f:sm83} and \ref{f:sm84}, right). At time $t\in [\tau-\ve, \, \tau+\delta+\ve]$, this perturbed boundary is described by 
\bel{pbo}
\Tilde \gamma_2(t,s) ~=~\gamma_2(t,s) + \sigma(t)\vp_\rho(s-\bar s)\bfn_3,\eeq
where $\bfn_3$ is the unit normal vector at (\ref{unv}), while $\vp_\rho(s) = \vp(s/\rho)$ is the function considered at (\ref{presc}).
Assuming 
\bel{small}  0<\ve<\!<\rho<\!<\delta<\!<1,\eeq
by a suitable choice of $\sigma(\cdot)$ and $\rho>0$, we claim that this motion $\Tilde\Omega$  is admissible and yields a strictly 
lower cost than $\Omega$.
\v
{\bf 4.} To compare the effort  of the two strategies $\Tilde\Omega(\cdot) $ 
and $\Hat\Omega(\cdot)$,  we observe that
along the boundary arc $\gamma_2$ by continuity one has
\bel{arc2} \bigl| \bfn(t,s)-\bfn_3\bigr|\,=\,o (t-\tau) + o(s-\bar s),
\qquad\qquad \bigl| \gamma_{2,s}(t,s)-1\bigr|\,=\,o (t-\tau) + o(s-\bar s).\eeq
Recalling (\ref{vpr1}), by (\ref{pbo}) the difference in the effort of the two strategies in a 
neighborhood of $Q(\tau)=\gamma_2(\tau,\bar s) $  can be estimated by
\bel{eff5}
\int_{\bar s-\rho}^{\bar s+\rho} \bigl[\Tilde E(t,s) - E(t,s)\bigr]\, ds~=~
{\rho\over 2} \dot\sigma(t) + \rho \Big(\bigl|\sigma(t)\bigr|+\bigl|\dot\sigma(t)\bigr|\Big)
\cdot o(\rho+\ve+\delta).\eeq
We now define $\sigma:[\tau-\ve, \tau+\delta+\ve]\mapsto \R$  to be the unique 
 piecewise affine function
such that
\bel{sigdef}
\left\{ \bega{rl} 
\sigma(\tau-\ve)~=& \sigma(\tau+\delta+\ve)~=~0,\\[1mm]
\dot \sigma\,=\, -3c_2\ve/ \rho\qquad &\hbox{if} ~~t\in [\tau-\ve, \tau],\\[1mm]
\dot \sigma\,=\,c_1\ve/2\rho\qquad &\hbox{if} ~~t\in [\tau, \tau+\delta],\\[1mm]
\dot\sigma\,=\,-c_3\ve/\rho\qquad &\hbox{if} ~~t\in [\tau+\delta, \tau+\delta+\ve],
\enda\right.
\eeq
for some constant $c_3>2c_2$, uniquely determined by the terminal condition
$\sigma(\tau+\delta+\ve)=0$.

Using (\ref{eff5}), one now checks that, with the choice of the constants $\ve,\delta,\rho$ as in (\ref{small}),
one has $\Tilde\E(t)\leq \E(t)=M$ for all $t$.   Hence this new strategy is admissible.
\v
{\bf 5.}
It remains to check that the perturbed strategy $\Tilde \Omega$ yields a strictly
lower total cost. Since the perturbation is effective only over the time interval 
$[\tau-\ve,\tau+\delta+\ve]$, we have to show that
\bel{TE2} 
\int_{\tau-\ve}^{\tau+\delta+\ve} \Big[\caL^2\bigl(\Omega(t)\bigr) -\caL^2\bigl(\Tilde\Omega(t)\bigr) 
\Big]dt~>~0.\eeq
At any time $t$, the triangular region with vertices $P, P_1, P_2$ has area $\leq \ve^2$.

Moreover, comparing the portions  of the sets $\Omega(t)$ and $\Tilde \Omega(t)$ contained in 
a  suitable ball $B(Q,r)$ centered at the point $Q= \gamma_2(\tau,\bar s)$, we find
\bel{TE3} \caL^2\bigl(\Omega(t)\cap B(Q,r)\bigr) -\caL^2\bigl(\Tilde\Omega(t)\cap B(Q,r)\bigr) ~=~
{\rho\over 2} \Big( \sigma(t) + o\bigl(\sigma(t)\bigr)\Big).\eeq
We now observe that
\begi
\item On the interval $[\tau-\ve,\tau]$, one has $\sigma(t)=\O(1)\cdot\ve^2/\rho$.
\item On the middle interval, choosing $\ve>0$ much smaller than the other parameters, we have
$$\bega{l} \ds\int_\tau^{\tau+\delta} {\rho\over 2} \Big( \sigma(t) + o\bigl(\sigma(t)\bigr)\Big) ~\geq ~ {\rho\over 2}  \int_\tau^{\tau+\delta}\Big[\O(1)\cdot {\ve^2 \over\rho} + {c_1\ve\over 3\rho} (t-\tau)\Big]dt\\[4mm]
\qquad\ds~=~{\rho\over 2} \left[\O(1)\cdot {\delta \ve^2\over \rho} + c_1{\ve \delta^2\over 6\rho}\right]~>~0.
\enda
$$
\item On the final  interval $[\tau-\ve,\tau+\delta, \tau+\delta+\ve]$, one has $\sigma(t)\geq 0$.
\endi
Combining  the above estimates we conclude that $J(\Tilde \Omega) < J(\Omega)$, 
showing that the original strategy was not optimal.\endproof

\section{Optimality conditions at boundary points}
\label{s:11}
\setcounter{equation}{0}

\begin{theorem} \label{t:111} Let $t\mapsto\Omega(t)$ be an optimal strategy for the 
minimization problem {\bf (OP)}, with $\kappa_1>0$.
Assume that, for $t$ in a neighborhood of $\tau$,   the boundary 
$\partial \Omega(t)$ contains a $\C^1$ arc $\gamma=\gamma(t,s)$
which meets the boundary $\partial V$ at one of its endpoints, 
say at $\gamma(t,0)$. Then
\begi
\item[(i)] either the junction is perpendicular, 
\item[(ii)] or else at the junction point the control effort vanishes, i.e., the inward normal speed is $\beta(t,s)=-1$. 
\endi 
\end{theorem}

{\bf Proof.} The main argument is similar to the one used in the proof of Theorem~\ref{t:101}.
Assuming that the intersection is not perpendicular and the control effort does not vanish,
we first perform a perturbation $\Hat \Omega$ which decreases the length of the boundary, so that  the 
total effort is strictly reduced: $\Hat \E(t)<M$, for $t\in [\tau, \tau+\delta]$.   
We then use the additional available effort to remove from $\Hat \Omega(t)$ a region in a neighborhood
of some other point $Q=\gamma(\tau, \bar s)$.   This will yield an admissible motion $t\mapsto \Tilde \Omega(t)$
with strictly smaller total cost.

\v
{\bf 1.} 
To fix ideas, assume that, near the boundary of $V$,  the curve $s\mapsto \gamma(t,s)$ 
is parameterized by arc-length,
so that $\gamma(t,0)\in \partial V$, $\gamma(t,s)\in \partial \Omega(t)$ for $s\geq 0$ small.

\begin{figure}[ht]
\centerline{\hbox{\includegraphics[width=16cm]{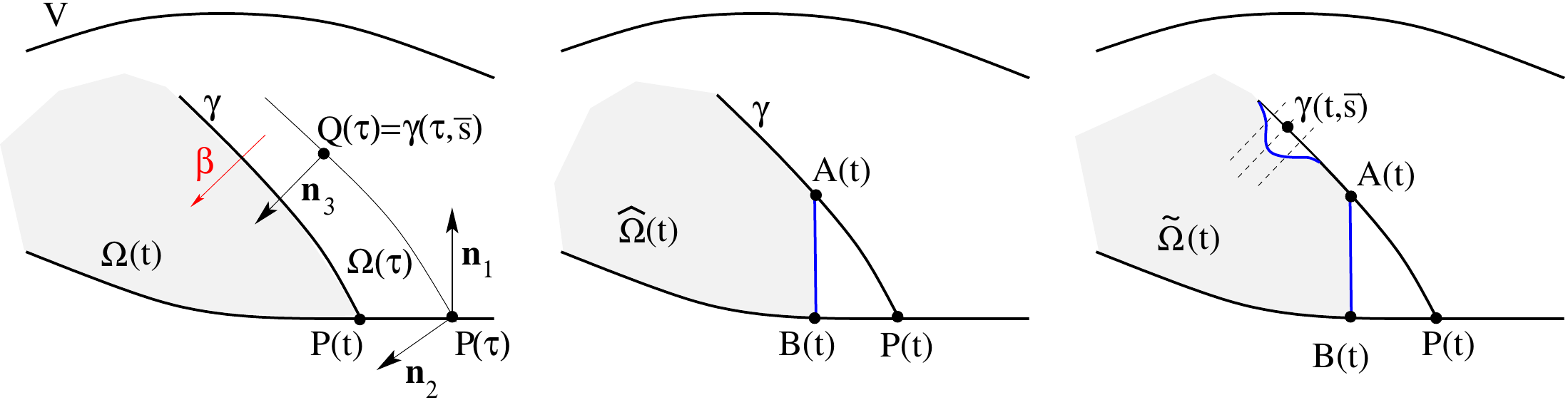}}}
\caption{\small Left: the original sets $t\mapsto\Omega(t)$. Center: the sets $\Hat\Omega(t)$ obtained by 
removing the triangular region $\Hat{ PAB}$. Right: the sets $\Tilde \Omega(t)$ obtained by further perturbing the boundary in a neighborhood of the point  $Q(\tau)$ }
\label{f:sm86}
\end{figure}

As shown in Fig.~\ref{f:sm86}, left, consider the points $P(\tau)= \gamma(\tau, 0)$ and choose some other 
point $Q(\tau)=\gamma(\tau,\bar s)$ with $\bar s>0$.   
The following notation will be used:
\begi
\item $\bfn_1$ is the inward pointing, unit vector, perpendicular to the boundary $\partial V$ at $P(\tau)$.
\item $\bfn_2= {\gamma_s(\tau, 0)^\perp\over \bigl|\gamma_s(\tau, 0)\bigr|} $ is the  unit vector perpendicular 
to $\gamma$ at $P(\tau)$.
\item $\bfn_3$  is the inward pointing, unit vector, perpendicular to the boundary $\partial \Omega(\tau) $ at $Q(\tau)$.
\endi
\v

{\bf 2.} Let $\Omega$ be an optimal strategy, and assume, by contradiction, that the intersection
at $P(\tau)$ is not perpendicular.   To derive a contradiction, we  first consider the case where
$\langle \bfn_1,\bfn_2\rangle <0$, so that the tangent cone to $\Omega(\tau)$ at $P(\tau)$ has an opening
$\theta <\pi/2$.

As shown in Fig.~\ref{f:sm86}, center, 
given $0<\ve<\!<\delta$, for $t\in [\tau-\ve, \tau+\delta+\ve]$ we consider the points
\bel{Atdef} A(t) ~=~\left\{ \bega{cl} \gamma\bigl(t, \, t-\tau+\ve\bigr) \quad & \hbox{if} ~~t\in [\tau-\ve,\, \tau],\cr
\gamma(t, \ve) \quad &\hbox{if}~~ t\in [\tau, \,\tau+\delta],\cr
\gamma\bigl(t, \, \tau +\delta +\ve-t) \bigr)\quad& \hbox{if}~~ t\in [\tau+\delta,\, \tau+\delta+\ve],\enda\right.\eeq
We then define $B(t)$ to be the projection of $A(t)$ on the boundary of $V$.   More precisely, $B(t)$ 
is the unique point which satisfies
\bel{Bdef} B(t)~=~A(t) + h(t) \bfn_1\,,\qquad\qquad B(t)\in \partial V.\eeq
By the implicit function theorem, such a point is uniquely determined.

Finally, we let $\Hat\Omega(t)$ be the set obtained by removing from $\Omega(t)$ the triangular set
with vertices $P(t), A(t), B(t)$ (see Fig.~\ref{f:sm86}, center).    

Performing a similar analysis as in step {\bf 2} of the proof of Theorem~\ref{t:101}, 
for $\tau-\ve\leq t\leq  \tau+\delta+\ve$ we obtain the bounds on the total effort:
\bel{HEE2} \left\{\bega{rll} 
\Hat \E(t)-\E(t)&\leq \,- c_1\ve\qquad&\hbox{for}~~ t\in [\tau, \tau+\delta],\\[2mm]
\Hat \E(t)-\E(t)&\leq ~ c_2\ve\qquad&\hbox{for}~~ t\in [\tau-\ve, \tau]\cup [\tau+\delta,\tau+\ve+\delta],
\enda\right.
\eeq
for suitable constants $c_1, c_2>0$.

\begin{figure}[ht]
\centerline{\hbox{\includegraphics[width=11cm]{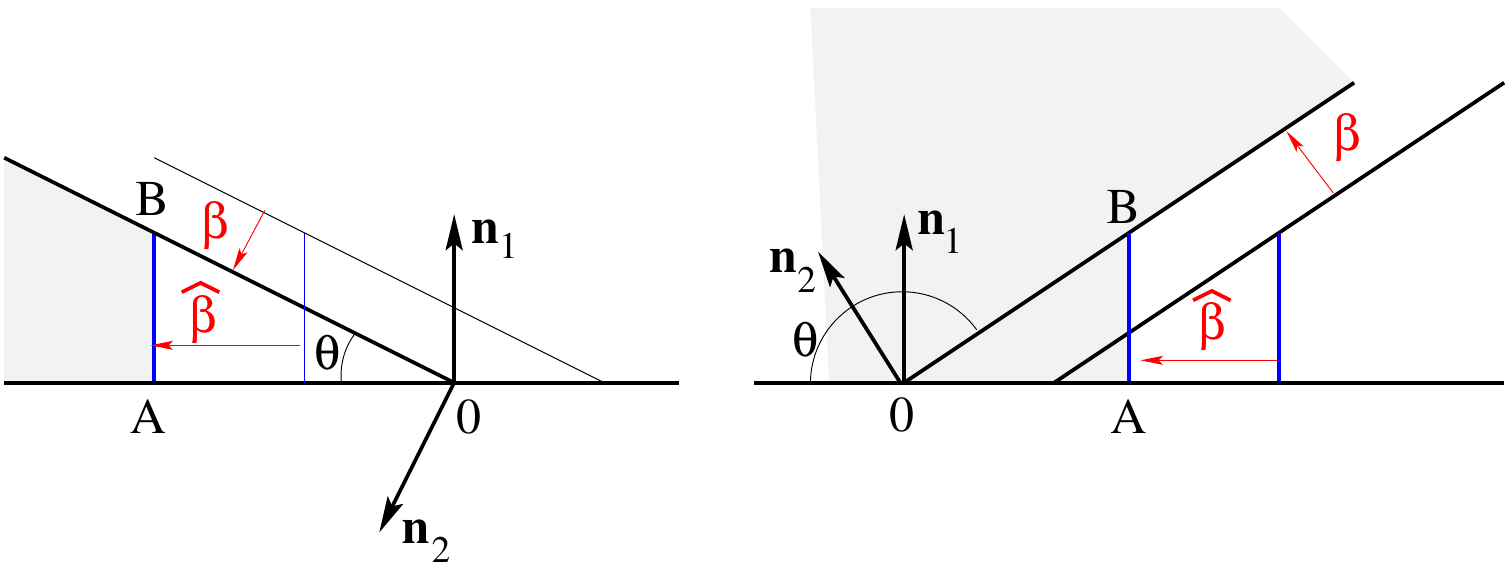}}}
\caption{\small Computing the reduction in the total effort as in (\ref{effred}). Left: the case where $0<\theta<\pi/2$.
Right: the case where $\pi/2<\theta<\pi$. }
\label{f:sm90}
\end{figure}

Indeed,  referring to Fig.~\ref{f:sm90} left, for $t\in [\tau, \tau+\delta]$ 
to leading order the difference in the instantaneous effort
is measured by 
$ E(\Hat \beta) \cdot |AB|-E(\beta) \cdot |0B| $.
We have 
$$|AB| \,=\, |0B|\sin\theta, \qquad \qquad\Hat\beta\,=\, {\beta\over\sin\theta},$$
\bel{effred}\bega{l}\ds E(\Hat \beta) \cdot |AB|-E(\beta) \cdot |0B|~=~\max\left\{ 1+ {\beta\over\sin \theta},0\right\}
\sin \theta \cdot  |0B|-(1+\beta) \cdot  |0B|\\[4mm]
\qquad\ds =~\max \{ \sin\theta+ \beta,0\} \cdot  |0B|-(1+\beta) \cdot  |0B|~<~0,\enda
\eeq
because we are assuming $\beta>-1$. This yields the first inequality in (\ref{HEE2}). The second one is clear,
choosing $c_2$ sufficiently large.

\v
{\bf 3.} The remainder of the proof is the same as for the previous theorem.
Without loss of generality, we can assume that, for  $|t-\tau|\leq\rho$ and $s\in [\bar s-\rho, \bar s+\rho]$, the arc $\gamma$ is parameterized according to 
a system of coordinates with coordinate axes parallel to $\bfn_3, \bfn_3^\perp$.
This means:
\bel{nco2}
\langle \gamma_{t} ,\bfn_3^\perp\rangle\,=\,0~~\hbox{and}~~
 \langle \gamma_{s} ,\bfn_3^\perp\rangle\,=-1\qquad \hbox{for}~ |t-\tau|\leq\rho,~|s-\bar s|\leq\rho.
\eeq
We now construct a further set $\Tilde \Omega(t)$, adding a perturbation
in a neighborhood of the point $Q(\tau)=\gamma_2(\tau, \bar s)$ where the boundary is smooth
(see Fig.~\ref{f:sm86}, right). At time $t\in [\tau-\ve, \, \tau+\delta+\ve]$, this perturbed boundary is described by 
\bel{pbo2}
\Tilde \gamma(t,s) ~=~\gamma(t,s) + \sigma(t)\vp_\rho(s-\bar s)\bfn_3,\eeq
where $\vp_\rho(s) = \vp(s/\rho)$ is the function at (\ref{presc}), 
and $\sigma:[\tau-\ve, \tau+\delta+\ve]\mapsto \R$ is the
 piecewise affine function
defined at (\ref{sigdef}).   Repeating the analysis in the last two steps of the proof of Theorem~\ref{t:101},
 we conclude that the set motion $t\mapsto \Tilde \Omega(t)$ is admissible and yields a strictly smaller
 total cost.   
\v
\begin{figure}[ht]
\centerline{\hbox{\includegraphics[width=16cm]{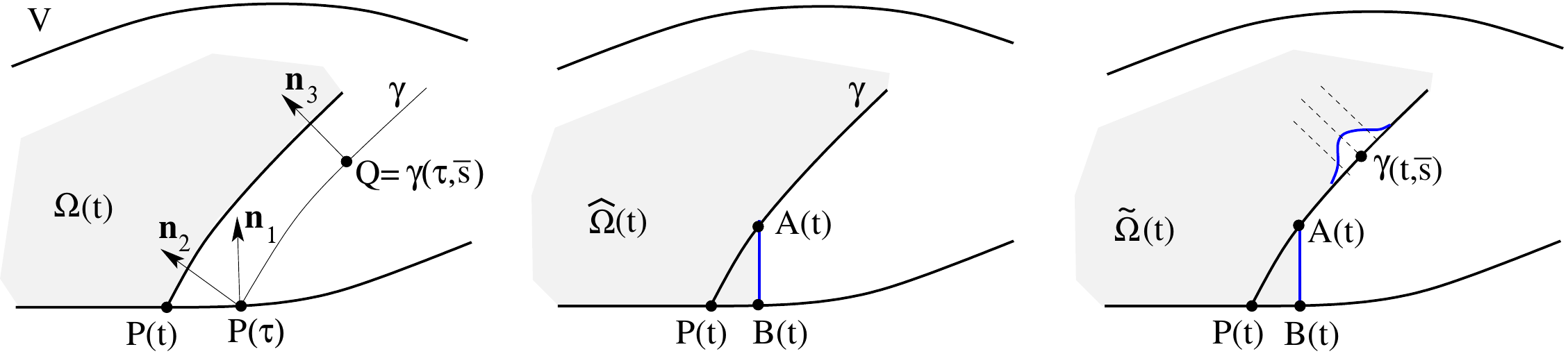}}}
\caption{\small The sets $\Omega(t)$, $\Hat\Omega(t)$ and  $\Tilde\Omega(t)$ in the case where
the tangent cone to the set $\Omega(\tau)$ at $P(\tau)$ covers an angle $\theta>\pi/2$. }
\label{f:sm87}
\end{figure}

{\bf 4.} It remains to consider the case where $\langle \bfn_1,\bfn_2\rangle >0$, 
so that the tangent cone to $\Omega(\tau)$ at $P(\tau)$ has an opening
$\theta>\pi/2$, as shown in Fig.~\ref{f:sm87}.  The points $A(t)$, $B(t)$ are defined as before.
The only difference is that the set $\Hat\Omega(t)$ is obtained by adding to $\Omega(t)$
the triangular region $\Hat{PAB}$.  During the time interval $[\tau, \tau+\delta]$ this reduces 
the length of the boundary $\partial \Hat\Omega(t)\cap V$, hence reducing the effort: $\Hat \E(t)<\E(t)=M$.

As before, the additional available effort is used  to remove from $\Hat \Omega(t)$ a region in a neighborhood
of some other point $Q=\gamma(\tau, \bar s)$.   This will yield an admissible motion $t\mapsto \Tilde \Omega(t)$
with strictly smaller total cost.
\endproof

\begin{remark}
{\rm In the above construction, the sets $\Tilde\Omega (t)$ coincide with 
$\Omega(t)$ for $t\notin [\tau-\ve, \tau+\delta+\ve]$.
Hence the modified strategy would achieve eradication within the same minimum time.

However, the same orthogonality  condition can be shown to hold for the minimum time problem.
Indeed, for $t\in [\tau, \tau+\delta]$ we can use the additional available effort $M-\Tilde \E(t)>0$ 
to construct a perturbed strategy $\Omega^\sharp (t)$ so that $\Omega(\tau+\delta) = \Omega^\sharp (\tau')$
for some $t'<\tau+\delta$.  

More in detail, consider a new strategy obtained by shifting time:
$\Omega^\sharp(t) = \Tilde\Omega\bigl(t^\sharp(t)\bigr)$.
The inward normal speed, at any boundary point $x\in \partial \Omega^\sharp(t^\sharp(t))$ is 
$$\beta^\sharp(t^\sharp,x)\,=\,\Tilde \beta(t,x)\cdot {dt^\sharp\over dt}\,.$$
Hence the corresponding total effort is 
\bel{effsharp}\E^\sharp\bigl(t^\sharp\bigr)~=~\int_{\partial \Omega^\sharp(t^\sharp)}  \max\left\{1+ \beta^\sharp(t,x),~0\right\}\, d\H^1(x)~=~
\int_{\partial\Tilde  \Omega(t)}  \max\left\{1+ \Tilde\beta(t,x)\cdot {dt^\sharp\over dt}\,,~0\right\}\, d\H^1(x).\eeq
We recall that, for $t\in [\tau, \tau+\ve]$, the previous construction yields
$$\Tilde\E(t)~=~\int_{\partial\Tilde  \Omega(t)}  \max\left\{1+ \Tilde\beta(t,x),~0\right\}\, d\H^1(x)~<~M.$$
By continuity, we can thus choose a function $t\mapsto t^\sharp(t)$ with 
$$\left\{ \bega{rl} dt^\sharp/ dt\,>\,1\quad&\hbox{for}~ \tau<t<\tau+\ve,\\[1mm]
dt^\sharp/ dt\,=\,1\quad&\hbox{for}~ t>\tau+\ve,\enda\right.
\qquad\qquad t^\sharp(\tau)=\tau,$$
and such that the right  hand side of (\ref{effsharp}) remains $\leq M$ for $t\in [\tau, \tau+\ve]$.
 This contradicts the optimality of the original strategy, for the minimum time problem.
}
\end{remark}

\begin{figure}[ht]
\centerline{\hbox{\includegraphics[width=15cm]{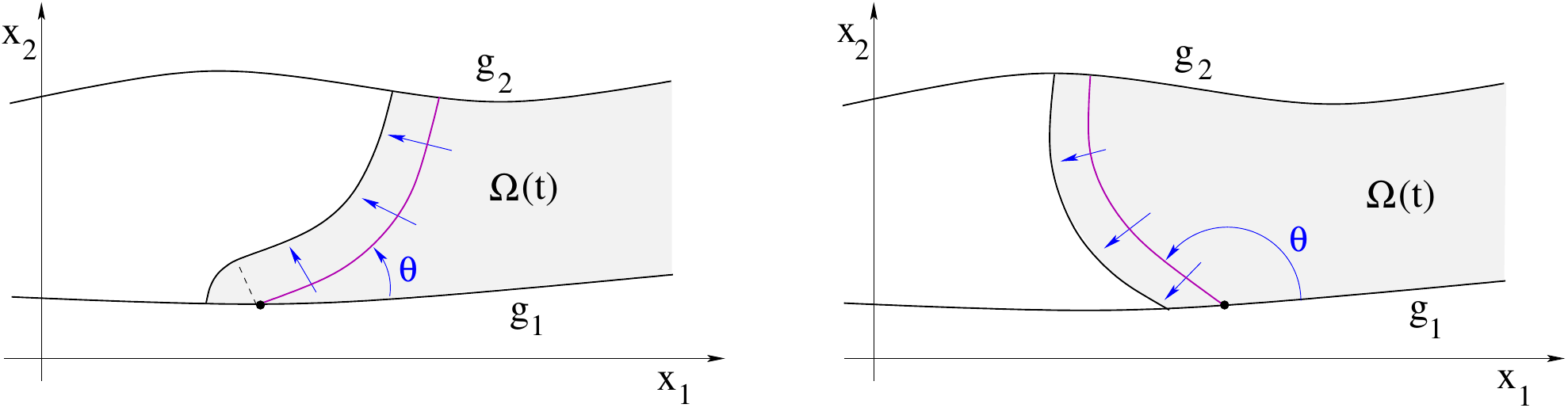}}}
\caption{\small   Left: if no control is active, an angle $\theta<\pi/2$ is immediately replaced by a perpendicular intersection.
Right: An angle $\theta>\pi/2$ can persist, even if no control is active.}
\label{f:tg20}
\end{figure}

\begin{remark}
{\rm If the control is not active on a portion of the boundary of  $\Omega(t)$ which touches the boundary of $V$,
the angle $\theta$
between the two boundaries $\partial \Omega(t)$ and $\partial V$ immediately becomes 
$\geq \pi/2$ (see Fig.~\ref{f:tg20}).  However, orthogonality may fail.
}
\end{remark}

%
%
%

\section{Optimal strategies determined by the necessary conditions}
\label{s:12}
\setcounter{equation}{0}
In the remaining sections of this paper we study
motions $t\mapsto \Omega(t)\subseteq V$
that satisfy all the necessary conditions derived in the previous sections.  Of course, these motions 
will be natural candidates for optimality.
However, we should point out that our necessary conditions are all based on regularity 
assumptions that may not be always verified.
It is only in the setting of Theorem~\ref{t:41} that our constructions are guaranteed to 
yield the globally optimal solutions.

Throughout the following we consider an admissible strategy $t\mapsto\Omega(t)\subseteq V$ such that the following holds.
\begi
\item[{\bf (A6)}] {\it For every $t\in [0,T]$, the relative boundary $\partial \Omega(t)\cap V$ is the concatenation
of finitely many $\C^2$ arcs 
$$\gamma_i(t)~=~\Big\{ \gamma_i(t,\xi)\,;~~\xi\in \bigl[a_i(t), b_i(t)\bigr]\Big\}, \qquad\qquad i=1,\ldots,N.$$
 Each $\gamma_i$ can be
\begi\item either a {\bf free arc}, where the inner normal velocity is $\beta=-1$ at every point,
\item or a {\bf controlled arc}, where $\beta>-1$ at every point.
\endi
}
\endi
Notice that, along a free arc, no control effort is present. Along this portion of its boundary, the set $\Omega(t)$
thus freely expands with unit speed.   On the other hand, along a controlled arc, the evolution of the boundary
depends on the pointwise control effort $E= 1+\beta(t,\xi)>0$.

Given an initial set $\Omega_0\subseteq V$,  we shall seek admissible set motions $t\mapsto \Omega(t)$
which satisfy the regularity assumptions {\bf (A6)} together with the optimality conditions derived in previous sections.
\begi
\item[{\bf (A7)}] {\it At a.e.~time $t\in [0,T]$ the following optimality conditions hold.
\begi
\item[(i)] The area of the set $\Omega(t)$ varies according to 
\bel{dareat}{d\over dt} \caL^2\bigl(\Omega(t)\bigr) ~=~\H^1\bigl( \partial \Omega(t)\cap V\bigr) -M.\eeq
\item[(ii)] All controlled arcs $\gamma_i$ have same constant curvature: $\omega_i(t,\xi)=\omega(t)$ for all 
$\xi$.
\item[(iii)] Controlled arcs join tangentially with free arcs, at their endpoints in the interior of~$V$.
\item[(iv)] At a point $Q= \gamma_i(t, \xi)$ where a controlled arc touches the boundary $\partial V$, 
 either the junction is perpendicular, or else the effort vanishes: $E\bigl(\beta(t,\xi)\bigr)=0$.\endi
}
\endi

\begin{figure}[ht]
\centerline{\hbox{\includegraphics[width=6cm]{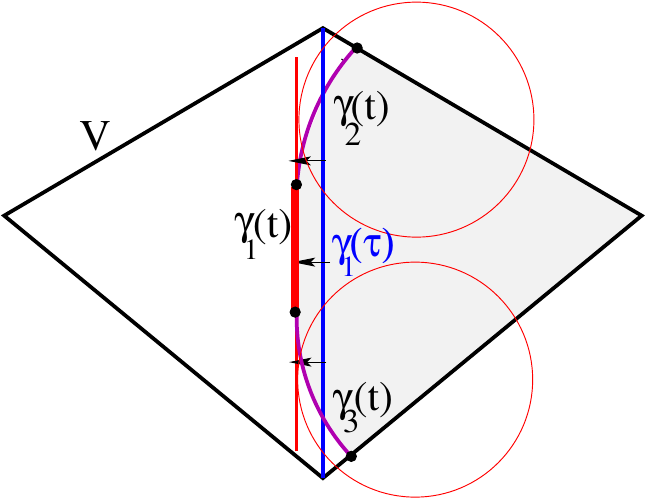}}}
\caption{\small  Here the segment $\gamma_1(\tau)$ is a maximal free arc.   At the later time $t>\tau$, the boundary
$\partial \Omega(t)\cap V$ is the union of a free arc $\gamma_1$ and the controlled arcs, $\gamma_2,\gamma_3$.
We first construct $\gamma_1(t)$ by shifting the free arc $\gamma_1(\tau)$  in the amount $t-\tau$ in the orthogonal direction.
The controlled arcs $\gamma_2(t), \gamma_3(t)$ are then two arcs of circumferences, with the same radius 
$r(t)$, tangent to $\gamma_1(t)$ and perpendicular to $\partial V$ at the endpoints.
There is a 1-parameter family of such arcs, depending on the radius $r$. A unique radius $r(t)$ can be determined
by imposing the identity (\ref{dareat}) on the rate of decrease of the area. }
\label{f:sm80}
\end{figure}

Assume that, at some time $\tau$, a configuration 
$$\partial \Omega(\tau)=\gamma_1(\tau)\cup\gamma_2(\tau)\cup
\cdots\cup\gamma_N(\tau)$$
is given. Based on the above conditions {\it (i)--(iv)},
for $t$ in a neighborhood of $\tau$ the motion $t\mapsto \Omega(t)$ can then be  determined as follows
(see Fig.~\ref{f:sm80}).

\begi
\item As a first step, we determine the  portion of the boundary $\partial\Omega(t)$
covered by free arcs.  Since the inward normal speed is $\beta=-1$, these free arcs $\gamma_i(t)$
can be constructed by shifting the points
$\gamma_i(\tau,\xi)$ in the normal direction $\bfn_i(t,\xi)$, in the amount $\tau-t$.
\item As a second step, we construct the controlled arcs $\gamma_j(t)$.    These are arcs of circumferences,
all with the same radius $r(t)$. At their endpoints they are tangent to the free arcs, and perpendicular to 
the boundary $\partial V$.  The common radius $r(t)$ is uniquely determined by the area identity (\ref{dareat}),
which yields an ODE for the time derivative $\dot r(t)$.
\endi

In our approach, a key role is played by  maximal extended free arcs.
\begin{definition}\label{d:121}
We say say that a free arc $\gamma_i(\tau)$ is a {\bf maximal free arc} if, for all times $t$ in a neighborhood of $\tau$ one has
\bel{mfa}
\gamma_i(t)~\subseteq~\Big\{ \gamma_i(\tau,\xi) - (t-\tau)\bfn_i(\tau,\xi)\,;~~\xi\in [a_i(\tau), b_i(\tau)]\Big\}.\eeq
Here $\bfn_i(t,\xi)$ denotes the unit inner normal vector at the point $\gamma_i(\tau,\xi)\in \partial \Omega(t)
\cap V$.
\end{definition}
 
If the maximal free arcs $\gamma_i(\tau)$ can be found, in turn these determine the free arcs 
$\gamma_i(t)$ at all other times $t\in [0,T]$.  In a further step, we can then construct the family of 
circumferences $\gamma_j(t)$ tangent to the free arcs and perpendicular to the boundary $\partial V$, thus
completely solving the optimal set motion problem (see Fig.~\ref{f:sm80}).
In the remaining sections we shall focus in more detail
on some aspects of this construction.   The next simple example illustrates the main ideas.

\begin{figure}[ht]
\centerline{\hbox{\includegraphics[width=15cm]{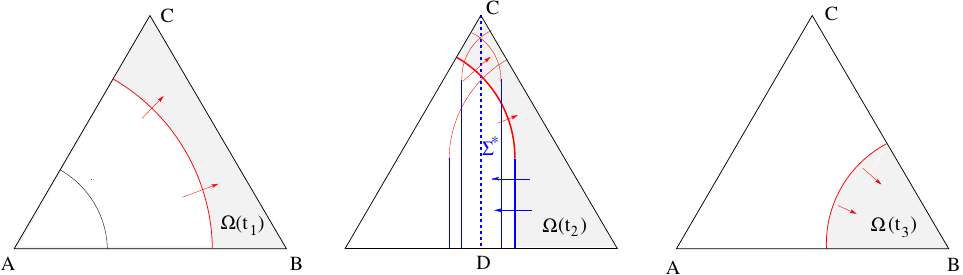}}}
\caption{\small  The (conjectured) optimal strategy, for the minimum time 
eradication problem, where the set $\Omega(t)$ is confined within the triangle $ABC$.
Left: the set $\Omega(t_1)$ for $0<t_1<t^*$.  Center: the set $\Omega(t_2)$ for  $t^*<t_2<T/2$.   Right:
the set $\Omega(t_3)$ for  $ T-t^*<t_3<T$.
}
\label{f:sc26}
\end{figure}
\begin{example}\label{e:121} {\rm In the $x$-$y$ plane, let $V$ be an equilateral triangle with sides of unit length, and vertices
$A=(-1/2,0)$, $B=(1/2,0)$, $C=(0,\sqrt 3/2)$.
In this case, 
 there is no way to 
satisfy the sufficient conditions stated in (\ref{opc}).    
A possible eradication strategy $t\mapsto \Omega(t)$, $t\in [0,T]$ is the one
shown in Fig.~\ref{f:sc26}.  For a suitable $t^*\in\,]0, T/2[\,$ this strategy can be described as follows.
\begi
\item[(i)] For $t\in [0, t^*]$ one has $\Omega(t)= V\setminus B\bigl(A, r(t)\bigr)$, where by (\ref{dareat}) 
the radius $r(t)$ 
evolves according to
\bel{dr1}\dot r(t)~=~{3\over \pi}{M\over r(t)} - 1\,.\eeq
\item[(ii)] For $t^*< t < T/2$, the boundary $\partial \Omega(t)$ is the union of a vertical segment
along the line where $x=x(t)=r(t^*)  - (t-t^*)$, and an arc of circumference.
\item[(iii)] For $t=T/2$, the set $\Omega(T/2)=\bigl\{ (x,y)\in V\,;~~x>0\bigr\}$ is the right half of the triangle.
\item[(iv)] For $t\in [T/2, T]$, the set $\Omega(t)$ can be defined by the symmetry property
$$\Omega(t)~=~\big\{ (x,y)\in V\,;~~(-x,y)\notin \Omega(T-t)\bigr\}.$$
\endi
Calling $r^*\doteq r(t^*)$, we observe that  (\ref{dr1}) implies 
\bel{r*}r^*~<~{3M \over \pi}\,.\eeq
Otherwise $\dot r\leq 0$ and the set $\Omega(t)$ does not shrink.
Moreover, by (iii) it follows
\bel{t*p} {T\over 2} - t^*~=~r^* - {1\over 2}\,.\eeq

\begin{figure}[ht]
\centerline{\hbox{\includegraphics[width=6cm]{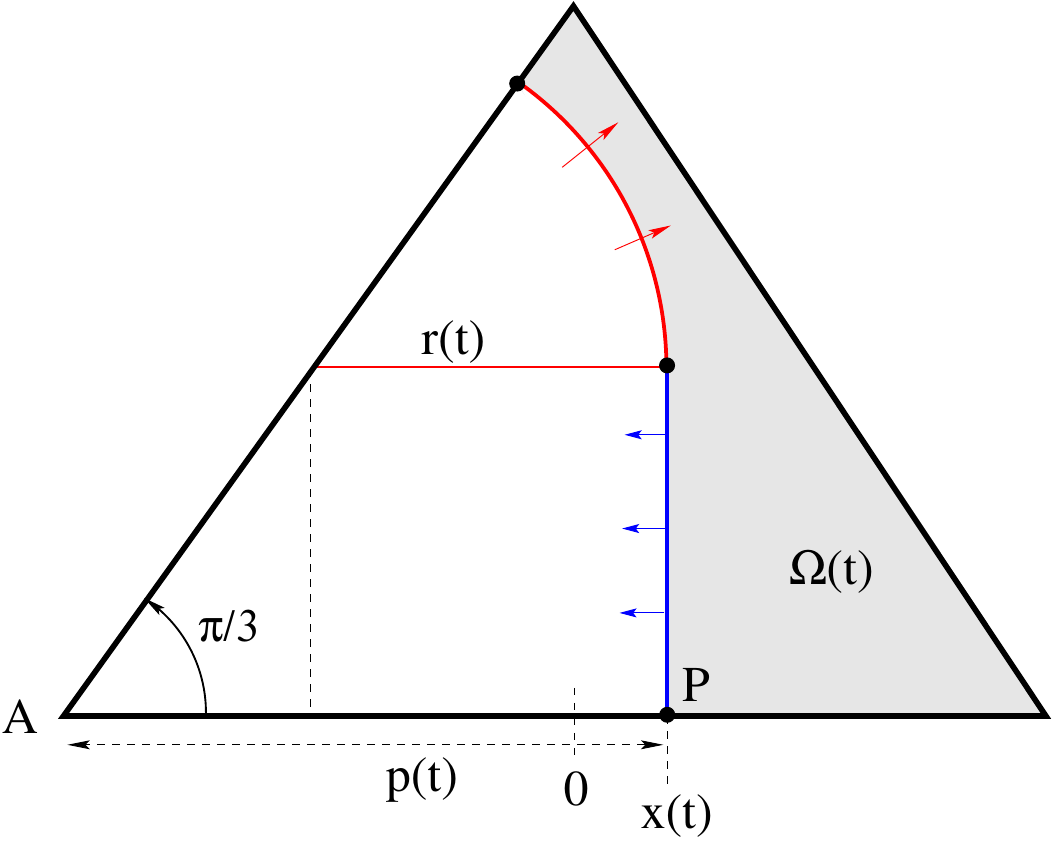}}}
\caption{\small Computing the time derivative of the radius $r(t)$, at (\ref{dr2}). }
\label{f:sc75}
\end{figure}

For $t^*<t<T/2$, using the rate of decrease of the area (\ref{dareat}),  one finds that the radius of the boundary arc shrinks at the rate
\bel{dr2}\dot r(t)~=~-1-\frac{M}{\Big(\sqrt{3}-{\pi\over 3}\Big)r(t)}\,.\eeq
Indeed, as shown in Fig.~\ref{f:sc75}, we have
$$\bega{rl} {{\caL}^2}\bigl(\Omega(t)\bigr) &\ds=~{{\caL}^2}(V) - \left[ {\pi\over 6} r^2(t) +  {\sqrt 3\over 2}\bigl(x(t)-r(t)\bigr)\cdot \bigl[x(t)+r(t)\bigr]\right]\\[4mm]&\ds=~{{\caL}^2}(V) -{\sqrt 3\over 2} x^2(t) +\left( {\sqrt 3\over 2}-{\pi\over 6 } 
\right) r^2(t),\enda
$$
$$\H^1\Big( V\cap \partial \Omega(t)\Big)~=~{\pi\over 3}\, r(t) + \sqrt 3 \bigl(x(t)-r(t)\bigr).$$
Since $\dot x(t)=-1$, from (\ref{dareat}) it follows
$$\left(\sqrt 3 - {\pi\over 3} \right) r(t) \dot r(t)~=~\left({\pi\over 3} - \sqrt 3\right) r(t) - M,$$
which yields (\ref{dr2}).

We seek a solution to (\ref{dr2}) with $r(t^*)= r^*$ and $r(T/2) = 0$, where $r^*$ satisfies (\ref{r*})  and 
$t^*, T$ are related by (\ref{t*p}).  Toward this goal, we observe that
the solution of  (\ref{dr2}) with initial condition
$r(0) = r^\sharp = { 3M / \pi}$ satisfies the implicit equation
\bel{ImR} t=(r^\sharp-r(t))+{1 \over \lambda }\ln\Big({\lambda r(t)+1\over \lambda r^\sharp+1 }\Big)\,, \qquad  \qquad \lambda\doteq {{1\over M}\left(\sqrt 3 - {\pi\over 3} \right) }\,. \eeq

Hence the time $t^\sharp$ at which $r(t^\sharp)=0$ is given by 
$$ t^\sharp= r^\sharp- {1 \over \lambda} \ln\big(  \lambda r^\sharp+1\big)\,. $$
The value of $M$ for which we have the identity $t^\sharp~=~r^\sharp - {1\over 2}$, is 
\bel{tr*}M^\sharp~=~ { (3\sqrt 3-\pi)\over 6}\Big[\ln{\Big({3\sqrt3\over \pi}\Big)}\Big]^{-1}~\approx 0.6805\,.\eeq
The previous analysis shows that   the eradication  problem for the equilateral triangle is solvable  if $M>M^\sharp$. Notice that $\kappa(V)<M^\sharp<K(V)$, in accordance with Theorem~\ref{t:41}.

We conjecture that the above eradication strategy is the optimal one, for the minimum time problem.
Indeed, it satisfies all the necessary conditions for optimality derived in the later sections.
However, at 
present  we cannot rule out the possibility that some other strategy (possibly with lower regularity)
may achieve eradication in a shorter time. }
  \end{example}

\begin{figure}[ht]
\centerline{\hbox{\includegraphics[width=12cm]{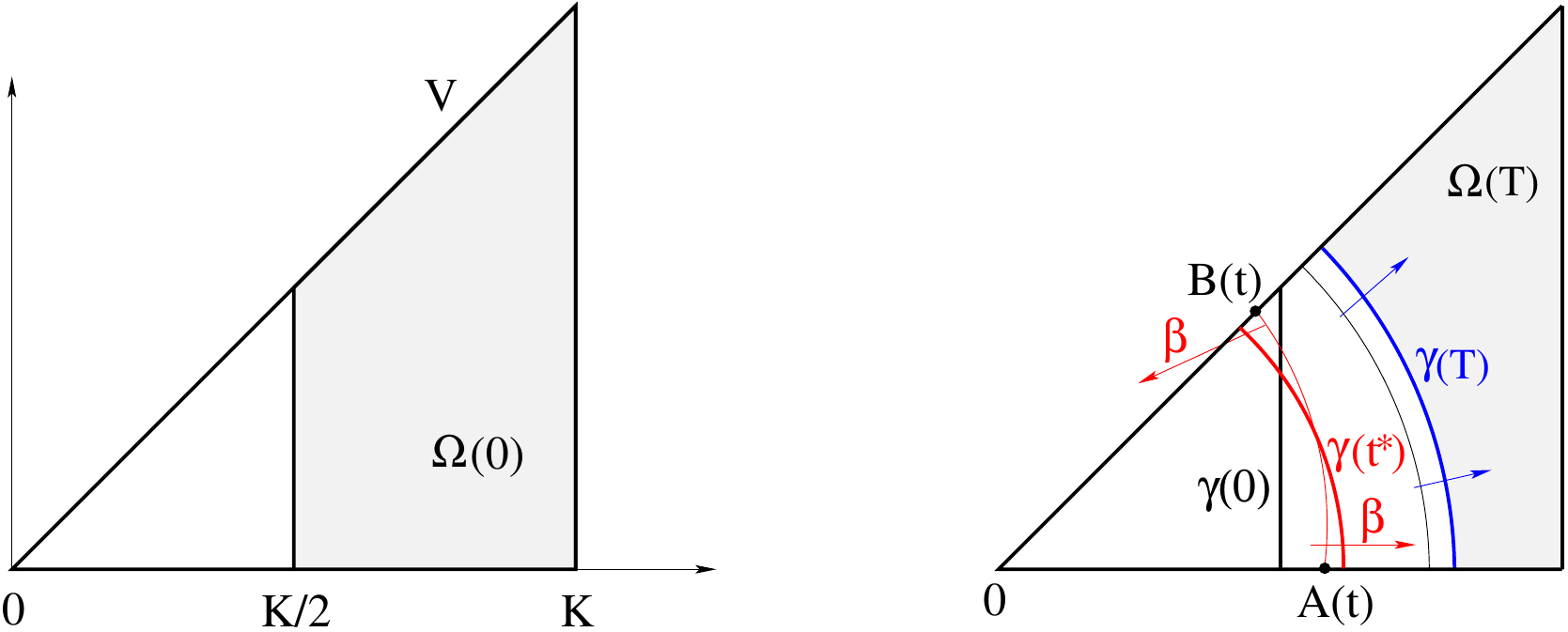}}}
\caption{\small  The optimal strategy in Example~\ref{e:122} }
\label{f:sm88}
\end{figure}

\begin{remark} {\rm If the boundary of the constraining set has a corner $C$, as in Example~\ref{e:121},
the best strategy has a curious behavior.   Namely, as the boundary of the
set $\Omega(t)$ crosses the corner, all the  effort  gets concentrated near the
singular point.  The portion of the boundary where the control is active
shrinks to a point, then grows again. This remains true even if in  (\ref{EM})  the 
constant $M$ is very large.  In other words, even if  there exist admissible eradication strategies
where the set $\Omega(t)$ is monotonically shrinking, these are never  optimal
if the domain $V$ has corners.}
\end{remark}
\v
\begin{example}
\label{e:122}{\rm Let $V$ be the triangle
$$V~=~\bigl\{ (x_1, x_2)\,;~~0<x_1< K, \quad 0<x_2<x_1\bigr\}.$$
As shown in Fig.~\ref{f:sm88}, left,  let the initial contaminated set be
$$\Omega(0)~=~\left\{ (x_1, x_2)\,;~~{K\over 2} < x_1< K, \quad 0<x_2<x_1\right\}.$$
As shown in Fig.\ref{f:sm88}, right.
an optimal strategy for {\bf (OP)} proceeds in two stages.  
For every $t\in [0,T]$, the boundary $\partial \Omega(t)\cap V$
is an arc of circumference $\gamma(t)$ with endpoints $A(t), B(t)$.
\begi
\item On an initial time interval  $t\in [0,t^*]$, this arc is perpendicular to $\partial V$ at $A(t)$,
while at $B(t)$ the control effort is zero.  This first stage  is continued up to the time
$t^*$ where the intersection at $B(t)$ becomes perpendicular as well.
\item For  $t\in [t^*, T]$, $\gamma(t)$ is an arc of circumference centered at the origin,
perpendicular to $\partial V$ at both endpoints.
\endi
Notice that this motion satisfies all the necessary conditions for optimality proved in the pervious sections.
The arcs $\gamma(t)$ are uniquely determined by the above conditions, together with the 
identity $\E(t)=M$ on the total effort (\ref{CEt}).
}\end{example}

%
%

\section{The initial stages of an optimal strategy}
\label{s:13}
\setcounter{equation}{0}
In this section we focus on the initial stages of an optimal eradication strategy.
Assuming that $V$ is an open set with piecewise smooth boundary, we seek an admissible
set motion $t\mapsto \Omega(t)\subseteq V$ with $\Omega(0)=V$, which satisfies the 
assumptions of Corollary~\ref{c:41}  on some initial interval $t\in [0, T_1]$.

\begin{figure}[ht]
\centerline{\hbox{\includegraphics[width=14cm]{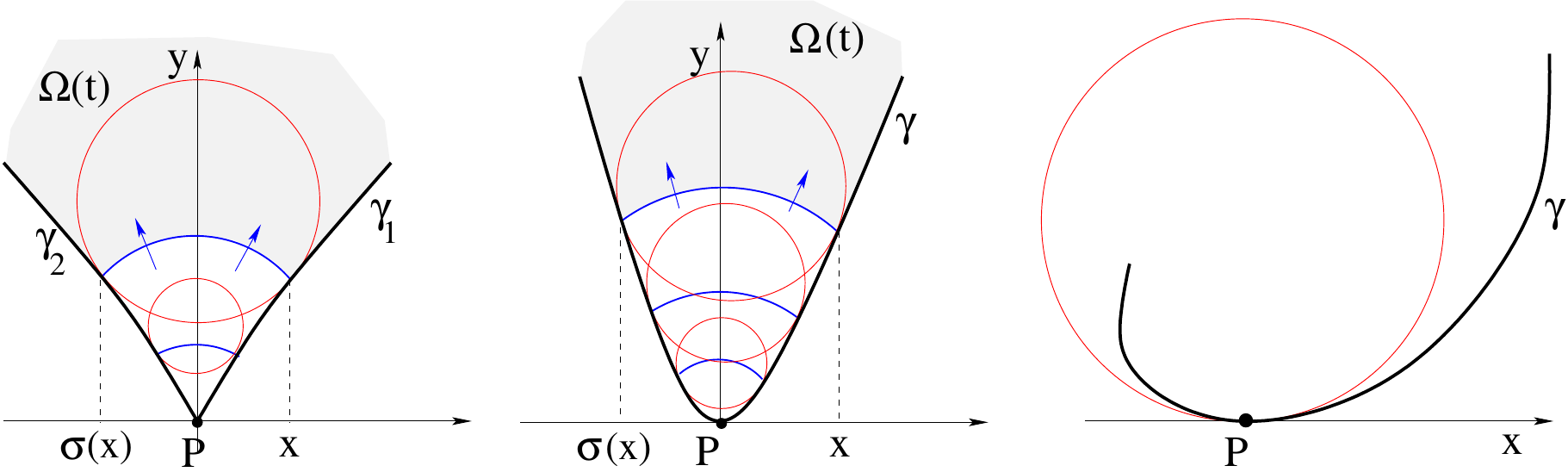}}}
\caption{\small  In order to construct a family of circumferences that cross perpendicularly
a curve $\gamma$, it suffices to construct a family of circumferences tangent to $\gamma$ at 
two distinct points. 
Left and center: this is possible near a point $P$ where $\gamma$ has a corner, or
where the curvature is maximum (or minimum). Right: this construction is not possible in case of a spiral-like curve $\gamma$, where the curvature is monotone increasing.}
\label{f:sm89}
\end{figure}

As shown in Fig.~\ref{f:sm89}, for $t>0$ small the complementary set $V\setminus\Omega(t)$ will be strictly
increasing, bounded by an arc of circumference which crosses the boundary $\partial V$
perpendicularly at both endpoints.

The next two propositions show that such a family of circumferences 
can be constructed in the neighborhood of a boundary point $P\in \partial V$ in two main cases, shown
in Fig.~\ref{f:sm89}, right and center:
\begi
\item[(i)] The boundary $\partial V$ has an outward corner at $P$. 
\item[(ii)] $P$ is a point where the curvature of the boundary $\partial V$ attains a local maximum.
\endi

\begin{proposition}\label{p:131}
Assume that the boundary $\partial V$ contains two $\C^2$ curves $s\mapsto \gamma_i(s)$,
$i=1,2$, 
joining at an angle
$\theta\in \,] 0, \pi[\,$ at a point $P= \gamma_1(0)=\gamma_2(0)$.   Then,  there exists 
$s_0>0$ and a family of circumferences
which cross perpendicularly the curves $\gamma_i$ at points $\gamma_1(s)$, $\gamma_2(\sigma(s))$, for all $s\in \,]0, s_0]$.
\end{proposition}

{\bf Proof.} {\bf 1.} As shown in Fig.~\ref{f:sm89}, left, using cartesian coordinates $x,y$
we can assume that $P=(0,0)$ and the 
curves  $\gamma_1$ and $\gamma_2$ are parametrized by 
$$\bega{rl} x\,\mapsto\, \gamma_1(x)\,=\,\bigl(x,\phi_1(x)\bigr), \qquad &\hbox{for} ~~x\in [0, x_0],\\[2mm] 
 x\,\mapsto\, \gamma_2(x)\,=\,\bigl(x,\phi_2(x)\bigr), \qquad &\hbox{for} ~~x\in [-x_0, 0],\enda$$   
for some $\C^2$ functions $\phi_1,\phi_2$ which satisfy
\bel{phiprop} \phi_1(0)\,=\,\phi_2(0)\,=\,0,\qquad\qquad \phi_1'(0)\,=\,- \phi_2'(0)>0.\eeq
To find a family of circumferences which cross perpendicularly the curves $\gamma_1, \gamma_2$, 
it clearly suffices to construct a family of circumferences which are tangent to both curves.

Toward this goal, for any $x>0$
small, consider a circle with radius $r>0$ tangent to $\gamma_1$ at the point $\gamma_1(x)$.   
We claim that, for a suitable choice of $r=r(x)$, this circle will be also tangent to $\gamma_2$ at some point 
$\gamma_2\bigl(\sigma(x)\bigr)$. 
\v
{\bf 2.} In the following,
tangent and normal vectors to the curves $\gamma_i$  at the points $\gamma_i(x)$, $i=1,2,$ 
are denoted by
$$\bft_1(x)=\big(1,\phi_1'(x)\big), \quad \bfn_1(x)=(-\phi_1'(x),1),
\quad \bft_2(x)=\big(1,\phi_2'(x)\big), \quad \bfn_2(x)=\big(-\phi_2'(x),1\big).$$
To uniquely determine the value of $\sigma(x)$, we  impose that the center of the circumference tangent to the two graphs lies
at the intersection of two lines: one parallel to $\bfn_1(x)$ through $(x, \phi_1(x))$ and the other 
parallel to to $\bfn_2\big(\sigma(x)\big)$ through $(\sigma, \phi_1(\sigma))$.
These can be parameterized as
$$t\,\mapsto\, \bigl(x, \phi_1(x)\bigr) + t \bigl(-\phi_1'(x),1),\qquad\quad s\,\mapsto\, \bigl(\sigma, \phi_2(\sigma)\bigr) + 
s \bigl(-\phi_2'(\sigma),1).$$
For any $x>0$ small, 
we thus seek a solution $t,s>0$, $\sigma<0$, of the system
\bel{eq-circ} 
\left\{\begin{array}{rl}
x-\phi'_1(x)t&=~\sigma-\phi'_2(\sigma)s,\\[2mm]
\phi_1(x)+t&=~\phi_2(\sigma)+s,\\[2mm]
\bigl((\phi_1'(x))^2+1\bigr)t^2&=~\bigl((\phi'_2(\sigma)\bigr)^2+1\big)s^2.
\end{array}\right.
\eeq
The above identities hold if and only if 
%
%
\bel{e33} \left\{\begin{array}{l}
\displaystyle s\,=\,\frac{x-\sigma+\phi_1'(x)\big(\phi_1(x)-\phi_2(\sigma)\big)}{\phi_1'(x)-\phi_2'(\sigma)},\\[10pt]
\displaystyle t\,=\,\frac{x-\sigma+\phi_2'(\sigma)\big(\phi_1(x)-\phi_2(\sigma)\big)}{\phi_1'(x)-\phi_2'(\sigma)},\\[10pt]
\sqrt{(\phi_1'(x))^2+1}\cdot t~=~
\sqrt{(\phi_2'(\sigma))^2+1}\cdot s.
\end{array}\right.\eeq
Inserting the two expressions for $t,s$ in the third equation, one obtains
$$
\sqrt{(\phi_1'(x))^2+1}\cdot \big(x-\sigma+\phi_2'(\sigma)\big(\phi_1(x)-\phi_2(\sigma)\big)~=~
\sqrt{(\phi_2'(\sigma))^2+1}\cdot\big(x-\sigma+\phi_1'(x)\big(\phi_1(x)-\phi_2(\sigma)\big).
$$
This holds if and only if
\begin{align*}
\Big(\sqrt{(\phi_1'(x))^2+1}-&\sqrt{(\phi_2'(\sigma))^2+1}\Big)\big(x-\sigma)\\
&+
\Big(\phi_2'(\sigma)\sqrt{(\phi_1'(x))^2+1}-\phi_1'(x)\sqrt{(\phi_2'(\sigma))^2+1}\Big)\big(\phi_1(x)-\phi_2(\sigma)\big)=0.
\end{align*}
Multiplying by $\Big(\sqrt{(\phi_1'(x))^2+1}+\sqrt{(\phi_2'(\sigma))^2+1}\Big)$
and dividing by $\big(\phi_1'(x)-\phi_2'(x)\big)$, we eventually obtain the equation
\bel{eq-F}\bega{l}
F(x,\sigma)~=~\big(\phi_1'(x)+\phi'_2(\sigma)\big)\big(x-\sigma\big)\\[2mm]
\qquad\qquad\ds+\big(\phi_1(x)-\phi_2(\sigma)\big)\Big(\phi'_1(x)\phi_2'(\sigma)-1-\sqrt{(\phi'_2(\sigma))^2+1}\sqrt{(\phi'_1(x))^2+1}\Big)~=~0.\enda
\eeq
By (\ref{phiprop}), at the origin
the partial derivatives are computed by 
{\small
$$\bega{rl}\ds {\partial F(x,\sigma)\over\partial x}\bigg|_{x=\sigma=0}&=~ \big(\phi_1'(0)+\phi'_2(0)\big)
+\phi_1'(0)\Big(\phi'_1(0)\phi_2'(0)-1-\sqrt{(\phi'_2(0))^2+1}\sqrt{(\phi'_1(0))^2+1}\Big)\\[4mm]
& \ds=~ \phi_1'(0)\Big( -\bigl( \phi_1'(0)\bigr)^2- 1 - \bigl[ (\phi_1'(0))^2 + 1\bigr]\Big)~<~0.
\enda
$$
$$\bega{rl}\ds {\partial F(x,\sigma)\over\partial\sigma}\bigg|_{x=\sigma=0}&=~- \big(\phi_1'(0)+\phi'_2(0)\big)
-\phi_2'(0)\Big(\phi'_1(0)\phi_2'(0)-1-\sqrt{(\phi'_2(0))^2+1}\sqrt{(\phi'_1(0))^2+1}\Big)\\[4mm]
& \ds=~ \phi_1'(0)\Big( -\bigl( \phi_1'(0)\bigr)^2- 1 - \bigl[ (\phi_1'(0))^2 + 1\bigr]\Big)~<~0.
\enda
$$}
By the implicit function theorem, the equation $F(x, \sigma)=0$ thus uniquely determines the function $\sigma(x)$,
for $x>0$ small. In fact,  $\sigma(0)=0$ and  ${d\over dx} \sigma(x)= -1$.
\endproof
%
%

\v
Next, we consider the case where the boundary $\partial V$ is smooth, and show that 
a family of circumferences 
can be constructed in a neighborhood of a point $P$ where the curvature of $\partial V$
attains a local maximum.

\begin{proposition}\label{p:132} Assume that a portion of the boundary $\partial V$ is a $\C^5$ curve $s\mapsto \gamma(s)$,
whose curvature $\omega(s)$ attains a positive, 
strict local maximum at a point $P= \gamma(\bar s)$. 
Namely, $\omega'(\bar s)=0$, $\omega''(\bar s)<0$.

Then there exists $\ve_0>0$ and a family of circumferences $\gamma^\ve$
with 
radii $r(\ve)$ depending continuously on $\ve\in [0,\ve_0]$, with the following properties.
Each arc $\gamma^\ve$ crosses the boundary $\partial V$
perpendicularly at two points $Q_1(\ve), Q_2(\ve)$.  Moreover, as $\ve\to 0$ one has the convergence $r(\ve)\to 0$ and 
$Q_1(\ve), Q_2(\ve)\to P$.
\end{proposition}

{\bf Proof.} Working in a cartesian coordinates, we parametrize the boundary 
$\gamma$ by $x\mapsto \bigl(x,\phi(x)\bigr)$.
Without loss of generality we can assume that $\bar s=0$, $\gamma(\bar s) = (0,0)$ and 
$\phi(0)=\phi'(0)=0$. Higher order derivatives will be denoted by  $ \phi^{(j)}(x)$, for $j=3,4,5$.
Notice that, since $\omega(x)=\tfrac{\phi''(x)}{\big(1+(\phi'(x))^2\big)^{3/2}}$, 
the assumption $\phi'(0)=0$ implies
\bel{omder} \omega(0)\,=\,\phi''(0)\,>\,0,\qquad \omega'(0)\,=\,\phi^{(3)}(0)\,=\,0,\qquad\omega''(0)=\phi^{(4)}(0)-3\big(\phi''(0)\big)^3\,<\,0.\eeq
Performing the same calculations as in the proof of  Proposition \ref{p:131}, to construct the family of circumferences 
we need to solve \eqref{eq-circ}, 
where now we simply define
$$\left\{\begin{array}{ll}
\phi_1(x) = \phi(x)& \text{ if } x\geq 0,\\
\phi_2(x) =\phi(x) & \text{ if } x<0.
\end{array}\right.$$ 
In view of (\ref{omder}) we have  the Taylor approximations
$$\phi(x)=\frac{\phi''(0)}2x^2+\frac{\phi^{(4)}(0)}{4!}x^4+\O(1)\cdot  x^5,\qquad 
\phi'(x)=\phi''(0)x+\frac{\phi^{(4)}(0)}{3!}x^3+\O(1)\cdot x^4,$$
The same expression for $F(x,\sigma)$ used at \eqref{eq-F} now yields
$$\bega{rl}
&\ds F(x,\sigma)=\ds\left(\phi''(0)\big(x+\sigma\big)+
{\phi^{(4)}(0)\over3!}\big(x^3+\sigma^3\big)+\O(1)\cdot (x^4+\sigma^4)\right)\big(x-\sigma\big)\\[4mm]
&\qquad \ds+\left(\frac{\phi''(0)}2\big(x^2-\sigma^2\big)+
\frac{\phi^{(4)}(0)}{4!}\big(x^4-\sigma^4\big)+\O(1)\cdot \bigl(|x|^5+|\sigma|^5\bigr)\right)\\[4mm]
&\ds\qquad\qquad \times
\left(-2+\frac{\phi''(0)^2}{2}\big(x-\sigma\big)^2+\O(1)\cdot (x^4+\sigma^4)\right)\\[4mm]
&\quad\ds=\frac{\phi^{(4)}(0)}{12}\Big(2\big(x^3+\sigma^3\big)\big(x-\sigma\big)-\big(x^4-\sigma^4\big)\Big)
-\frac{(\phi''(0))^3}4\big(x^2-\sigma^2\big)\big(x-\sigma\big)^2+R(x,\sigma)\\[4mm]
&\quad\ds=\frac1{12}\Big(\phi^{(4)}(0)-3(\phi''(0))^3\Big)\big(x-\sigma\big)^3\big(x+\sigma\big)+R(x,\sigma)\\[4mm]
&\quad\ds=\big(x-\sigma\big)^3\left(\frac1{12}\omega''(0)\big(x+\sigma\big)+\frac{R(x,\sigma)}{\big(x-\sigma\big)^3}\right).
\enda$$
Here the remainder $R$ is a $\C^5$ function that satisfies
$$R(x,\sigma)=\O(1)\cdot\bigl( |x|^5+|\sigma|^5\bigr).$$
The equation $F(x,\sigma)=0$ is equivalent to 
\bel{eqe}
\sigma~=~-x + {12 R(x,\sigma)\over \omega''(0)\cdot (x-\sigma)^3}~\doteq~G(x,\sigma)\,.\eeq
In the region where $0<x<1$, $\sigma\leq 0$, one has
$$G(x, -x)\,=\,-x+\O(1)\cdot x^2,
\qquad\qquad {\partial\over\partial\sigma} G(x,\sigma)~=~\O(1)\cdot \bigl( |x|+|\sigma|\bigr).$$
In this setting, $\sigma(x)$ is the fixed point of the map $\sigma\mapsto G(x,\sigma)$.
By the contraction mapping theorem,  a unique fixed point $\sigma(x)$ exists for all $x>0$ small enough.  Moreover, the map $x\mapsto \sigma(x)$ is $\C^5$ and 
satisfies 
$$\bigl| \sigma(x) + x\bigr|~=~\O(1)\cdot  x^2.$$
This establishes the existence of a family of circumferences tangent to the graph of $\phi$ at the points $x, \sigma(x)$.   In turn, this yields the existence of circumferences which
cross the graph of $\phi$ perpendicularly at the same points.
\endproof

\begin{remark} {\rm In Proposition~\ref{p:131}, one can still carry out the same construction
in the case where the angle satisfies $\theta\in \,]\pi, 2\pi[\,$.  
Such a strategy will satisfy the necessary conditions, but we do not expect
that it will be globally optimal.  Similarly, in
Proposition~\ref{p:132}, one can consider the case where the curvature has a local minimum.
The above construction will yield a strategy which satisfies the necessary conditions, but it will likely not be optimal.}
\end{remark}

\section{Maximally extended free arcs}
\label{s:14}
\setcounter{equation}{0}

This section is focused on the construction of maximally extended free boundary arcs.
Our goal is to derive a system of ODEs that describes these special curves.  
The  starting point is provided by the following remark.

\begin{figure}[ht]
\centerline{\hbox{\includegraphics[width=8cm]{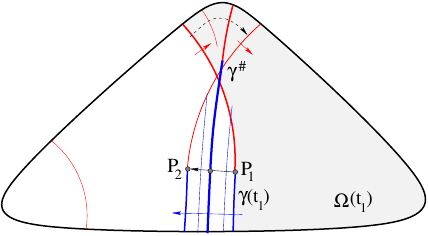}}}
\caption{\small The setting considered in Remark~\ref{r:141}.}
\label{f:sm91}
\end{figure}

\begin{remark}\label{r:141} {\rm
Fix $\xi$ and consider a trajectory $t\mapsto x(t,\xi)$, perpendicular to the 
boundary $\partial \Omega(t)$.  As shown in Fig.~\ref{f:sm91}, assume that 
\begi
\item[(i)] For $t\in [t_1, t_2]$, the point $x(t,\xi)$ lies on the free portion of the boundary,
where no control effort is present.
\item[(ii)] For $i=1,2$, the point  $P_i=x(t_i,\xi)$ lies at the edge of the arcs where the control is active.
\endi
The necessary conditions (\ref{max1}) imply 
\bel{nc5} Y(t_1, \xi)= Y^*(t_1),\qquad\qquad Y(t_2,\xi)= Y^*(t_2).\eeq
Calling  $\omega(t,\xi)$  the  curvature of the boundary 
$\partial\Omega(t)$ at the point $x(t,\xi)$, and calling 
$\omega^*(t)$  the  curvature of the portion of the boundary $\partial \Omega(t)$
where the control is active,  by (\ref{nc5}) and (\ref{Yeq}) it follows
\bel{nc6}
\int_{t_1}^{t_2} \omega(t,\xi)\, dt~=~\int_{t_1}^{t_2} \omega^*(t)\, dt.\eeq
}
\end{remark}

\begin{figure}[ht]
\centerline{\hbox{\includegraphics[width=6cm]{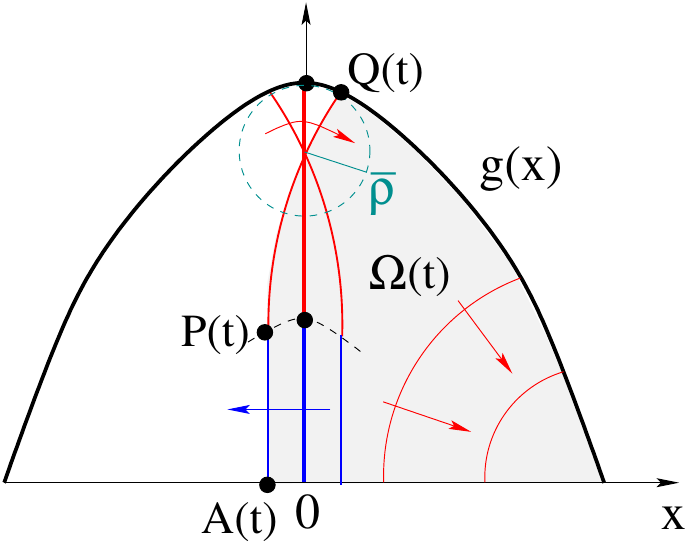}}\qquad
\hbox{\includegraphics[width=7cm]{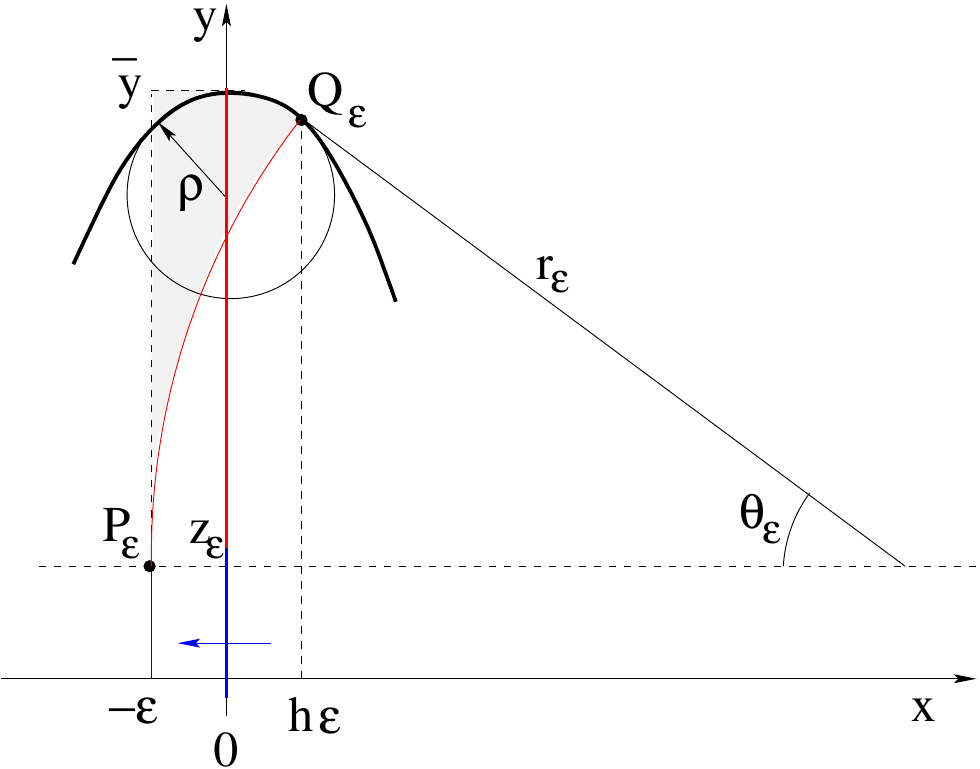}}}
\caption{\small  Left: the relative boundary $\partial \Omega(t)\cap V$ contains a free arc: the vertical segment with endpoints
$A(t), P(t)$, and a controlled arc: the arc of circumference with endpoints $P(t), Q(t)$.   Right: computing the active portion of the boundary,
at the time $t=0$ when this boundary is contained in the $y$-axis.
 }
\label{f:sm92}
\end{figure}

\subsection{The symmetric case.}\label{s:51}
As a first step, we consider 
the case where the domain $V$ is symmetric w.r.t.~the $y$-axis. More precisely (see Fig.~\ref{f:sm92}, left)
\bel{Vdef} 
V\,=\,\bigl\{ (x,y)\,;~~0<y<g(x)\bigr\},\eeq
for some $\C^2$ function $g:\R\mapsto\R$ such that 
$$g(x)= g(-x), \qquad g(0)>0,\qquad g''(x)<0\qquad\forall x\in \R.$$
In this setting, a natural guess is that the maximal free boundary $\ov\gamma$ 
should be a vertical segment 
contained in the axis of symmetry.   
For notational convenience we shift time, so that at time $t=0$ the relative boundary is contained in the $y$-axis:
$\Omega(0)= \bigl\{ (x,y)\,;~~0<y<g(x)\,,~~x>0\bigr\}$.

As a consequence, as shown in Fig.~\ref{f:sm92}, left, at each time $t$ the relative boundary $\partial\Omega(t)\cap V$ will be the union of a free portion: the vertical segment with endpoints $A(t), P(t)$, and a controlled arc: the arc of circumference with endpoints $P(t), Q(t)$. Notice that, for a fixed $t$, there is a 1-parameter family
of circumferences that are tangent to the vertical line $\{x=-t\}$ and cross perpendicularly the boundary $\partial V$.
A unique choice of this circumference is determined by the area identity (\ref{dareat}).

We observe that, in this symmetric setting, the identity (\ref{nc6}) is trivially satisfied.   Indeed,
the curvature of the free arc is $\omega(t,\xi)\equiv 0$.   On the other hand, by symmetry we have
$t_1(\xi)=-t_2(\xi)$, while $\omega^*(-t) = -\omega^*(t)$ for all $t$.  Hence the right hand side of (\ref{nc6}) 
vanishes as well.

Next, we wish to determine the lengths of the arcs $AP$ and $PQ$ at time $t=0$, when they are both vertical segments on the $y$-axis.
These should depend on:
\begi
\item The radius of curvature $\ov \rho$ of the boundary $\partial V$ at the point $Q$.
\item The bound $M$ on the total effort.
\endi
With reference to Fig.~\ref{f:sm92}, right, for any 
$\ve>0$, consider the point $Q_\ve= (x_\ve, y_\ve)$ where a circumference tangent to the vertical line $\{x=-\ve\}$
crosses the boundary $\partial V$ perpendicularly.  We denote by 
$P_\ve= (-\ve, z_\ve)$ the  point where this circumference has a vertical tangent.
 
As $\ve\to 0$, the area identity (\ref{dareat})  yields 
$$[\hbox{area of shaded region}] ~=~{1\over 2} (y_\ve-z_\ve) \cdot (\ve + x_\ve)+o(\ve)~=~ M\ve + o(\ve)$$
We now set $x_\ve = h\ve+ o(\ve)$ and denote by 
$\ell\doteq\lim_{\ve\to 0} (y_\ve-z_\ve)$ the length of the free arc at time $t=0$.
Letting $\ve\to 0$, we obtain the  identity
\bel{Meq}   {\ell\over 2}\, (1+h)~=~M.\eeq
Next, we impose the orthogonality condition at the boundary.  Setting $\bar y= g(0)$ and
calling $ \rho= -1/g''(0) >0$ the radius of curvature of the boundary $\partial V$ at the point $Q= \bigl(0, g(0)\bigr)$, 
recalling that $g$ is an even function, a Taylor expansion yields
$$y~=~ g(x)~=~\bar y - {1\over 2 \rho} x^2+o(x^3).$$
Referring to Fig.~\ref{f:sm92}, right, 
call $r_\ve$ the radius of the circumference through $P_\ve $ and $Q_\ve$, and $\theta_\ve>0$ the angle of the corresponding arc between $Q_\ve$ and $ P_\ve$.
By standard trigonometric identities we find
\bel{theq}r_\ve\sin \theta_\ve~=~\ell - {h^2\ve^2\over 2\rho} + o(\ve^2),\qquad\qquad r_\ve-r_\ve \cos\theta_\ve~=~(1+h)\ve +o(\ve),\eeq
\bel{id1}{1-\cos\theta_\ve\over\sin\theta_\ve}\,\ell~=~(1+h)\ve + o(\ve).\eeq
On the other hand, the orthogonality condition at $Q_\ve$ implies
\bel{id2}\theta_\ve~=~{h\ve\over\rho}+o(\ve).\eeq
Combining the two  identities (\ref{id1})-(\ref{id2}) one obtains
$${\theta_\ve^2/2\over \theta_\ve}\,\ell~=~ {\theta_\ve\over 2}\ell~=~{h\ve\over 2\rho}\,\ell+o(\ve)~=~(1+h)\ve + o(\ve).$$
Letting $\ve\to 0$ we conclude
\bel{eleq}h\ell~=~2\rho (1+h).\eeq
Together with the first identity (\ref{Meq}), this determines the two constants $\ell, h$:
\bel{sol} h\,=\,{2M\over \ell}-1,\qquad\qquad \ell\,=\,M- \sqrt{M^2 -4\rho M}\,.\eeq
Notice that the sign in front of the square root is consistent with the fact that, as the radius $\rho\to 0$, the length of the controlled arc also approaches zero:
$$\ell ~=~M\left(1-\sqrt{1-{4\rho\over M}}\right)~\approx~ 2\rho~\to ~0.$$

%
\subsection{Maximal free arcs: the non-symmetric case.}
\begin{figure}[ht]
\centerline{\hbox{\includegraphics[width=11cm]{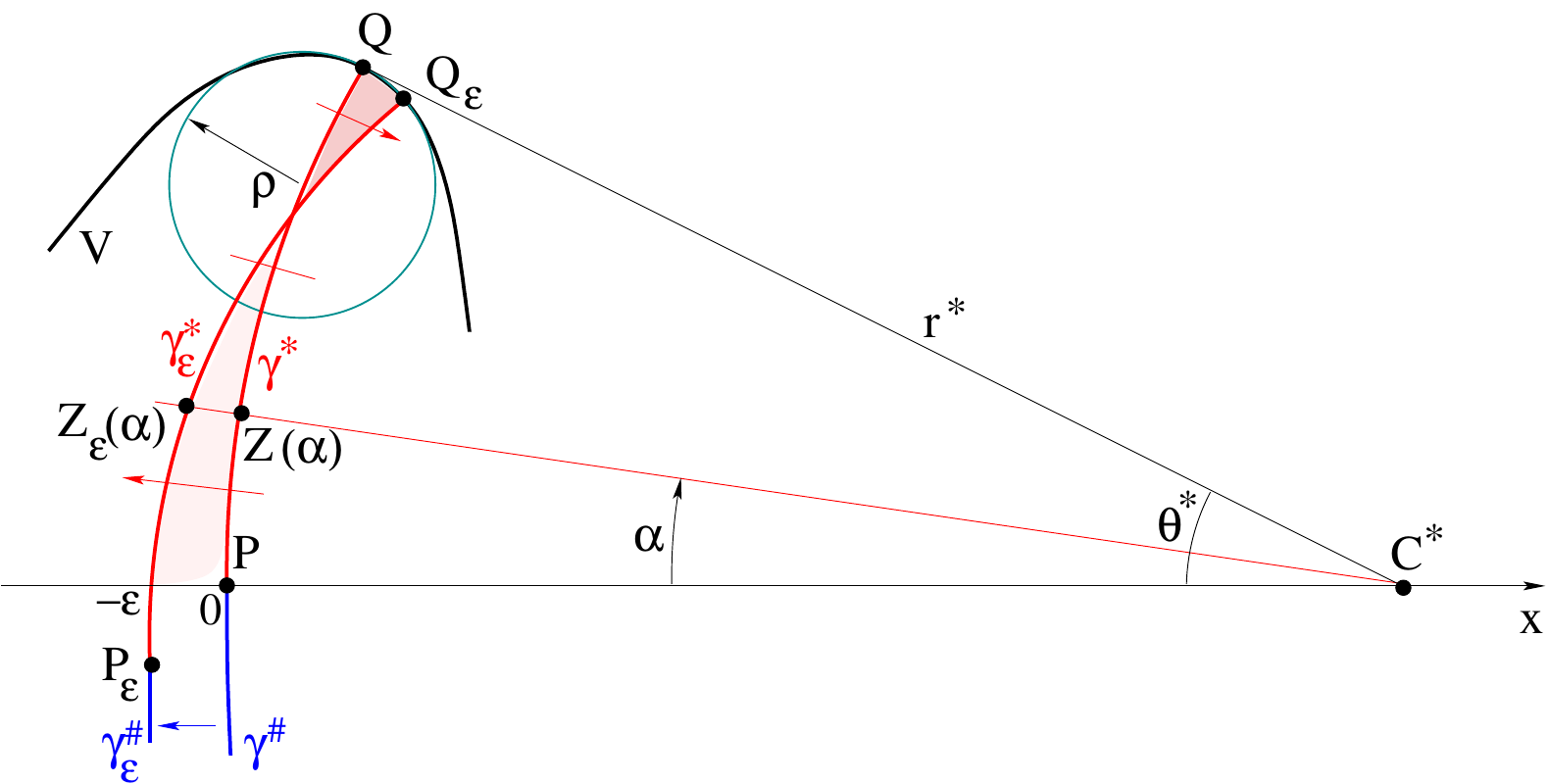}}}
\caption{\small   The shaded region yields the signed area swept by the moving arc of
circumference for $t\in [0,\ve]$, up to higher order infinitesimals.
 }
\label{f:sm93}
\end{figure}

If the set $V$ is not symmetric, we do not expect that a maximal free arc  $ \gamma^\sharp$ should be a straight line.
We seek a system of ODEs describing this arc. This will be achieved in several steps.
\v
{\bf 1.} Consider the configuration in Fig.~\ref{f:sm93}.
At time $t=0$ the relative boundary 
$\partial \Omega(0)\cap V= \gamma^*\cup \gamma^\sharp$, is the union of a 
controlled arc $\gamma^*$ and a free arc 
$\gamma^\sharp$  where the set $\Omega(t)$ expands with unit speed.
   By a suitable choice of coordinates, 
we assume that the arc of circumference $\gamma^*$ is tangent to $\gamma^\sharp$ at 
the point $P=(0,0)$, and perpendicular to $\partial V$ at the other endpoint $Q$.
We choose the $x$-axis so that it is perpendicular to $\gamma^\sharp$ and $\gamma^*$
at $P$. 

Let $\gamma^\sharp$ be parameterized in terms of arc-length, as $\xi\mapsto \gamma^\sharp(\xi)$, with $\gamma^\sharp(0)=P$,
and call $\bfn(\xi)$ the unit inward normal vector.
Then at time $t=\ve>0$ the corresponding free arc will be parameterized by 
\bel{gesharp} \gamma^\sharp_\ve(\xi)~=~ \gamma^\sharp(\xi) - \ve \bfn(\xi).\eeq
In addition, at time $\ve$ the controlled arc $\gamma_\ve^*$ will be a portion of circumference with endpoints
 $Q_\ve\in \partial V$ and $P_\ve =\gamma_\ve^*\bigl(\xi(\ve)\bigr)$, for some $\xi(\ve)$.
 We seek a formula for the time derivatives 
 \bel{dotPQ}\dot P\,=\,\lim_{\ve\to 0+} \,{P_\ve-P\over\ve}, \qquad\qquad \dot Q\,=\, \lim_{\ve\to 0+} \,{Q_\ve-Q\over\ve}\,.\eeq
 This will rely on three properties:
  \begi
\item The junction at $P_\ve$ is tangential.
\item The intersection at $Q_\ve$ is perpendicular.
\item By the  identity (\ref{dareat}),  the signed area swept by the moving arc of circumference 
equals $M - \H^1(\gamma^*)$.
 \endi

{\bf 2.} 
To fix ideas, assume that at time $t=0$ the controlled arc $\gamma^*$ has center at 
$C^*= (r^*,0)$,  radius $r^*$ and spans  
an angle $\theta^*$. 
As shown in Fig.~\ref{f:sm93}, for every angle $\alpha\in [0, \theta^*]$ let $Z(\alpha)$
and $Z_\ve(\alpha)$ be the intersections of the circumferences $\gamma^*,\gamma^*_\ve$
with the line through $C^*$, forming an angle $\alpha$ with the $x$-axis, and call 
$r^*_\ve(\alpha) = \bigl|Z_\ve(\alpha)-C^*\bigr|$.   Notice that $\bigl|Z(\alpha)-C^*\bigr|=r^*$
for all $\alpha$.

The signed area swept  (in the inward direction) by the moving arc of circumference during the time interval $[0,\ve]$ is computed by
\bel{Aep}A_\ve~=~\int_0^{\theta^*} \bigl( r^*_\ve(\alpha) - r^*\bigr) \, r^*d\alpha+o(\ve)    ~=~\ve\bigl(|\dot Q|-1\bigr){\theta^*r^*\over 2}+ o(\ve).\eeq
Indeed, the endpoint $P_\ve$ moves outward with unit speed, while the other endpoint
$Q_\ve$ moves inward with speed $\dot Q$.

The area identity (\ref{dareat}) now yields
\bel{Ae2} A_\ve~=~\ve\bigl[ M - \H^1(\gamma^*)\bigr] +o(\ve)~=~\ve\bigl[M- \theta^*r^*\bigr]+o(\ve).\eeq
Combining the two above formulas and letting $\ve\to 0$ we obtain 
\bel{Ae3} M-\theta^*r^*~=~\bigl(|\dot Q|-1\bigr){\theta^*r^*\over 2}\,,
\qquad \qquad |\dot Q|~=~{2M\over  \theta^*r^*}-1\,.\eeq

\begin{figure}[ht]
\centerline{\hbox{\includegraphics[width=11cm]{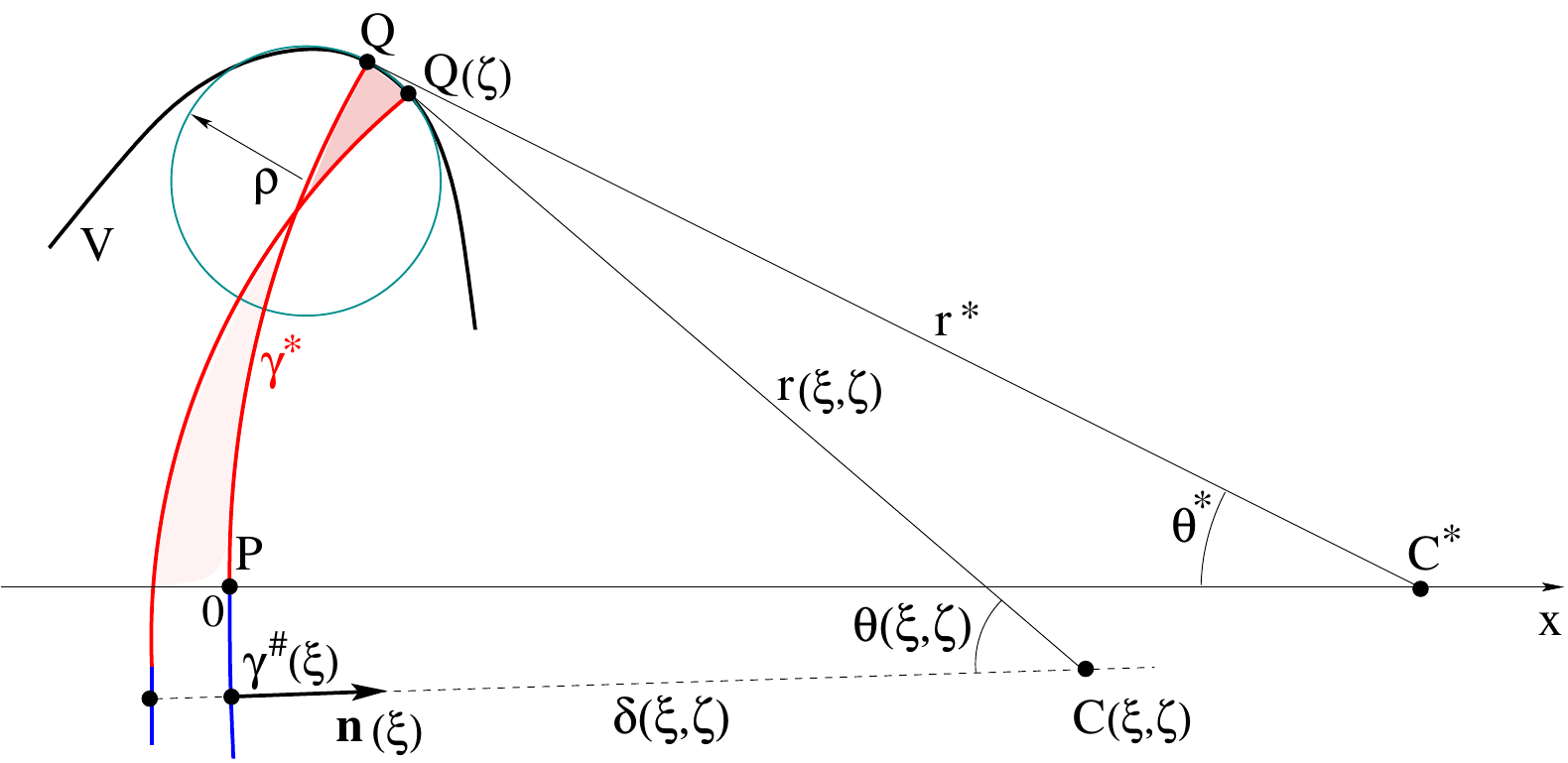}}}
\caption{\small   The shaded region is the area swept by the moving arc of
circumference for $t\in [0,\ve]$.
 }
\label{f:sm94}
\end{figure}

{\bf 3.} Next, referring to Fig.\ref{f:sm94}, let the boundary $\partial V$ be parameterized by 
arc-length as $\zeta\mapsto Q(\zeta)$, with $Q(0)=Q$.
For any pair $(\xi,\zeta)\approx(0,0)$, denote by $C(\xi,\zeta)$ the intersection of
\begi\item the perpendicular line to $\gamma^\sharp$ at the point $\gamma^\sharp(\xi)$,
\item the tangent line to $\partial V$ at the point $Q(\zeta)$.
\endi
The distances between $C(\xi,\zeta)$ and the points $\gamma^\sharp(\xi)$, $Q(\zeta)$
will be denoted respectively by
\bel{dist} \delta(\xi,\zeta) \,=\,\bigl| C(\xi,\zeta)-\gamma^\sharp(\xi)\bigr|,\qquad\qquad
r(\xi,\zeta) \,=\,\bigl| C(\xi,\zeta)-Q(\zeta)\bigr|.\eeq
When $\xi=\zeta=0$ we clearly 
have
\bel{dist*} \delta(0,0) \,=\, |C^*-P|\,=\,r^*\,=\,\bigl| C^*-Q|\,=\,r(0,0).\eeq

\begin{figure}[ht]
\centerline{\hbox{\includegraphics[width=10cm]{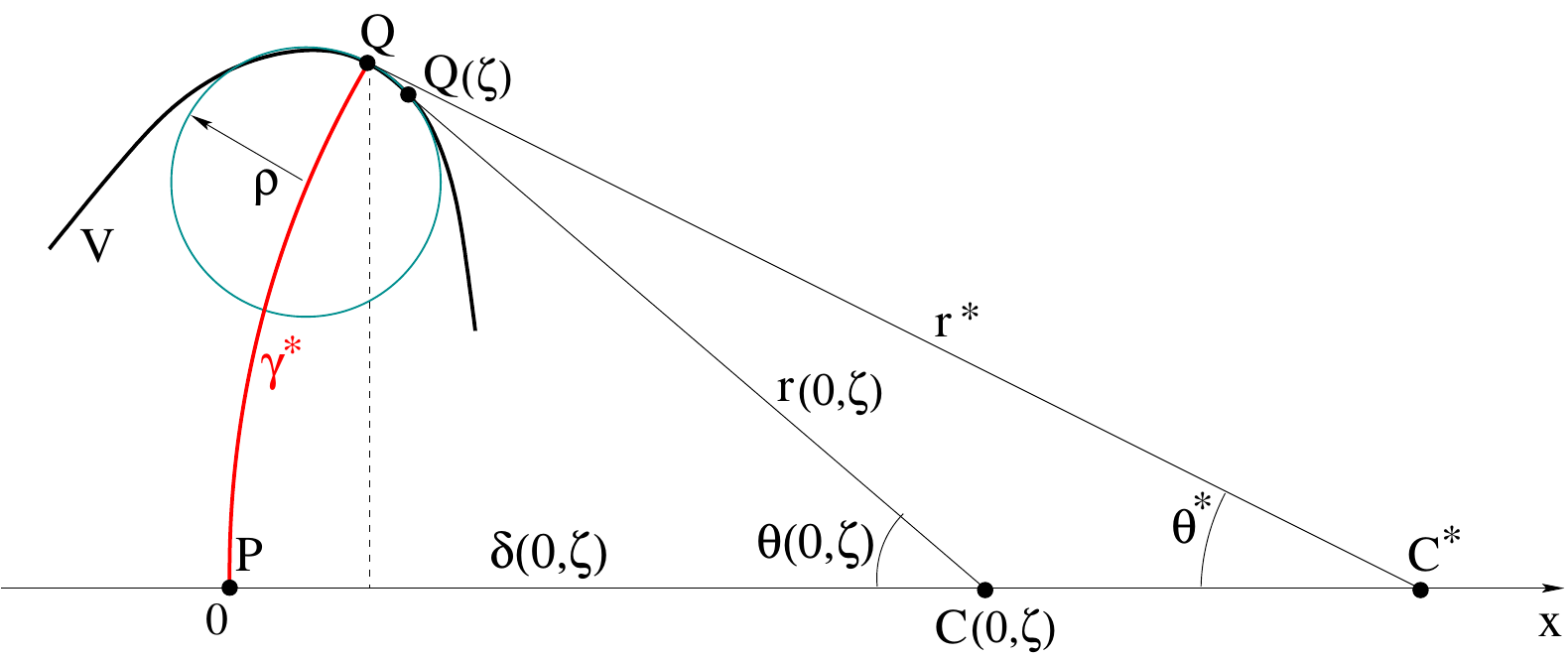}}}
\caption{\small  Computing the partial derivatives of the distances $r(\xi,\zeta)$ and 
$\delta(\xi,\zeta)$ w.r.t.~$\zeta$. }
\label{f:sm95}
\end{figure}

To compute the derivatives of these distances w.r.t.~$\xi$ and $\zeta$,
we let $\omega= 1/ \rho$ be the curvature of the boundary $\partial V$ at the point $Q$, 
and let $\omega^\sharp$ be the curvature of $\gamma^\sharp$ at the point $P$.
At the point $C(\xi,\zeta)$, the angle formed by the lines through $\gamma^\sharp(\xi)$ and 
$Q(\zeta)$ is thus
\bel{thxz}\theta(\xi,\zeta)~=~\theta^* + \omega^\sharp \xi+ \omega \zeta +o(\xi)+o(\zeta).\eeq
At $\xi=\zeta=0$, recalling that the point $Q(\zeta)$ is parameterized by arc-length, by elementary trigonometric identities (see Fig.~\ref{f:sm95}) one finds
\bel{prz}
{\partial \over \partial \zeta} r(0,\zeta)~=~- r^*\omega\,\cot \theta^*-1.\eeq
Moreover, observing that
$$\bigl|Q-C(0,\zeta)\bigr| \,\sin \theta(0,\zeta)~=~r^*\sin \theta^* + o(\zeta),\qquad \qquad 
{\partial\over\partial \zeta} \theta(0,\zeta)~=~\omega,$$
at $\xi=\zeta=0$ we obtain
\bel{pdz} {\partial\over\partial \zeta} \delta(0,\zeta)~
=~{\partial\over\partial \zeta} \bigl[ r^*\sin \theta^*\cot\theta(0,\zeta)\bigr]
~=~ r^*\sin \theta^*\cdot {-\omega\over\sin^2\theta^*}~=~- {r^* \omega\over \sin\theta^*}\,.
\eeq

\begin{figure}[ht]
\centerline{\hbox{\includegraphics[width=10cm]{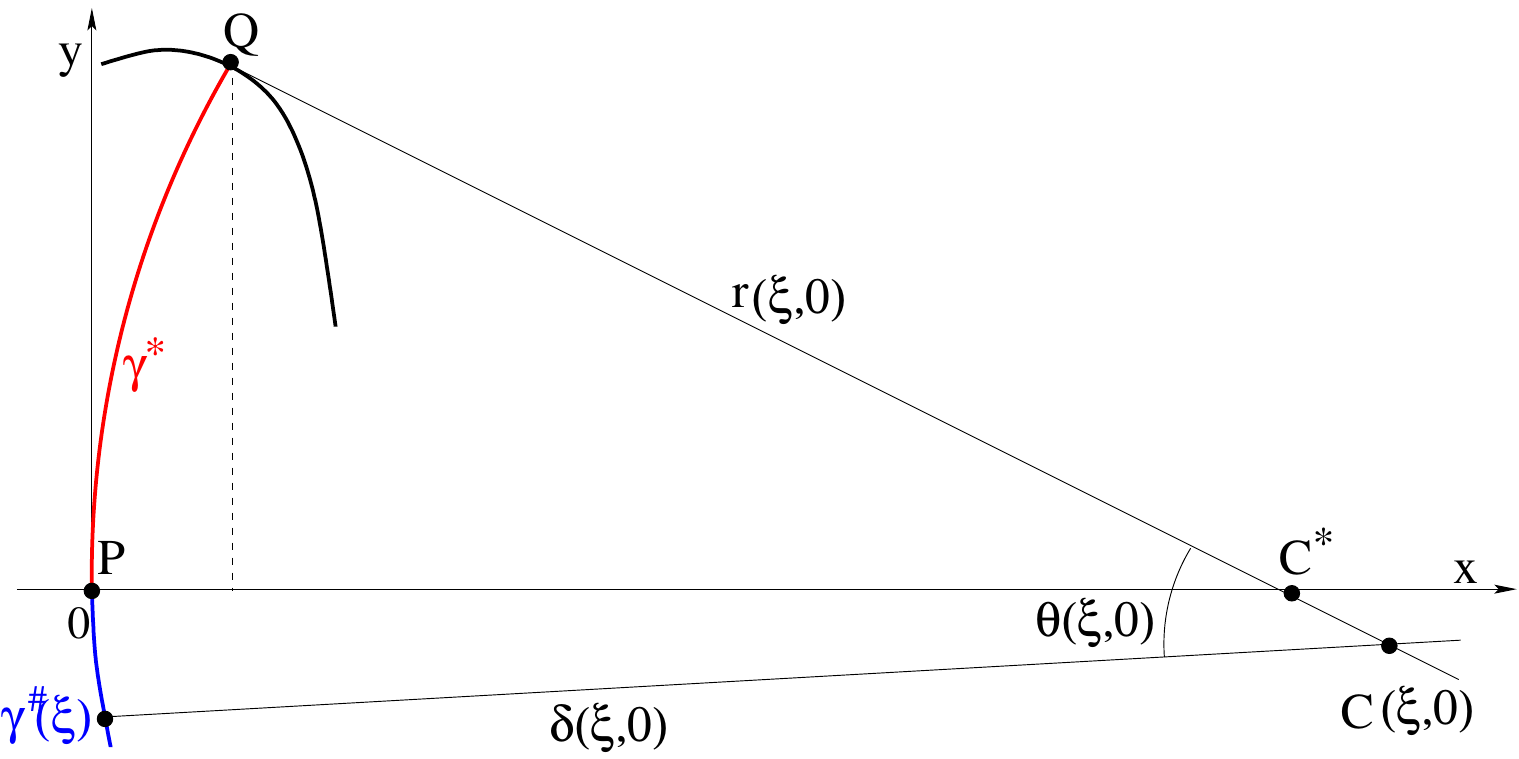}}}
\caption{\small  Computing the partial derivatives of the distances $r(\xi,\zeta)$ and 
$\delta(\xi,\zeta)$ w.r.t.~$\xi$. 
Here $P=(0,0)$, $C^*=(r^*,0)$. }
\label{f:sm96}
\end{figure}

Next, referring to Fig.~\ref{f:sm96}, we observe that the point $C(\xi,0)$ lies at the intersection
of two straight lines with equations
\bel{lines} y\,=\,- \tan\theta^*\,(x-r^*),\qquad\qquad y\,=\, -\xi + \omega^\sharp \xi \,x +o(\xi)\cdot x.\eeq
%
%
At the intersection point, this yields
$$\tan\theta^* (x-r^*) + (\omega^\sharp r^*-1)\xi~=~\omega^\sharp \xi (x-r^*) +o(\xi)\cdot x.$$
Observing that $x-r^*= \O(1)\cdot \xi$, we obtain
$C(\xi,0)~=~ \bigl(x(\xi), y(\xi)\bigr)$, 
where
\bel{xyx}
x(\xi)\,=\, r^* + {1-\omega^\sharp r^*\over\tan \theta^*} \,\xi +o(\xi),\qquad\qquad 
y(\xi)\,=\, (\omega^\sharp r^*-1)\,\xi+o(\xi).\eeq
Using (\ref{xyx}), the partial derivatives of $\delta$ and $r$ w.r.t.~$\xi$ at $\xi=\zeta=0$ are  computed by
\bel{derde}{\partial \over \partial \xi}  \delta(\xi,0)~=~{\partial \over \partial \xi}x(\xi)~=~  
{1-\omega^\sharp r^*\over\tan \theta^*} \,,\qquad\qquad 
{\partial \over \partial \xi} r(\xi,0)~=~{1\over \cos\theta^*} \,{\partial \over \partial \xi}x(\xi)~=~  
{1-\omega^\sharp r^*\over\sin\theta^*}\,.\eeq
 
{\bf 4.} 
For $t\approx 0$, denote by
$$Q(t)~=~Q\bigl(\zeta(t)\bigr), \qquad\qquad P(t)~=~\gamma^\sharp\bigl(\xi(t)\bigr)
- t\, \bfn\bigl(\xi(t)\bigr),$$ the endpoints
of the controlled arc of circumference $\gamma^*(t)$  at time $t$.
By the previous analysis we have
$${d\over dt}\zeta(t)\bigg|_{t=0} ~=~|\dot Q|~=~{2M\over  \theta^*r^*}-1.$$
We now use the identity
\bel{r=d}r\bigl(\xi(t), \zeta(t)\bigr) \,=\, t+\delta\bigl(\xi(t), \zeta(t)\bigr)\eeq
which is valid for all $t$, to compute the derivative
$d\xi(t)/dt $.
Inserting (\ref{prz}), (\ref{pdz}) and (\ref{derde}) in (\ref{r=d}),  we obtain
$${\partial r\over\partial \xi} {d\xi\over dt} +  {\partial r\over\partial \zeta} {d\zeta\over dt} 
~=~1+ {\partial \delta \over\partial \xi} {d\xi\over dt}+ {\partial \delta \over\partial \zeta} {d\zeta\over dt}\,,$$
\bel{dxep}\bega{rl}\ds
{d\over dt}\xi(t)\bigg|_{t=0}&\ds=~\left[ {\partial r\over\partial \xi}-{\partial \delta\over\partial \xi}
\right]^{-1} \left( 1+ {\partial \delta \over\partial \zeta} {d\zeta\over dt} -  {\partial r\over\partial \zeta}{d\zeta\over dt} \right)\\[4mm]
&\ds=~\left[(1-\omega^\sharp r^*) {1-\cos \theta^*\over\sin\theta^*} \right]^{-1}\cdot
\left[1 + \left( -{r^*\omega\over \sin\theta^*} + r^*\omega\cot\theta^*+1\right) \left( {2M\over \theta^*r^*} -1\right)\right].
\enda
\eeq
\v
{\bf 5.} We are now ready to derive a  second order system of ODEs satisfied by the maximal free arc $\gamma^\sharp$. As before, let 
 $\xi\mapsto \gamma^\sharp(\xi)$ be an arc-length parameterization, with $\gamma^\sharp(0) = P$.
For any $\xi>0$, let $t_1(\xi)<0<t_2(\xi)$  be the times where the points 
$$P_i(\xi)~=~\gamma^\sharp(\xi) - \bigl(t_i(\xi)-\tau\bigr)\bfn(\xi)$$ lie at the junction between the active and the free arcs (see Fig.~\ref{f:sm91}).

To derive an ODE for $\gamma^\sharp$, let $\omega^\sharp(\xi)$ be the curvature of $\gamma^\sharp$ at the point $\gamma^\sharp(\xi)$, and denote by $\omega(t)$ the curvature of the controlled arc at time $t$.
Differentiating (\ref{nc6}) w.r.t.~$\xi$ we obtain
%
\bel{curveq}
{\omega^{\sharp}(\xi)\over 1+ t_2\omega^{\sharp}(\xi)} {dt_2(\xi)\over d\xi} - {\omega^{\sharp}(\xi)\over 1+ t_1\omega^{\sharp}(\xi)} {dt_1(\xi)\over d\xi}~=~{\omega\bigl(t_2(\xi) \bigr)}\, {dt_2(\xi)\over d\xi} -{\omega\bigl(t_1(\xi) \bigr)}\, {dt_1(\xi)\over d\xi}\,.\eeq
Here the time derivatives $dt_i(\xi)/d\xi$ can be computed by inverting (\ref{dxep}).
Notice that this is a second order, highly nonlinear ODE.
Indeed, this is an implicit  equation that can be solved for the curvature $\omega(\xi)$ in terms of the other variables.
\v
\begin{remark} {\rm At a time $\tau$ when the free arc is maximal, assume that 
the boundary $\partial \Omega(\tau)\cap V$
is the union of the free arc $\gamma^\sharp$ together with  a controlled arc of circumference $\gamma^*$ with endpoints $P,Q$.
From the necessary conditions for optimality we deduce:
\begi
\item[(i)] The arc of circumference $\gamma^*$ is {\bf perpendicular} to the boundary 
$\partial V$ at the point $Q$.
\item[(ii)] The two arcs $\gamma^\sharp$ and $\gamma^*$  have a {\bf second order tangency} at the point $P$.
\endi
Indeed, (i) follows from Theorem~\ref{t:111}.   The first order tangency condition in (ii) follows from Theorem~\ref{t:101}.
Both of this facts remain valid at all times. On the other hand, the  
second order tangency is a unique property of maximal arcs.  
Indeed, let 
$$P(t,\xi)~=~\gamma^\sharp(\xi) -(t-\tau) \bfn(\xi)$$
denote a point on the free arc for $t\in \bigl[t_1(\xi), \, t_2(\xi)\bigr]$, 
which coincides with the endpoint of the controlled arc for $t=t_1(\xi)$ and $t=t_2(\xi)$, with 
$$t_1(\xi)<\tau<t_2(\xi),\qquad\qquad \lim_{\xi\to \bar \xi-} t_1(\xi)~=~\lim_{\xi\to \bar \xi-} t_1(\xi)~=~\tau.$$
By (\ref{nc6}) the curvatures are related by
$${1\over  t_2(\xi)- t_1(\xi)} \cdot \int_{t_1(\xi) }^{t_2(\xi)} \omega(t,\xi)\, dt~=~{1\over  t_2(\xi)- t_1(\xi)} \cdot \int_{t_1(\xi) }^{t_2(\xi)} \omega^*(t)\, dt.$$
Taking the limit as $\xi\to \bar \xi-$ one obtains (ii).
}\end{remark}
%
%
%
\section{The free interface near a corner point}
\label{s:15}
\setcounter{equation}{0}

\begin{figure}[ht]
\centerline{\hbox{\includegraphics[width=9cm]{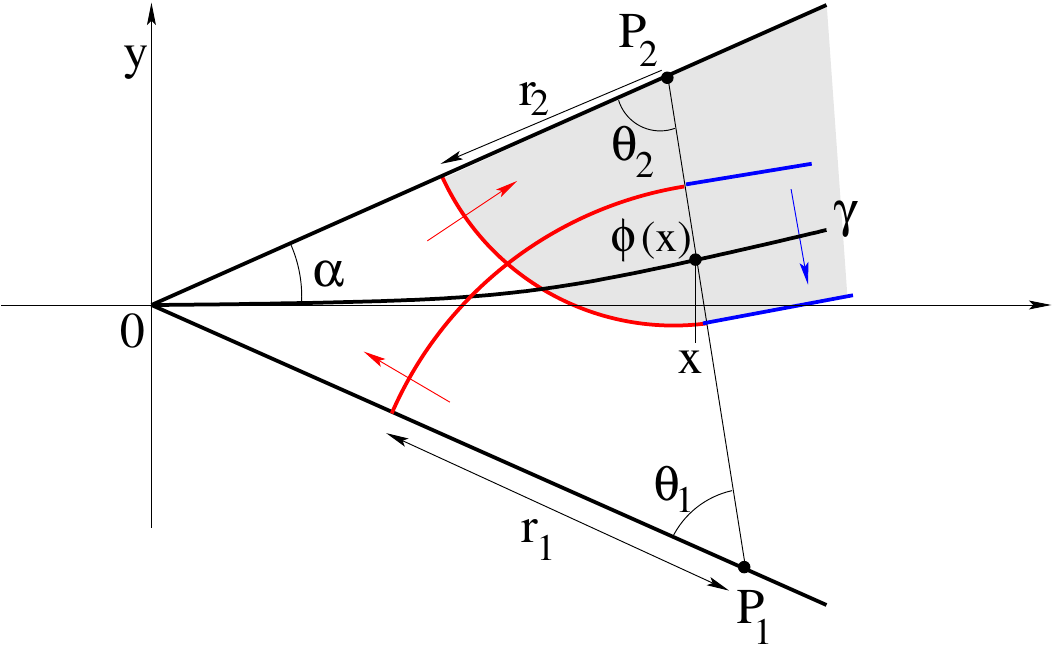}}}
\caption{\small Evolution of the free and the controlled portions of the 
boundary $\partial\Omega(t)$. Here $y=\phi(x)$ is the equation of the maximal free interface.}
\label{f:tg29}
\end{figure}

In the $x$-$y$ plane let the motion take place within a set $V$ having a corner point.
 After a change of coordinates,
as shown in Fig.~\ref{f:tg29} we assume
\bel{V1}
V~=~\bigl\{ (x,y)\,;~ |y|\leq   x \tan \alpha \bigr\},\eeq
for some $0<\alpha<\pi/2$.
%
Our goal is to describe the maximal free interfaces.

To fix ideas, assume that at time $t=0$ a smooth free interface 
\bel{fi0} y\,=\, \phi(x)\eeq
is given.   Then for small times $t\in [-\delta, \delta]$  we can uniquely
determine a family a circumferences, tangent to the free curves at points
$\bigl(x(t), \, \phi(x(t))\bigr)$.   Indeed, the radii of these circumferences will be determined
by the tangency requirement, together with the swept area equation (\ref{dareat}).  
In addition, the optimality conditions require that
the identities (\ref{nc6}) be satisfied.

Our goal is to show that 
there exists a 1-parameter family of admissible maximal interfaces, all satisfying
$\phi(0)=\phi'(0)=0$.

As shown in Fig.~\ref{f:tg29}, 
let $x\mapsto \gamma(x)=(x,\phi(x))$ be a parameterization of the
maximal free interface $\gamma$, with $\phi(0)=0$.  
For each $x>0$, call $P_1(x)$, $P_2(x)$ the points where the 
straight line perpendicular to $\gamma$ at the point $\gamma(x)=(x,\phi(x))$
crosses the lower and the upper boundary of
$V$, respectively.  Moreover, we call $\alpha$ the angle between the $x$-axis 
and the upper boundary of $V$, and set
$\beta=\frac{\pi}{2}-\alpha$

\begin{figure}[ht]
\centerline{\hbox{\includegraphics[width=9cm]{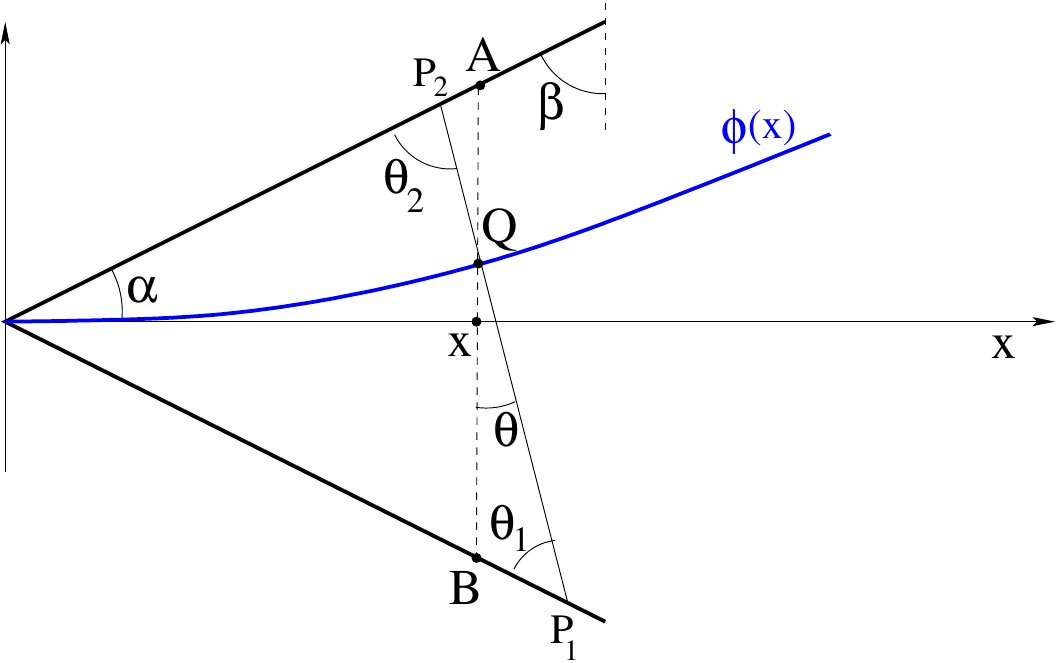}}}
\caption{\small  Deriving the identities in (\ref{lsine}).}
\label{f:tg30}
\end{figure}

Call
\bel{the}\theta(x)~=~\arctan \phi'(x),\eeq
so that 
\bel{cos} \sqrt{1+[\phi'(x)]^2} ~=~{1\over\cos\theta(x)}\,.\eeq
 The angles at $P_1(x)$ and $P_2(x)$ are then computed by 
\bel{the12}\left\{ \bega{rl}\theta_1(x)&=~\beta-\theta(x),\\[2mm]
\theta_2(x)&=~\beta+\theta(x).\enda\right.\eeq
With reference to Fig.~\ref{f:tg30}, according to the law of sines one has
\bel{lsine}{|P_1-Q|\over \sin (\pi-\beta)}~=~
{|B-Q|\over \sin(\beta-\theta)}\,,\qquad\qquad {|P_2-Q|\over \sin \beta}~=~
{|A-Q|\over \sin(\pi-\beta-\theta)}\,.
\eeq
Let $t_1=t_1(x)<0$ be the time when the center of the active circumference is at $P_1(x)$, and call $r_1= r_1(x)$ the radius of this circumference. Similarly, let $t_2=t_2(x)>0$ be the time when the center of the active circumference is at $P_2(x)$, and let $r_2= r_2(x)$
be the radius.
Using (\ref{lsine}) we obtain
\bel{r12}\left\{
\bega{rl}r_1(x, t_1)&\ds=~|P_1-Q| - t_1(x)~=~{\sin\beta\bigl(x\tan \alpha +\phi(x)\bigr) \over \sin
\bigl(\beta-\theta(x)\bigr)} - t_1,
\\[4mm]
r_2(x, t_2)&\ds=~|P_2-Q| + t_2(x)~=~{\sin\beta\bigl(x\tan \alpha -\phi(x)\bigr) \over \sin
\bigl(\beta+\theta(x)\bigr)} + t_2.\enda\right.\eeq

\begin{figure}[ht]
\centerline{\hbox{\includegraphics[width=9cm]{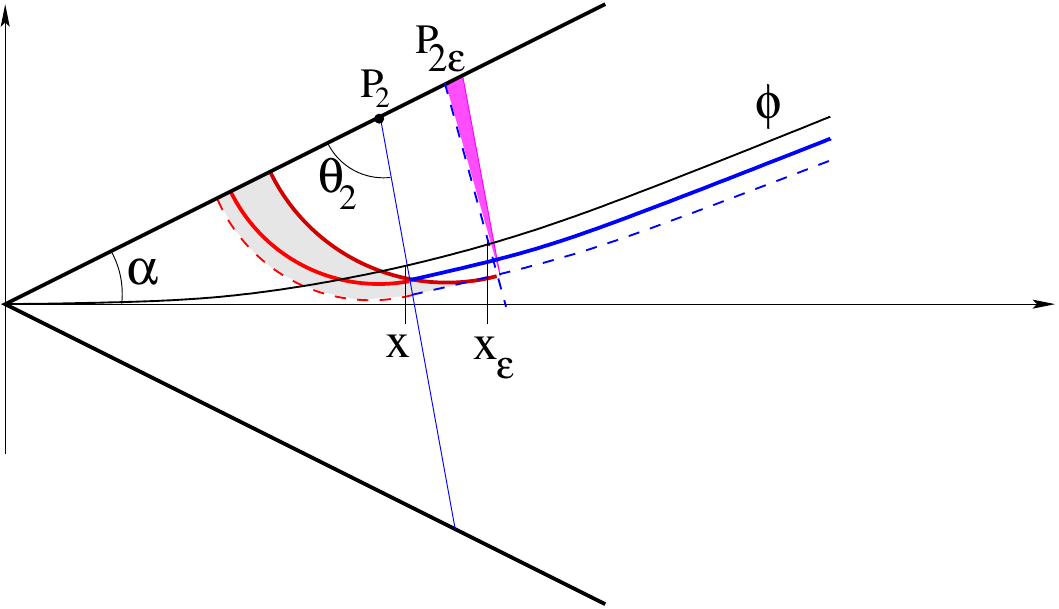}}}
\caption{\small  During the time interval $[t,\, t+\ve]$ the center of the active circumference moves from $P$ to $P_\ve$.  The shaded region has area $M\ve$.}
\label{f:tg31}
\end{figure}
Imposing the area identity (\ref{dareat}), we obtain
(see Fig.~\ref{f:tg31})
$${1\over 2} \theta_{2}  (r_2 +\ve)^2+ (x_\ve - x) \sqrt{ 1 + [\phi'(x)]^2}r_2  - (\theta_{2,\ve}-\theta_2){r_2^2\over 2} 
~\approx~
\theta_{2,\ve} {1\over 2} r_{2,\ve}^2 + M\ve.$$
Denoting the derivatives w.r.t.~time by an upper dot, letting $\ve\to 0$ we obtain
\bel{a3}\theta_2 r_2  + \dot x r_2 \sqrt{ 1 + [\phi'(x)]^2}  - \dot\theta_2 {r_2^2\over 2} 
~=~{1\over 2} \dot \theta_2 r_2^2 +\theta_2 r_2 \dot r_2 + M.\eeq
Observing that
$$\dot \theta_2~=~{d\over dx}\Big[ \arctan \phi'(x)\big] \dot x~=~{\phi''(x)\over 1+[\phi'(x)]^2}\, \dot x\,,$$
and dividing both sides of (\ref{a3}) by $r_2$, one finds
\bel{14}\theta_2  (\dot r_2 -1) +{M\over r_2}~=~\left( \sqrt{ 1 + [\phi'(x)]^2} - {\phi''(x)\over 1 + \bigl[ \phi'(x)\bigr]^2} 
\,r_2
\right)\, \dot x.
\eeq
Recalling (\ref{the12}) and (\ref{r12}), since $\tan\alpha=\cot\beta$  we obtain
\bel{dr21}\dot r_2 - 1~=~\left[ {\sin\beta \bigl(\cot\beta -\phi'(x)\bigr) \over
 \sin\theta_2(x)}
-    {\sin\beta \bigl(x\cot\beta -\phi(x)\bigr)\over \sin^2
 \theta_2(x)}\cos\bigl(\theta_2(x)\bigr){\phi''(x)\over 1 + \bigl[ \phi'(x)\bigr]^2}\right]\dot x\,. \eeq
Combined with (\ref{14}) and (\ref{the12}), this yields
\bel{15}\bega{l}\ds {M\over r_2} ~\ds=~\left( \sqrt{ 1 + [\phi'(x)]^2} - {\phi''(x)\over 1 + \bigl[ \phi'(x)\bigr]^2} 
\,r_2
\right)\dot x\\[4mm]
\quad \ds -  \bigl(\beta +\theta(x) \bigr) \left[ {\sin\beta \bigl(\cot\beta -\phi'(x)\bigr) \over \sin
 \bigl(\beta +\theta(x) \bigr)}
-    {\sin\beta \bigl(x\cot\beta -\phi(x)\bigr)\over \sin^2
 \bigl(\beta +\theta(x) \bigr)}\cos\bigl( \beta+\theta(x)\bigr)
{\phi''(x)\over 1 + \bigl[ \phi'(x)\bigr]^2}
\right]\dot x\,.\enda
\eeq
Solving for $\dot x = dx/dt_2$ we obtain
\bel{T2}\bega{l}\ds
{dt_2\over dx}~=~\T_2(x, t_2,\phi,\phi',\phi'')\\[4mm]
\ds \doteq~
{r_2\over M} \Bigg\{ \left( \sqrt{ 1 + [\phi'(x)]^2} - {\phi''(x)\over 1 + \bigl[ \phi'(x)\bigr]^2} 
\,r_2
\right) \\[4mm]
\ds\quad
-  \bigl(\beta +\theta(x) \bigr) \left[ {\cos\beta -\phi'(x)\sin\beta \over \sin
 \bigl(\beta +\theta(x) \bigr)}
-    { x\cos\beta -\phi(x)\sin\beta\over \sin^2
 \bigl(\beta +\theta(x) \bigr)}\cos\bigl(\beta+\theta(x)\bigr)
{\phi''(x)\over 1 + \bigl[ \phi'(x)\bigr]^2}
\right]\Bigg\}.
\enda
\eeq
An entirely similar computation can be done for $t<0$.  
To treat this case, it is simpler to reverse time and replace $\phi$ by $-\phi$.
$$\theta_1 r_1  - \dot x r_1 \sqrt{ 1 + [\phi'(x)]^2}+{1\over 2} \dot \theta_1 r_1^2
~=~-{1\over 2} \dot \theta_1 r_1^2 -\theta_1 r_1 \dot r_1 + M.$$
$$\theta_1 (\dot r_1 +1) -{M\over r_1}~=~\left( \sqrt{ 1 + [\phi'(x)]^2} + {\phi''(x)\over 1 + \bigl[ \phi'(x)\bigr]^2} 
\,r_1
\right)\, \dot x.$$
$$\dot r_1 + 1~=~\left[ {\sin\beta \bigl(\cot\beta +\phi'(x)\bigr) \over
 \sin\theta_1(x)}
-    {\sin\beta \bigl(x\cot\beta +\phi(x)\bigr)\over \sin^2
 \theta_1(x)}\cos\bigl(\theta_1(x)\bigr){-\phi''(x)\over 1 + \bigl[ \phi'(x)\bigr]^2}\right]\dot x\,,$$
$$\bega{l}\ds -{M\over r_1} ~\ds=~\left( \sqrt{ 1 + [\phi'(x)]^2} + {\phi''(x)\over 1 + \bigl[ \phi'(x)\bigr]^2} 
\, r_1
\right)\dot x\\[4mm]
\quad \ds -  \bigl(\beta -\theta(x) \bigr) \left[ {\sin\beta \bigl(\cot\beta +\phi'(x)\bigr) \over \sin
 \bigl(\beta -\theta(x) \bigr)}
+{\sin\beta \bigl(x\cot\beta +\phi(x)\bigr)\over \sin^2
 \bigl(\beta -\theta(x) \bigr)}\cos\bigl( \beta-\theta(x)\bigr)
{\phi''(x)\over 1 + \bigl[ \phi'(x)\bigr]^2}
\right]\dot x\,.\enda
$$
Eventually we obtain
\bel{T1}\bega{l}\ds
{dt_1\over dx}~=~\T_1(x, t_1,\phi,\phi',\phi'')\\[4mm]
\ds \doteq~-{r_1\over M} \Bigg\{ \left( \sqrt{ 1 + [\phi'(x)]^2} +{\phi''(x)\over 1 + \bigl[ \phi'(x)\bigr]^2} 
\,r_1
\right) \\[4mm]
\ds\quad 
-  \bigl(\beta -\theta(x) \bigr) \left[ {\cos\beta +\phi'(x)\sin\beta \over \sin
 \bigl(\beta -\theta(x) \bigr)}
+ { x\cos\beta +\phi(x)\sin\beta\over \sin^2
 \bigl(\beta -\theta(x) \bigr)}\cos\bigl( \beta-\theta(x)\bigr)
{\phi''(x)\over 1 + \bigl[ \phi'(x)\bigr]^2}
\right]\Bigg\}.
\enda
\eeq

\begin{remark}
{\rm
In the case where $\phi(x)\equiv 0$, the identity (\ref{15}) reduces to
$${M\over r_2}~=~\dot x -(\beta\cot \beta) \dot x~=~
[ 1-\beta\cot\beta] (\dot r_2-1)\tan\beta ~=~[\tan\beta-\beta](\dot r_2-1).$$
}
\end{remark}

\subsection{Necessary conditions.} 
Next, we impose the necessary condition (\ref{nc6}).
The curvature  of the maximal free interface $\gamma$ at the point 
$\gamma(x)= \bigl(x, \phi(x)\bigr)$ 
is computed by
\bel{curvature}
\omega(x)={1\over \rho(x)}~=~\frac{\phi''(x)}{\big(1+(\phi'(x))^2\big)^{3/2}}
\,.\eeq
\v
By (\ref{nc6}) it now follows that, for every $x>0$,
\bel{nc7} \int_0^{t_2(x)} {dt\over r_2(t)}-\int_{t_1(x)}^0 {dt\over r_1(t)}  ~=~\int_{t_1(x)}^{t_2(x)}
{dt\over \rho(x) + t}\,.\eeq
\v
We wish to use these identities to obtain a second order ODE
describing the curve $\gamma(\cdot)$.
Differentiating (\ref{nc7}) w.r.t.~$x$, one obtains
\bel{ode1}
{1\over r_2} {d t_2\over d x} +
{1\over r_1}{d t_1\over d x}~=~
{\omega \over 1+t_2\omega}\cdot {d t_2\over dx} - 
{\omega\over 1+ t_1\omega}\cdot {d t_1\over d x}.
\eeq

\begin{remark}{\rm
We choose the time $t=0$ as the time when the free interface is maximal, touching the
vertex $C$.  In this case, it is precisely the graph of $\phi(\cdot)$.

 We have
$${dt_1(x)\over dx} ~<~0,\qquad\qquad {dt_2(x)\over dx} ~>~0.$$

The radii $r_1, r_2$ are always taken to be positive.   However, for
$t>0$ the curvature is $\omega = {1\over r_2} $, while for $t<0$ it is $\omega = - {1\over r_1}$. 

There is a huge cancellation of the terms on the left hand sides of (\ref{nc7}) or equivalently (\ref{ode1}), 
while there is no cancellation on the right hand sides.}
\end{remark}

%
%
%
In view of (\ref{T2})-(\ref{T1}), multiplying both sides of (\ref{ode1}) by $M$ we obtain an equation of the form
\bel{F=G}
F(x, t_1, t_2, \phi, \phi',\phi'')~=~G(x, t_1, t_2, \phi, \phi',\phi'').\eeq
Setting
\bel{rtheta} \left\{
\bega{rl}r_1&\ds=~{x\cos\beta +\phi \sin\beta \over \sin
\bigl(\beta-\theta)} - t_1\,,
\\[4mm]
r_2&=~\ds{x\cos\beta -\phi\sin\beta \over \sin
\bigl(\beta+\theta\bigr)} + t_2\,,\enda\right.
\qquad\quad \theta\,=\,\arctan \phi', \qquad \omega~=~{\phi''\over \bigl(1+ |\phi'|^2\bigr)^{3/2}}\,,
\eeq
the left hand side of (\ref{F=G}) is computed by
\bel{Fdef} \bega{l}\ds
F(x, t_1, t_2, \phi, \phi',\phi'')~=~ \Bigg\{ \left( \sqrt{ 1 + |\phi'|^2} - {\phi''\over 1 +| \phi'|^2} \,
r_2
\right) \\[4mm]
\ds\qquad\qquad 
- (\beta +\theta) \left[ {\cos\beta -\phi'\sin\beta \over \sin
 (\beta +\theta)}
-    { x\cos\beta -\phi \sin\beta\over \sin^2 (\beta +\theta)}\cos(\beta+\theta)
{\phi''\over 1 + |\phi'|^2}
\right]\Bigg\}\\[4mm]
 \qquad \quad \ds -\Bigg\{ \left( \sqrt{ 1 + |\phi'|^2} +{\phi''\over 1 +|\phi'|^2} 
\, r_1
\right) \\[4mm]
\ds\qquad\qquad 
-  (\beta -\theta) \left[ {\cos\beta +\phi'\sin\beta \over \sin(\beta -\theta)}
+  { x\cos\beta +\phi\sin\beta\over \sin^2(\beta -\theta)}\cos( \beta-\theta)
{\phi''\over 1 + | \phi'(x)|^2}
\right]\Bigg\}.
\enda
\eeq
After some cancellations, one obtains
\bel{FF}\bega{l}\ds F(x, t_1, t_2, \phi, \phi',\phi'')~=~- {\phi''\over 1 +| \phi'|^2} \cdot
(r_1+r_2)
\\[4mm]
\ds\qquad\qquad 
-  {\beta +\theta\over  \sin(\beta +\theta) }
  \left[ (\cos\beta -\phi'\sin\beta)
-    { x\cos\beta -\phi \sin\beta\over \tan(\beta +\theta)}\cdot
{\phi''\over 1 +|\phi'|^2} 
\right]\\[6mm]
\ds\qquad\qquad 
+  {\beta -\theta \over \sin(\beta -\theta)}
 \left[ (\cos\beta +\phi'\sin\beta)+  { x\cos\beta +\phi \sin\beta\over \tan
 (\beta -\theta)}\cdot 
{\phi''\over 1 +|\phi'|^2} \right].
\enda
\eeq
On the other hand, recalling (\ref{T2})-(\ref{T1}), the right hand side is computed by
\bel{GG} \bega{l}\ds
G(x, t_1, t_2, \phi, \phi',\phi'')~=~ M\omega\left({dt_2\over dx} - {dt_1\over dx}\right) - 
M\omega\left[{t_2\omega \over 1+t_2\omega} {dt_2\over dx} -  {t_1\omega \over 1+t_1\omega} {dt_1\over dx}\right]\\[4mm]
\ds
 =~ M\omega\left[{t_1\omega \over 1+t_1\omega} \T_1 -  {t_2\omega \over 1+t_2\omega} \T_2\right]
+ {\phi''\over 1 +|\phi'|^2}\,(r_1+r_2)+{\phi''\over 1 +|\phi'|^2}\omega(r_1^2-r_2^2)
\\[4mm]
\ds\qquad- 
 {\beta +\theta\over\sin (\beta +\theta) } \left[ \cos\beta -\phi'\sin\beta
-    { x\cos\beta -\phi\sin\beta\over \tan (\beta +\theta)}\cdot 
{\phi''\over 1 +|\phi'|^2}
\right] \omega \,r_2 \\[4mm]
 \qquad \ds 
-{\beta -\theta\over \sin(\beta -\theta)}
 \left[\cos\beta +\phi'\sin\beta +  { x\cos\beta +\phi\sin\beta\over \tan (\beta -\theta)}\cdot {\phi''\over 1 +|\phi'|^2}
\right] \omega \,r_1\,.
\enda
\eeq

%
%

Summarizing, we need to solve a system of three ODEs for the three
functions $\phi(x), t_1(x), t_2(x)$.    Namely:
the two first order ODEs (\ref{T2})-(\ref{T1}) for the functions $t_2, t_1$, respectively,
together with the second order ODE (\ref{F=G}) for $\phi$.
These are
supplemented by the identities in  (\ref{rtheta}).
They are to be solved with boundary conditions
\bel{bc3}
t_1(0)~=~t_2(0)~=~\phi(0)~=~0.\eeq
Since the equations also involve the second order derivative $\phi''(x)$,
we expect to find a 1-parameter family of solutions.

The main difficulty in constructing solutions  stems from the fact that the above equations
(\ref{T2})-(\ref{T1}) and (\ref{F=G}) are implicit, and singular at the initial point $x=0$.
To gain some insight, we make the guess
\bel{guess} \bigl|\phi(x)\bigr|~<\!<~x\qquad\hbox{as} \quad x\to 0.\eeq
Starting with (\ref{T2})-(\ref{T1}) and neglecting higher order infinitesimals, we are led to
\bel{T22}
{dt_2\over dx}~\approx~{r_2\over M}\Bigg\{1-\beta\cot\beta -\phi'(x)\bigl(\cot\beta-\beta-\beta\cot^2\beta\bigr)-
x\phi''(x)\left(\cot\beta-\beta\cot^2\beta\right)\Bigg\},
\eeq
\bel{T11}{dt_1\over dx}~
\approx~-{r_1\over M}\Bigg\{1-\beta\cot\beta+\phi'(x)\bigl(\cot\beta-\beta-\beta\cot^2\beta\bigr)+
x\phi''(x)\left(\cot\beta-\beta\cot^2\beta\right)\Bigg\}.
\eeq
In view of (\ref{r12}),  this yields
\bel{tr12}
t_i(x)~=~\O(1) \cdot x^2, \qquad\qquad r_i(x)~=~x\cot \beta + o(x),\qquad i=1,2.\eeq
As usual, $o(x)$ denotes a higher order infinitesimal.
Inserting these expressions 
in (\ref{ode1}) 
we obtain
\bel{Fapp}
F(x,t_1,t_2,\phi, \phi',\phi'')~\approx~ \Bigg\{-2(\cot\beta-\beta-\beta\cot^2\beta)\phi'-2
(\cot\beta-\beta\cot^2\beta)x\phi'' \Bigg\}.
\eeq
On the other hand, the right hand side of (\ref{ode1}) can be approximated by
\bel{Gapp}
G(x,t_1,t_2,\phi,\phi',\phi'')~\approx~
{2} {\tan\beta-\beta\over\tan^2\beta}x\phi''
~=~{2} ( \cot\beta-\beta\cot^2\beta)\, x\phi''.
\eeq
To leading order, (\ref{Fapp}) and (\ref{Gapp}) yield the approximate ODE
\bel{ode6}
x \phi''(x)~=~\sigma \phi'(x),\qquad\qquad \phi(0)\,=\,0,\eeq
where
\bel{sigma}
\sigma~=~\sigma(\beta)~=~2{\cot\beta -\beta-\beta\cot^2 \beta\over
4\beta\cot^2\beta-4\cot\beta}~=~{\beta+\beta\tan^2\beta-\tan\beta\over
2\tan\beta- 2\beta}~=~{\beta-\sin\beta\cos\beta\over 2\cos\beta(\sin\beta-\beta\cos\beta)}\,.
\eeq
We claim that $\sigma(\beta)>1$ for every $\beta\in\,]0,\,\pi/2]$. 
Indeed, a Taylor approximation yields
\bel{siglim}\lim_{\beta\to 0+ }\sigma(\beta) =~~\lim_{\beta\to 0+ } {\beta-\big(\beta-\beta^3/6\big)\big(1-\beta^2/2\big)\over
2\big(\beta-\beta^3/6)-\beta\big(1-\beta^2/2\big)}~=~1.\eeq
Moreover, by elementary differentiations one obtains
$${d\over d\beta} \Big(\beta\tan^2\beta+\beta -\tan \beta\Big)~=~
{d\over d\beta} \left({\beta\over\cos^2\beta} -{\sin \beta\over\cos\beta}\right)
=~\beta\,{2\sin\beta\over\cos^3\beta}\,,$$
hence 
\bel{sigder}\bega{rl}\ds
{d\sigma(\beta)\over d\beta}&= ~\ds
{\beta{2\sin\beta\over\cos^3\beta}(\tan\beta- \beta)-(\beta+\beta\tan^2\beta-\tan\beta)\tan^2\beta\over 2(\tan\beta- \beta)^2}\\[4mm]
&=~\ds {\sin\beta\over2(\tan\beta- \beta)^2\cos^4\beta}
\big(2\beta(\sin\beta- \beta\cos\beta)-\sin\beta(\beta-\sin\beta\cos\beta)\big)\\
&>~\ds {\sin^2\beta\over2(\tan\beta- \beta)^2\cos^4\beta}
\big(2(\sin\beta- \beta\cos\beta)-(\beta-\sin\beta\cos\beta)\big)>0,
\enda
\eeq
in view of the fact that
$$ {d\over d\beta} \big(2(\sin\beta- \beta\cos\beta)-(\beta-\sin\beta\cos\beta)\big)~=~
2\beta\sin\beta-2\sin^2\beta~>~0\,.$$
Together, (\ref{siglim}) and (\ref{sigder}) imply that $\sigma(\beta)>1$.
%

The general solution to the linear equation (\ref{ode6}) is
\bel{sol6}
\phi(x)~=~c x^{\sigma+1},
\eeq
where $c$ is an arbitrary constant.
Since $\sigma>1$, this 
yields a 1-parameter family of solutions, for which the asymptotic condition (\ref{guess}) holds.

\section{Local solutions to the equations of maximal free interfaces}
\label{s:16}
\setcounter{equation}{0}
Aim of this section is to prove a local existence theorem for solutions to the singular, implicit Cauchy problem
\bel{new} \left\{ \bega{l} H(x,t_1, t_2, \phi,\phi',\phi'')~=~0,\\[2mm]
 t_1'~=~\T_1(x, t_1, \phi,\phi',\phi''),\\[2mm]
  t_2'~=~\T_2(x, t_2, \phi,\phi',\phi''),\enda
\right.\qquad\qquad \left\{ \bega{rl} \phi(0)&=~0,\cr
\ds\lim_{x\to 0+} 
{\phi(x)\over x^{\sigma+1}}&=~c,\cr
t_1(0)&=~0,\cr
t_2(0)&=~0.\enda\right.
 \eeq
Here $\T_2,\T_1$ are the functions introduced at (\ref{T2})-(\ref{T1}), while
$H=G-F$, with $F,G$ as in (\ref{FF})-(\ref{GG}).

\begin{theorem} \label{t:161}
Let $0<\beta<\pi/2$ and  $c\in \R$ be given, and let  $\sigma>0$ be the constant in (\ref{sigma}).
Then there exists $x^\dagger>0$ and a local solution $x\mapsto \bigl(\phi(x),t_1(x), t_2(x)\bigr)$
to
the initial value problem (\ref{new}), defined for $x\in [0,x^\dagger]$.
\end{theorem}

{\bf Proof.} {\bf 1.}
In view of the previous analysis, to leading order the equation (\ref{F=G}) reduces to (\ref{ode6}).
We thus expect that
solutions will satisfy
$$\phi(x)\,=\,\O(1)\cdot x^{\sigma+1}, \quad \phi'(x)\,=\,\O(1)\cdot x^{\sigma},\quad \phi''(x)\,=\,\O(1)\cdot x^{\sigma-1},\qquad t_1(x), \,t_2(x)\,=\, \O(1)\cdot x^2.$$
We thus consider a domain of the form
\bel{Dp1} \D~=~\Big\{ (t_1, t_2,
\phi,\phi')\,;~~\bigl| t_i(x)\bigr| \leq  C_0\, x^2,
\quad \bigl|\phi(x)\bigr|\leq  C_0 \,x^{\sigma+1}  ,~~ \bigl|\phi'(x)\bigr|\leq  C_0 \,x^{\sigma} 
 \Big\},\eeq
for a suitable constant $C_0$.

Within this domain,  the local solution will be obtained 
as the fixed point of
a Picard-type operator
\bel{Pic}
\Big(\P(t_1,t_2, \phi,\phi')\Big)(x)~=~(\Tilde t_1,\Tilde t_2,\Tilde\phi,\Tilde\phi')(x),
\eeq
where the functions on the right hand side are defined as follows.
We begin by constructing the function
$\Ups=\Ups(x, t_1, t_2,\phi,\phi')$, implicitly defined by the identity
\bel{Ups}H(x,t_1,t_2, \phi,\phi',\Upsilon)~=~0.\eeq
Details of the construction of $\Ups$ will be worked out in step 2.
We then define
\bel{Ti} T_i(x,t_1,t_2, \phi,\phi')~=~\T_i\bigl(x,t_i, \phi,\phi',\Ups(x, t_1, t_2,\phi,\phi')\bigr),\eeq
and set
\bel{Pict}\Tilde t_i(x)~\doteq~\int_0^x T_i(x,t_1,t_2,\phi,\phi')\, dy,\qquad\quad i=1,2.
\eeq
Next, we write the equation  $H=0$  
in the equivalent form
\bel{odeK}
x \phi''(x ) - \sigma \phi'(x)~=~K(x,t_1,t_2,\phi,\phi',\phi'').\eeq
for a suitable function $K$.  Solving (\ref{odeK}) with $\phi''=\Ups$, in view of the second boundary condition
at (\ref{new})
this leads to 
\bel{Pic1}\Tilde\phi'(x)~\doteq~x^\sigma \left((\sigma+1)c+\int_0^x y^{-\sigma-1}\, K(y,t_1,t_2,\phi,\phi',\Ups)\, dy
\right)\eeq
 \bel{Pic0}\Tilde \phi(x)~\doteq~\int_0^x \Tilde \phi'(y) dy.  \eeq
 In the remainder of the proof, we will show that the above formulas (\ref{Pict}), (\ref{Pic1}), (\ref{Pic0})
yield a strictly contractive  map, on an interval $[0, x^\dagger]$ small enough.
\v
{\bf 2.} In this step, using the implicit function theorem we construct the function $\Ups$
and provide estimates on the function $K$ in (\ref{odeK}).
Computing the partial derivatives of $H=G-F$, in view of  (\ref{FF})-(\ref{GG}) and (\ref{rtheta})  we obtain
\bel{pHp''}{\partial \over\partial \phi''} H~=~
4x(\cot \beta -  \beta\cot^2\beta)  +
\O(1)\cdot x^{1+\ve},\eeq
\bel{pHp'} {\partial \over\partial \phi'}H~=~ -2(\beta+\beta\cot^2\beta-\cot\beta) + \O(1)\cdot x^\ve,
\eeq
\bel{pHp} {\partial \over\partial \phi} H~=~\O(1)\cdot |\phi''|+
\O(1)\cdot 
\bigl(x+|\phi|\bigr)|\phi''|^2~=~\O(1)\cdot x^{\sigma-1}.
\eeq
\bel{pHt} {\partial \over\partial t_i} H~=~\O(1)\cdot |\phi''|+
\O(1)\cdot 
\bigl(x+|\phi|\bigr)|\phi''|^2~=~\O(1)\cdot x^{\sigma-1}.\eeq
Since $\partial H/\partial \phi'' >0$, this yields the local existence of the function $\Ups$ implicitly defined by 
(\ref{Ups}). 
 Here, $\ve=\min\{1,\sigma\}$.

 Moreover, we have
\bel{pU1}{\partial \Ups\over\partial \phi'}~=~ -{ \partial H/\partial \phi'\over
\partial H/\partial \phi''}~=~{2(\beta+\beta\cot^2\beta-\cot\beta) + \O(1)\cdot x^\ve\over
4x(\cot \beta -  \beta\cot^2\beta) + \O(1)\cdot x^{1+\ve}}~=~\sigma +\O(1)\cdot x^{\ve-1}.
\eeq
Similarly, we obtain
\bel{pU2}{\partial \Ups\over\partial \phi}~=~ -{ \partial H/\partial \phi\over
\partial H/\partial \phi''}~=~\O(1) \cdot x^{\sigma-2},\qquad\qquad 
{\partial \Ups\over\partial t_i}~=~ -{ \partial H/\partial t_i\over
\partial H/\partial \phi''}~=~\O(1) \cdot x^{\sigma-2}.\eeq
Next, we observe that the equation $H=0$ is equivalent to 
$$\phi''(x)~=~\Ups\bigl(x, t_1, t_2,\phi(x),\phi'(x)\bigr).$$
In turn, this yields 
$$x\phi''~=~x\Ups,\qquad\qquad x\phi'' - \sigma \phi'~=~K,$$
where
\bel{K}K(x, t_1, t_2, \phi,\phi')~=~x\,\Ups(x, t_1, t_2, \phi,\phi')-\sigma\phi'.
\eeq
In view of (\ref{pU1})-(\ref{pU2}), the partial derivatives of $K$ thus satisfy the bounds
\bel{pK}{\partial K\over\partial \phi'}~=~ O(1)\cdot x^{\ve},\qquad\qquad
{\partial K\over\partial \phi}~=~\O(1) \cdot x^{\sigma-1},\qquad\qquad 
{\partial K\over\partial t_i}~=~\O(1) \cdot x^{\sigma-1}.\eeq
\v
{\bf 3.} We now provide estimates on the functions $T_1, T_2$ introduced at (\ref{Ti}).
Recalling (\ref{T2})-(\ref{T1}) and (\ref{rtheta}), we begin by  estimating the partial 
derivatives
\bel{pT4}
{\partial \T_i\over\partial \phi}~=~\O(1),\qquad\qquad{\partial \T_i\over\partial t_i}~=~\O(1),
\eeq
\bel{pT5}
{\partial \T_i\over\partial \phi'}~=~\O(1)\cdot x,
\qquad\qquad 
{\partial \T_i\over\partial \phi''}~=~\O(1)\cdot x^2.
\eeq

By (\ref{T2})-(\ref{T1}) and (\ref{pU1})-(\ref{pU2}) it follows
\bel{dT6}
{\partial T_i\over\partial \phi}~=~{\partial \T_i\over\partial \phi} + {\partial \T_i\over \partial \phi''}\,{\partial \Ups\over\partial \phi}~=~\O(1)
\eeq
\bel{dT7}
{\partial T_i\over\partial \phi'}~=~{\partial \T_i\over\partial \phi'} + {\partial \T_i\over \partial \phi''}\,{\partial \Ups\over\partial \phi'}~
=~\O(1)\cdot ( x + x^2 x^{\ve -1})~=~\O(1)\cdot x,
\eeq
\bel{dT8}
{\partial T_i\over\partial t_j}~=~{\partial \T_i\over\partial t_j} + {\partial \T_i\over \partial \phi''}\,{\partial \Ups\over\partial t_j}~=~\O(1) \cdot (1 + x^2\, x^{\sigma-2})~=~\O(1).
\eeq

\v{\bf 4.} We need to prove that the transformation 
$\Big(\P(t_1,t_2, \phi,\phi')\Big)(x)=(\Tilde t_1,\Tilde t_2,\Tilde\phi,\Tilde\phi')(x)$, as in (\ref{Pic}), 
maps the domain (\ref{Dp1}) into itself.
Let $(t_1, t_2, \phi, \phi')\in\D$. Then, recalling (\ref{pK}), we obtain
\begin{align}\label{D1}
\ds \Tilde\phi'(x)&=x^\sigma(\sigma+1)c+x^\sigma \int_0^x y^{-\sigma-1}
\Bigg\{ \int_0^1{\partial K\bigl(\lambda(t_1,t_2, \phi,\phi')(y)\bigr) \over \partial \phi}
\,d\lambda\cdot \phi(y)\nonumber\\[4mm]
&\quad \ds +\int_0^1{\partial K\bigl(\lambda(t_1,t_2, \phi,\phi')(y)\bigr) \over \partial \phi'}
\,d\lambda\cdot \phi'(y) +\sum_{i=1,2} \int_0^1{\partial K\bigl(\lambda(t_1,t_2, \phi,\phi')(y)\bigr) \over \partial t_i}
\,d\lambda\cdot t_{i}(y) \Bigg\}dy\nonumber\\[4mm]
&=x^\sigma(\sigma+1)c+x^\sigma O(1)\big\{x^\ve+x^{\sigma}+x\big\}=(\sigma+1)cx^\sigma+O(1)x^{\sigma+\ve}.
\end{align}
Hence, 
$$\Tilde\phi(x)=\int_0^x\Tilde \phi'(y)dy=x^{\sigma+1}c+O(1)x^{\sigma+1+\ve}\,.$$
For $i=1,2$, working as in (\ref{D1}) and using  (\ref{dT6})--(\ref{dT8}), one gets
\begin{align*}
\ds \Tilde t_i(x)&=\int_0^x \Bigg\{ \int_0^1{\partial T_i\bigl(\lambda(t_1,t_2, \phi,\phi')(y)\bigr) \over \partial \phi}
\,d\lambda\cdot \phi(y)\\[4mm]
&\quad \ds +\int_0^1{\partial T_i\bigl(\lambda(t_1,t_2, \phi,\phi')(y)\bigr) \over \partial \phi'}
\,d\lambda\cdot \phi'(y) +\sum_{j=1,2} \int_0^1{\partial T_i\bigl(\lambda(t_1,t_2, \phi,\phi')(y)\bigr) \over \partial t_j}
\,d\lambda\cdot t_{j}(y) \Bigg\}dy\\[4mm]
&=O(1)\big\{x^{\sigma+2}+x^3\big\}=O(1)x^{2+\ve}.
\end{align*}
Taking  $C_0\geq2(\sigma+1)c$ and $x^{\dagger}$ small enough, we obtain that, for any $x\in[0,x^{\dagger}]$, 
$\Big(\P(t_1,t_2, \phi,\phi')\Big)(x)\in\D$. 

\v{\bf 5.}
To achieve the contraction property, consider the norm
 \bel{norm}
\bigl\| (t_1, t_2, \phi,\phi')\bigr\|~=~\sup_{0<x<x^\dagger} \max\left\{ {\bigl|\phi(x)\bigr|\over x^{\sigma+1}}\,,~
{\bigl|\phi'(x)\bigr|\over x^\sigma}\,,
~{\bigl| t_1(x)\bigr|\over x^2}\,,~~{\bigl| t_2(x)\bigr|\over x^2}\right\}\,.
\eeq
Given two sets of functions 
$(t_{j1}, t_{j2}, \phi_j,\phi'_j)$,  $j=1,2$, such that
\bel{d12}
\Big\| (t_{11}, t_{12}, \phi_1,\phi'_1)-(t_{21}, t_{22}, \phi_2,\phi'_2)\Big\|~=~\delta~>~0,\eeq
we will prove that the corresponding Picard iterates satisfy
\bel{dP12}\Big\| (\Tilde t_{11}, \Tilde t_{12}, \Tilde \phi_1,\Tilde \phi'_1)-
(\Tilde t_{21}, \Tilde t_{22}, \Tilde \phi_2,\Tilde \phi'_2)\Big\|~\leq~{\delta\over 2}\,.\eeq
\v
{\bf 6.} We are now ready to prove the contraction property (\ref{dP12}), with distance induced by the norm
(\ref{norm}).  
If (\ref{d12}) holds, then
\bel{dista} \bigl|\phi_1(x)-\phi_2(x)\bigr|\,\leq\, \delta x^{\sigma+1},\qquad\bigl|\phi'_1(x)-\phi'_2(x)\bigr|\,\leq\, \delta 
x^\sigma,
\qquad
\bigl|t_{1i}(x) - t_{2i}(x)\bigr|~\leq~\delta x^{2} \,,\quad i=1,2.\eeq
Consider the intermediate points
$$Z(\lambda,x)~\doteq~\lambda (t_{11}, t_{12}, \phi_1,\phi'_1)(x) + (1-\lambda)(t_{21}, t_{22}, \phi_2,\phi'_2)(x).$$
By (\ref{Pic1}), the difference between the two values of $\Tilde \phi'$ is  computed by
\bel{Tppp}\bega{l}
\ds \Tilde\phi'_1(x)-\Tilde\phi'_2(x)~=~x^\sigma \int_0^x y^{-\sigma-1}\Bigg\{ \int_0^1{\partial K\bigl(Z(\lambda,y)\bigr) \over \partial \phi}
\,d\lambda\cdot (\phi_1(y)-\phi_2(y))\\[4mm]
\qquad \ds +\int_0^1{\partial K\bigl(Z(\lambda,y)\bigr) \over \partial \phi'}
\,d\lambda\cdot (\phi'_1(y)-\phi'_2(y)) +\sum_{i=1,2} \int_0^1{\partial K\bigl(Z(\lambda,y)\bigr) \over \partial t_i}
\,d\lambda\cdot (t_{i1}(y)-t_{i2}(y)) \Bigg\}dy.
\enda
\eeq
Recalling (\ref{pK}) and (\ref{dista}), one obtains
\bel{Te4}\bega{l}
\ds \Tilde\phi'_1(x)-\Tilde\phi'_2(x)~=~ \O(1)\cdot
\Big\{x^{\sigma-1} \,\delta x^{\sigma+1} + x^\ve \, \delta x^\sigma + x^{\sigma-1} \delta x^2\Big\}
~=~\O(1)\cdot \delta x^{\sigma+\ve}.
\enda\eeq
In turn, this yields
\bel{Te5}
 \Tilde\phi_1(x)-\Tilde\phi_2(x)~=~ \int_0^x \bigl[ \Tilde\phi'_1(y)-\Tilde\phi'_2(y)  \bigr]\, dy
~=~\O(1)\cdot \int_0^x \delta y^{\sigma+\ve}\, dy~=~\O(1) \cdot \delta x^{\sigma+1+\ve}.\eeq
Finally we estimate the change in the $t_i$ variables. 
For $i=1,2$, the same approach as in (\ref{Tppp}) now yields
$$\bega{l}
\ds \Tilde t_{1i}(x)-\Tilde t_{2i}(x)~=~\int_0^x \Bigg\{ \int_0^1{\partial T_i\bigl(Z(\lambda,y)\bigr) \over \partial \phi}
\,d\lambda\cdot (\phi_1(y)-\phi_2(y))\\[4mm]
\qquad \ds +\int_0^1{\partial T_i\bigl(Z(\lambda,y)\bigr) \over \partial \phi'}
\,d\lambda\cdot (\phi'_1(y)-\phi'_2(y)) +\sum_{j=1,2} \int_0^1{\partial T_i\bigl(Z(\lambda,y)\bigr) \over \partial t_j}
\,d\lambda\cdot (t_{j1}(y)-t_{j2}(y)) \Bigg\}dy.
\enda
$$
Using (\ref{dT6})--(\ref{dT8}) and recalling (\ref{dista}), we now obtain 
\bel{Te6}
\ds \Tilde t_{1i}(x)-\Tilde t_{2i}(x)~=~\int_0^x \O(1)\cdot \bigl[1\cdot \delta y^{\sigma+1} + y\cdot \delta y^\sigma
+ 1\cdot y^2\bigr]\, dy~=~\O(1)\cdot x^{2+\ve}.
\eeq

As $x$ ranges in an interval  $ [0, x^\dagger]$ small enough, the three inequalities (\ref{Te4})--(\ref{Te6}) 
yield the contractive property (\ref{dP12}).
\v
{\bf 6.} 
The unique fixed point of the contractive transformation defined at (\ref{Pict}), (\ref{Pic1}) and (\ref{Pic0}) 
now provides the desired solution to the Cauchy problem (\ref{new}).

\endproof

\section{Concluding remarks}
\label{s:17}
\setcounter{equation}{0}
In this paper we analyzed a family of set-motion problems, discussing existence of optimal solutions and necessary
conditions toward their explicit computation.
We like to think of these as ``time dependent isoperimetric problems". Indeed, by the area formula (\ref{dareat}),
to minimize the area one should reduce the perimeter as fast as possible.  In the optimal motions
considered in Theorem~\ref{t:41}, at each time $t\in [0,T]$ the set $\Omega(t)$ has the shortest 
relative boundary $\partial \Omega(t)\cap V$, 
compared with all other subsets with the same area.

Among several problems which are left open, a key issue is the regularity of optimal solutions.  While the existence
results proved in Section~\ref{s:5} yield optimal solutions in the class of BV functions, it seems reasonable to conjecture that for a.e.~$t\in [0,T]$ the set $\Omega(t)$ should 
satisfy an interior ball condition.   That is: $\Omega(t)$ should be the union of discs of radius $r(t)>0$.
Such regularity would likely suffice for deriving necessary conditions for optimality.

A second goal would be to use the analysis in Sections~\ref{s:13}--\ref{s:16} in order to construct
set motions that satisfy all optimality conditions {\em globally in time}.
Two cases appear to be particularly interesting, where either (i) $V$ is a polygon, or (ii) $V$ is a convex set with smooth boundary, whose curvature has finitely many local maxima and minima.

Both of these problems are left for future investigation.

\v
{\bf Aknowledgements.}  
The second author is a member of the Gruppo Nazionale per l’Analisi Matematica, la Probabilit\`a e le loro Applicazioni (GNAMPA), of the Istituto Nazionale di Alta Matematica (INdAM), and is supported by the GNAMPA project CUP\_E53C23001670001.

The research of the third author was partially supported by Portuguese funds through the Center for Research and Development in Mathematics and Applications (CIDMA), within project UID/MAT/04106/2019.

\end{document}